\documentclass[11pt]{amsart}
\usepackage{amsmath,amssymb,amsthm,tikz,mathrsfs,upgreek,dsfont,bm}

\RequirePackage{color}
\RequirePackage[colorlinks, urlcolor=my-blue,linkcolor=my-red,citecolor=my-green]{hyperref}
\definecolor{my-blue}{rgb}{0.0,0.0,0.6}
\definecolor{my-red}{rgb}{0.5,0.0,0.0}
\definecolor{my-green}{rgb}{0.0,0.5,0.0}
\definecolor{nicos-red}{rgb}{0.75,0.0,0.0}
\definecolor{light-gray}{gray}{0.6}
\definecolor{really-light-gray}{gray}{0.8}
\definecolor{sussexg}{rgb}{0.0,0.5,0.5}
\definecolor{sussexp}{rgb}{0.5,0.0,0.5}

\usepackage{nicefrac}
\newtheorem{theorem}{\textcolor{sussexp!90}{\bf THEOREM}}[section]
\newtheorem{lemma}[theorem]{\textcolor{sussexp}{\bf Lemma}}
\newtheorem{proposition}[theorem]{\textcolor{sussexp}{\bf Proposition}}
\newtheorem{corollary}[theorem]{\textcolor{sussexp}{\bf Corollary}}
\newtheorem{assumption}[theorem]{\bf Assumption}
\newtheorem{definition}[theorem]{\textcolor{purple}{\bf Definition}}
\numberwithin{equation}{section}
\theoremstyle{remark}
\newtheorem{remark}[theorem]{\textcolor{nicos-red}{\bf Remark}}

\newcommand{\be}{\begin{equation}}
\newcommand{\ee}{\end{equation}}

\newcommand{\fl}[1]{\lfloor{#1}\rfloor} 
\newcommand{\ce}[1]{\lceil{#1}\rceil}

\def\bN{\mathbb{N}}
\def\bP{\mathbb{P}}

\def\bR{\mathbb{R}}
\def\bZ{\mathbb{Z}}

\def\e{\varepsilon}
\def\Leb{{\textrm{Leb}}}

 \def\Z{\bZ} 
 
\def\R{\bR}
\def\N{\bN}

\def\P{\bP} %% environment measure 

\definecolor{darkgreen}{rgb}{0.0,0.5,0.0}
\definecolor{darkblue}{rgb}{0.0,0.0,0.3}
\definecolor{nicosred}{rgb}{0.65,0.1,0.1}
\definecolor{light-gray}{gray}{0.7}
\allowdisplaybreaks[1]

\RequirePackage{datetime} %comment out when arxiving or submitting

\begin{document}
\usdate
\title[Space-time discontinuous TASEP]
{Hydrodynamic limits for TASEP with space-time discontinuities}

\author{Jacob Butt}
\address{Jacob Butt, University of Rome `La Sapienza', Rome, Italy}
\email{jacob.butt@uniroma1.it}
\author{Nicos Georgiou}
\address{Nicos Georgiou, University of Sussex, Brighton, UK.}
\email{N.Georgiou@sussex.ac.uk}
\author{Enrico Scalas}
\address{Enrico Scalas, University of Rome `La Sapienza', Rome, Italy}
\email{enrico.scalas@uniroma1.it}

\date{\today}

\begin{abstract}
We develop a hydrodynamic theory for a height-dependent version of the totally
asymmetric simple exclusion process in which the jump rate at a growth site is sampled from a
macroscopic two-dimensional speed function evaluated at the spatial coordinate and the current
height level. The speed function is allowed to have discontinuities along locally finitely many curves.
Through the TASEP height-function representation, the process is coupled to an inhomogeneous
directed last-passage percolation model whose exponential rates vary discontinuously in the two
macroscopic LPP coordinates. Combining the law of large numbers for this last-passage model
with an extension of the variational coupling method, we prove a hydrodynamic limit for the
height function and for the associated particle density.

The limiting current is characterised by a Lax--Oleinik type variational formula built from the
discontinuous last-passage shape function. We then identify the first-order PDE structure selected
by the microscopic dynamics. At points of differentiability of the limiting current and continuity
of the sampled coefficient, the current solves a Hamilton--Jacobi equation whose Hamiltonian
depends discontinuously on the spatial variable and on the value of the solution itself. At discontinuities, the variational formula leads to a natural envelope-based discontinuous
viscosity formulation, and we prove that the limiting current satisfies this formulation. Finally, when the coefficient has only spatial discontinuities, we prove uniqueness
of the Hamilton--Jacobi solution in the natural class of nondecreasing Lipschitz currents, and
identify its spatial derivative as the maximal-current weak solution of the associated scalar
conservation law with discontinuous flux.
\end{abstract}

\maketitle

\setcounter{tocdepth}{1}
\tableofcontents

\maketitle

\section{Introduction}

Interacting particle systems and random growth models provide one of the most successful meeting points between probability, mathematical physics, variational analysis, and nonlinear PDE. A basic example is the totally asymmetric simple exclusion process (TASEP), in which particles on the one-dimensional integer lattice jump only one step to the right, subject to the target position being empty. The same dynamics can be viewed through its integrated current, or height function, and this height function can in turn be coupled to a directed last-passage percolation (LPP) model. Thus, particle densities, growing interfaces, maximizing lattice paths, and Hamilton--Jacobi variational formulas are different representations of the same underlying object.

In the homogeneous exponential setting, this correspondence is especially fruitful.  
Consider the lattice corner $\Z^2_+$. The environment $\tau = \{ \tau_{i,j} \}_{(i,j) \in \Z^2_+}$ is a collection of independent and identically distributed (i.i.d.) exponential random variables of rate 1. In the classical model, the last passage time from ${\bf u}$ to ${\bf v}$, ${\bf u, v} \in \Z^2_+$, is defined simply as 
\[
G_{\bf u, \bf v}=\max_{\pi\in\Pi_{\bf u, \bf v}}\bigg\{\sum_{(i,j)\in\pi}\tau_{i,j}\bigg\},
\]
where the maximum is over ordered paths $\pi$ of vertices that connect $\bf u$ to $\bf v$ and consecutive vertices differ only by a standard unit vector:
\[ \pi_{\bf u, \bf v} = \Big\{(i_0, j_0), (i_1, j_1), \ldots,  (i_m, j_m) : (i_{k+1}, j_{k+1}) - (i_{k}, j_{k}) \in \{ {\bf e}_1, {\bf e}_2 \}, {\bf u} =(i_0, j_0), {\bf v} =(i_m, j_m)\Big\}.\]
The set of all such admissible paths is denoted by 
$\Pi_{ \bf u, \bf v} = \{ \pi_{\bf u, \bf v}\}.$

Since the environment is i.i.d.~ with finite moments, the subadditive ergodic theorem \cite{Liggett1985} gives a law of large numbers (LLN): Set $\bf u = 0$ and for $(x,y) \in \bR^2_+$, set ${\bf v}_n = (\ce{nx}, \ce{ny})$. Then the law of large numbers is the existence of an almost surely (a.s.) continuous, concave and 1-homogeneous {\em limiting shape function} $\gamma$ on $\bR_+^2$ 
\be\label{thm:classicalLPP}
\lim_{n\to \infty} \frac{G_{{\bf 0}, {\bf v}_n}}{n} =: \gamma(x,y).
\ee   
The a.s.\ result is under the i.i.d.~measure $\P$ of the exponential environment. 

The shape function $\gamma(x,y)$ was identified to be 
\be \label{eq:limshap}
\gamma(x,y) = (\sqrt{x} + \sqrt{y})^2
\ee
via the duality between the corner growth model and the height function for the TASEP in \cite{rost1981}, 
while Johansson connected the shape
fluctuations to random-matrix laws \cite{Johansson2000Shape}.  

The bridge between the last-passage and particle-system pictures is provided by the TASEP height function. In this setting \cite{rost1981} studied the hydrodynamic limit. That is, the limiting (macroscopic) height viewed as a function and the PDE it solves. The microscopic exclusion rule becomes a nonlinear conservation law after Euler scaling.  If $\rho_0$ is the initial macroscopic particle density then, in the homogeneous rate-one process, the limiting density $\rho(x,t)$ is governed by the inviscid Burgers equation
\[
        \rho_t+\bigl(\rho(1-\rho)\bigr)_x=0,
        \qquad
        \rho(x,0)=\rho_0(x),
\]
while the integrated current, or height function, is described by the corresponding Hamilton--Jacobi equation
\[
        v_t+v_x(1-v_x)=0,
        \qquad
        v(x,0)=v_0(x),
        \qquad
        v_x=\rho .
\]  
Rost \cite{rost1981} proved the rarefaction-fan density profile and local
Bernoulli equilibrium for the totally asymmetric exclusion process, and also
gave a geometric growth interpretation on the lattice in which the rescaled occupied region
has boundary
\[
\sqrt{s_1}+\sqrt{s_2}=1.
\]
In modern last-passage terminology, this is equivalent to the exponential
corner-growth shape \eqref{eq:limshap},
although Rost's formulation is in terms of the exclusion process and its
associated growth model rather than the last-passage notation used above. For a detailed account of the exponential corner-growth model, including its
stationary structure, variational formulas, Busemann functions, and fluctuation
exponents, see \cite{Seppalainen2018CornerGrowth} for a survey.

 Exponential LPP and TASEP sit at the centre of the Kardar--Parisi--Zhang universality class,
introduced in \cite{KardarParisiZhang1986DynamicScaling}; see also
\cite{Corwin2012KPZ,KriecherbauerKrug2010Pedestrian} for surveys.  They provide
both exactly solvable examples and a route to continuum objects such as the KPZ fixed point and the
directed landscape \cite{DauvergneOrtmannVirag2022DirectedLandscape,MatetskiQuastelRemenik2021KPZFixedPoint}.  
The solvable-particle-system viewpoint also extends beyond continuous-time
TASEP.  In step-initial TASEP, the connection with corner growth and
exponential last-passage percolation gives a classical route from particle
systems to deterministic limit shapes and Tracy--Widom fluctuations
\cite{Johansson2000Shape}.  

Moreover, exact solvability is not limited to total
asymmetry: Bethe-ansatz formulas for ASEP and weak-asymmetry limits link
exclusion processes to the KPZ equation
\cite{BertiniGiacomin1997KPZ,TracyWidom2008ASEP}.  
At the continuum end of this picture, the solution theory for the KPZ equation of \cite{Hairer2013}
made precise a canonical singular SPDE limit for one-dimensional random growth models.

A parallel continuous-space picture is provided by
Hammersley's interacting particle process and the longest increasing
subsequence problem, where hydrodynamic limits connect microscopic particle
dynamics to Burgers-type conservation laws and Ulam-type asymptotics
\cite{AldousDiaconis1995Hammersley,Seppalainen1996BurgersLIS,
BaikDeiftJohansson1999LIS}. A related variational approach appears in the study of randomly forced Burgers equations, where global solutions, invariant measures, and stationary regimes are
constructed through Langrangian or action-minimizing curves \cite{EKhaninMazelSinai2000, Bakhtin2013, BakhtinCatorKhanin2014}. 

More recent algebraic
frameworks, including Macdonald processes, $q$-TASEP, and stochastic vertex
models, further show that particle systems, growth models, polymer models, and
vertex models often share common integrable structures
\cite{BorodinCorwin2014Macdonald,BorodinCorwinSasamoto2014Duality,
BorodinCorwinGorin2016SixVertex}.

The preceding examples all exploit a homogeneous or algebraically controlled
microscopic environment.  From the hydrodynamic point of view, however, it is
natural to ask what remains of this picture when the jump rates are no longer
constant, and away from the solvable models. Inhomogeneous rates arise both as models of defects, media with
spatially varying transport capacity, and growth environments whose local
speed depends on the macroscopic position.  They also test the robustness of
the variational correspondence: once the weights are no longer identically
distributed, the limiting shape is no longer expected to be given by a single
explicit homogeneous formula, but the same `last-passage to height-function
mechanism' may still identify the correct macroscopic variational problem.

\subsection{Inhomogeneous LPP and particle systems}
The simplest form of inhomogeneity is spatial.  Already in this setting, the
variational viewpoint is useful because a microscopic defect or spatially
varying rate need not pass directly to the limit in a naive pointwise way.
In the slow-bond model, for instance, a single microscopic defect gives rise to
a macroscopic variational problem with an effective rate at the origin, rather
than simply reproducing the microscopic rate
and the
question of whether an arbitrarily weak slow bond changes the macroscopic
current has been a longstanding and delicate problem
\cite{Seppalainen2001SlowBond, BasuSidoraviciusSly2014SlowBond}.

Hydrodynamic limits for spatially inhomogeneous particle systems have a substantial history beyond the homogeneous TASEP or the TASEP with microscopic defects. \cite{CovertRezakhanlou1997NonconstantSpeed} treated particle systems with nonconstant speed parameters, deriving the hydrodynamic equation in the continuously varying one-dimensional totally asymmetric case and obtaining bounds for various regimes.  \cite{Bahadoran1998HeterogeneousExclusion} proved hydrodynamic limits for spatially heterogeneous asymmetric exclusion processes with smooth macroscopic velocity profiles, while \cite{Rezakhanlou2002GrowthModels} derived Hamilton--Jacobi continuum limits for a class of growth models.

For TASEP with spatially discontinuous jump rates $c(i/n)$, the limiting profiles are again described by a variational formula, and in the two-phase case they agree with entropy solutions of the scalar conservation law
\[
        \rho_t+\bigl(c(x)\rho(1-\rho)\bigr)_x=0
\]
under the appropriate discontinuous-flux entropy condition \cite{GeorgiouKumarSeppalainen2010DiscontinuousTASEP}.  

Related hydrodynamic problems for zero-range and other attractive systems with discontinuous speed parameters also lead to conservation laws with discontinuous flux, and illustrate how a microscopic model can select a particular admissibility condition at the macroscopic interface \cite{ChenEvenKlingenberg2008Hydrodynamic}. 

Inhomogeneity can enter these models in several different ways, and the form of the inhomogeneity strongly affects what kind of macroscopic object one should expect.  One possibility is to keep the underlying lattice dynamics but allow the microscopic rates to vary with position.  For example, in exactly solvable inhomogeneous corner growth models, the exponential rate at $(i,j)$ may be chosen from row and column parameters, such as $a_i+b_j$, leading to variational formulas for the limit shape and, in some regimes, flat pieces or other non-strictly-concave behaviour \cite{Emrah2016LimitShapes,EmrahJanjigianSeppalainen2021Flats}.  Recent work on inhomogeneous corner growth further shows that these macroscopic irregularities are reflected in the geometry of Busemann functions, competition interfaces, and semi-infinite geodesics \cite{EmrahJanjigianSeppalainen2025Anomalous}. Exact formulas for inhomogeneous TASEP point in a complementary direction,
building on the broader determinantal/non-intersecting path tradition in solvable growth
and tiling models \cite{Johansson2002NonintersectingPaths}: particle- and time-inhomogeneous
discrete-time TASEP admits non-intersecting path and Fredholm determinant representations
\cite{BisiLiaoSaenzZygouras2023InhomogeneousTASEP}, while related TASEP variants with
inhomogeneous speeds and memory lengths can also be solved through explicit biorthogonal
kernels \cite{MatetskiRemenik2025InhomogeneousSpeeds}.

Macroscopic inhomogeneity in last-passage percolation has also been studied directly.  Rolla and Teixeira considered last-passage percolation in macroscopically inhomogeneous media and obtained a variational problem for the scaling limit, together with consequences for the associated TASEP picture \cite{RollaTeixeira2008InhomogeneousMedia}.  Calder treated directed last-passage percolation with discontinuous macroscopic weights through a continuum Hamilton--Jacobi formulation in the viscosity sense, also developing a dynamic-programming approach to the corresponding asymptotic shapes \cite{Calder2015DiscontinuousWeights}.

Another possibility is to change the geometry of the particle system itself: for instance, TASEP on Galton--Watson trees studies asymmetric exclusion away from a root, where the branching structure and edge rates replace the one-dimensional lattice as the source of inhomogeneous behaviour \cite{GantertGeorgiouSchmid2021TreeTASEP}. Hydrodynamic limits in disordered or nonstandard geometries provide another related direction: random-rate asymmetric systems and particlewise-disordered TASEP were studied in \cite{BenjaminiFerrariLandim1996RandomRates,SeppalainenKrug1999ParticlewiseDisorder}, quenched Euler hydrodynamics for attractive systems in random ergodic environments was developed in \cite{BahadoranGuiolRavishankarSaada2014RandomEnvironment}, diffusive exclusion processes with singular conductances were treated in
\cite{FrancoLandim2010} and multilane or ladder exclusion processes lead to scalar conservation laws coupled with relaxation limits for hyperbolic balance systems \cite{AmirBahadoranBusaniSaada2025MultilaneHydro,AmirBahadoranBusaniSaada2025Invariant}. A complementary weakly asymmetric direction studies fluctuations of exclusion
systems with spatial inhomogeneity around their macroscopic behaviour such as the SPDE limit for weakly inhomogeneous ASEP on the torus obtained in
\cite{CorwinTsai2020WeaklyInhomogeneousASEP}.

The inhomogeneity considered here is closer in spirit to directed last-passage percolation with a macroscopic speed function.  In \cite{ciech2021}, the exponential environment is no longer identically distributed: the rate at the lattice point $(i,j)$ is obtained from a discretisation of a function $c(i/n,j/n)$, where $c$ may have discontinuities along a locally finite family of curves.  The corresponding last-passage times satisfy a law of large numbers whose limit is given by the variational formula
\be \label{eq:LPPform}
        \Gamma_c(x,y)
        =
        \sup_{\mathbf{x}\in\mathcal H(x,y)}
        \int_0^1
        \frac{\gamma(\mathbf{x}'(s))}
             {c(\mathbf{x}(s))}
        \,ds,
\ee
where the supremum is over non-decreasing macroscopic paths from $(0,0)$ to
$(x,y)$. This result supplies the last-passage input for the present paper. The new step is to feed such a two-dimensional discontinuous environment through the TASEP height-function correspondence. The jump rate seen by the height function depends on both the spatial coordinate and the current level, producing a genuinely space-time inhomogeneous hydrodynamic problem. Thus the microscopic rates are fixed through the height-coordinate environment, but the macroscopic coefficient becomes space-time dependent once it is sampled along the evolving height profile.

\subsection{Discontinuous Hamilton--Jacobi equations and Scalar Conservation laws.}

On the PDE side, first-order equations with discontinuous coefficients have been studied from two closely related viewpoints: scalar conservation laws with discontinuous flux and Hamilton--Jacobi equations with discontinuous Hamiltonians. For scalar conservation laws, the presence of an interface destroys the direct applicability of the classical Kru{\v z}kov theory \cite{kruzkov1970} (see also \cite{Evans2010}), and uniqueness depends on choosing an admissibility condition that correctly encodes the behaviour across the discontinuity set. 

Early contributions such as \cite{AdimurthiGowda2003} already showed that discontinuous-flux problems require a more delicate entropy framework, while \cite{AudussePerthame2005} introduced adapted entropies and obtained uniqueness in a form that avoids stronger assumptions such as global BV estimates, convexity of the flux, or an explicit trace theory. 
This point of view was later developed in more systematic ways, including discontinuous-flux
formulations motivated by traffic and network models \cite{GaravelloNataliniPiccoliTerracina2007DiscontinuousFlux}:
the $L^1$-dissipative solver framework of \cite{AndreianovKarlsenRisebro2011} provides a robust uniqueness theory tied to admissibility and vanishing viscosity, while \cite{AndreianovMitrovic2015} revisits entropy conditions in a multidimensional setting and identifies classes of discontinuous-flux problems for which uniqueness can still be recovered. In a different but related direction, \cite{CrastaDeCiccoDePhilippis2015} proves uniqueness for BV solutions by means of a kinetic formulation.

For Hamilton--Jacobi equations, the issue is the validity of a comparison principle when the Hamiltonian is discontinuous in the space-time variables or across interfaces. Comparison and uniqueness were then established in several important settings, for example for eikonal-type equations with discontinuous right-hand side in \cite{DeckelnickElliott2004}, and for Hamilton--Jacobi equations with discontinuous coefficients in \cite{CocliteRisebro2007, DeZanSoravia2010}, where the discontinuity set is treated in the spirit of an internal boundary. Further comparison results for discontinuous Hamiltonians were obtained in \cite{GigaGorkaRybka2011}. 
% The comparison result of Chen and Hu \cite{ChenHu2008} is  relevant
% for the spatially discontinuous case considered in the uniqueness part of the present article. 
More geometric interface problems, especially those motivated by traffic and network models, were studied in \cite{ImbertMonneauZidani2013}, while the broader control-theoretic and stratified-domain viewpoint is surveyed in \cite{BarlesChasseigne2015}. 

The connection with traffic follows, since exclusion-type dynamics and related first-order
conservation laws have long served as idealized models for driven transport, traffic flow, and
self-driven particle systems \cite{ChowdhurySantenSchadschneider2000Traffic,Helbing2001Traffic}. Taken together, these works show that, in discontinuous first-order problems, uniqueness is inseparable from the choice of the correct notion of solution at the discontinuity set, whether this is expressed in entropy, kinetic, solver, or viscosity language.

The earliest discontinuous viscosity framework was introduced by Ishii \cite{Ishii1985}. This is a natural framework for
Hamilton--Jacobi equations with discontinuous Hamiltonians. In that framework,
subsolutions are tested against the upper semicontinuous envelope of the
Hamiltonian, while supersolutions are tested against the lower semicontinuous
envelope. This is appropriate for many comparison arguments, but it is not the only possible way to encode a discontinuity, and it does not by itself capture the selection imposed by the TASEP variational formula in our setting here.
In interface and network problems
one often has to prescribe the condition at the discontinuity set separately.
This is the point of view behind flux-limited viscosity solutions and related
junction conditions where the equation away from the interface is supplemented by an
interface Hamiltonian, or flux limiter, which represents the selection imposed
by the underlying variational, control, or microscopic problem \cite{CocliteRisebro2007, ImbertMonneau2017}.

In the present paper, we work within this general circle of ideas, but the structure of our equation and the way the discontinuity enters the Hamiltonian require a formulation tailored to our specific model rather than a direct application of the standard discontinuous-flux theories.
\subsection{Contributions} The starting point is the observation that the
discontinuous last-passage model of \cite{ciech2021} has a natural TASEP
counterpart.  In the usual spatially inhomogeneous TASEP, the jump rate depends only on the macroscopic
position of the particle.  Here, the rate depends on the position of the particle and on the height level of
the TASEP interface. 

The first main result is a hydrodynamic limit for this height-dependent, discontinuous TASEP.  Starting
from an arbitrary sequence of initial configurations with a macroscopic density profile, we prove that the
rescaled height functions converge to a deterministic current $v(x,t)$.  The associated particle density is
then obtained as the spatial derivative $\rho=v_x$, exactly as in the classical TASEP picture.  The proof
keeps the variational-coupling philosophy of homogeneous TASEP, but the auxiliary step-initial processes
now converge to level curves of a discontinuous, two-dimensional last-passage shape.  In this way the
law of large numbers for discontinuous-rate last-passage percolation is converted into a hydrodynamic
law of large numbers for the particle system.

The second main result is the variational description of the limiting current.  The microscopic envelope
formula for TASEP survives in the limit and gives
\be \label{eq:macroenv}
        v(x,t)=\sup_{q\in\mathbb R}\{v_0(q)-g^q(x-q,t)\},
\ee
where $g^q$ denotes the macroscopic level curve of the last-passage model started from the point
$(q,-v_0(q))$.  We also rewrite this representation in a Lax--Oleinik form.  This is useful for two
reasons: it identifies the variational principle selected by the microscopic dynamics, and it provides the
bridge from the probabilistic hydrodynamic limit to the Hamilton--Jacobi equation.

The third main result is the PDE identification. 
The LPP and particle-system descriptions use slightly different coordinates. Passing from the height-function coordinates to the last-passage coordinates amounts to the shear
\[
\tilde c(x,y)=c(x+y,y),
\]
where $c$ is used in the discontinuous last passage formulation \eqref{eq:LPPform}.
The coefficient sampled by the limiting current is then
\[
\tilde c(x,-v(x,t)).
\]
At points where the limiting current is differentiable and this sampled
coefficient is continuous, the current solves
\be\label{eq:limitngPDE}
        v_t+\tilde c(x,-v)v_x(1-v_x)=0.
\ee
At discontinuity points, however, the microscopic variational formula leads to viscosity
inequalities expressed through the lower and upper envelopes of the discontinuous coefficient.
This gives a natural envelope-based interpretation of the Hamilton--Jacobi equation across the
discontinuity set.

Equivalently, after scaling, the rate seen by the limiting height function is governed
by a two-dimensional speed function and appears through the composition
\[
        (x,t)\longmapsto \tilde c(x,-v(x,t)).
\]
This makes the model space-time inhomogeneous from the point of view of the evolving particle system, even though the randomness itself is still generated by independent Poisson clocks.

Our central point is twofold. We wish to obtain the limiting PDE \eqref{eq:limitngPDE}, and also identify the
discontinuous Hamiltonian selected by the same microscopic variational mechanism that gives the
hydrodynamic limit.
  The height process admits an envelope representation in terms of auxiliary processes with step initial data; these auxiliary processes are coupled to last-passage growth models, and the microscopic envelope becomes, after scaling, formula \eqref{eq:macroenv}. This variational-coupling approach, developed in particular by Sepp\"al\"ainen, is useful because it can prove hydrodynamic limits from a law of large numbers for the initial profile, without requiring an explicit description of invariant measures \cite{Seppalainen1999KExclusion,Seppalainen2001SlowBond}. 
  
  The limiting PDE \eqref{eq:limitngPDE} generalises the result in \cite{Rezakhanlou2002GrowthModels} where the macroscopic equation has the form
\[
u_t+\lambda(x,u)H(u_x)=0,
\]
with the clock rate depending on both position and height; there, the
coefficient is assumed bounded, uniformly positive, continuous, and Lipschitz in
the height variable. In the present work the analogous coefficient is sampled along the graph
$(x,-v(x,t))$, but the underlying speed function may be discontinuous and uses only local regularity conditions. 

Finally, when the coefficient is spatial only, $\tilde c(x,y)=\tilde c(x)$, we prove a uniqueness
statement.  In this case the limiting current is the unique viscosity solution in the natural class of
nondecreasing Lipschitz currents.  The corresponding density $\rho=v_x$ solves the scalar conservation
law
\[
        \rho_t+\bigl(\tilde c(x)\rho(1-\rho)\bigr)_x=0
\]
in the weak sense, and the TASEP solution is characterised by a maximal-current selection principle.
Thus, in the spatially discontinuous case, the microscopic particle dynamics select a distinguished weak
solution of the discontinuous-flux conservation law.

\subsection{Structure of the paper}
The rest of the paper is organised as follows. In Section \ref{sec:results} we introduce the discontinuous last-passage environment, construct the corresponding height-dependent TASEP sequences, and state the main hydrodynamic and PDE results. Section \ref{sec:Law of large numbers} proves the law of large numbers for the height functions by extending the variational coupling between TASEP and last-passage percolation to the present discontinuous setting. In Section \ref{sec:visco} we identify the Hamilton–Jacobi structure of the limiting current and prove that it is a selected viscosity solution for the discontinuous Hamiltonian chosen by the microscopic variational formula utilising a Lax--Oleinik representation, proven in \ref{sec:altrep}. Section \ref{sec:uniqueness} treats the spatially inhomogeneous case, where we prove uniqueness of the Hamilton–Jacobi solution and characterise the associated density through a maximal-current weak solution of the discontinuous-flux scalar conservation law. In Appendix \ref{app:LLN} we record supporting results for the discontinuous last-passage shape function that are not available in \cite{ciech2021}.

\subsection{Acknowledgments} 
Nicos Georgiou acknowledges partial support from the Dr Perry James (Jim) Browne Research Center at the Department of Mathematics, University of Sussex, and acknowledges financial
support provided by Sapienza University of Rome through the programme Professori Visitatori 2025.

Enrico Scalas acknowledges financial support provided by Sapienza University of Rome 000317 24 RICERCA UNIV
2023 PROG MEDI SCALAS - RICERCA ATENEO 2023 - SCALAS PROGETTI MEDI. The title of the
project is “Approximation of stochastic processes by means of sums of random telegraph processes”.

\section{Model and Results} 
\label{sec:results}

\subsection{Space-time discontinuous last passage percolation}

Our model breaks the identical distribution assumption by not assuming homogeneity in the rates of the environment. However, because we would like the shape function to exist even when the weights are non-stationary, we will pick rates for the realisations utilising a macroscopic function 
\[
c: \R^2 \to \R_{>0}.
\]   

For any $n \in \N$ the passage times themselves will be coupled through a common realisation of the exponential random variables $\{\tau_{i,j}\}$. However, for a fixed $n$, the rates of these random variables will be chosen according to an {\em a priori} discrete approximation of the function $c$ which we will call the {\em speed function} $c(x,y)$.

To be precise, for any $n \in \N$, we have that the mean of the environment variable $\tau^{(n)}_{i,j}$, for any $(i,j) \in \Z^2_+$ will be given by  
\be\label{eq:DSF}
r_{i,j}^{(n)} = c\Big(\frac{i}{n},\frac{j}{n}\Big)^{-1},\qquad (i,j)\in\bZ^2_+,
\ee
and the $n$-scaled, inhomogeneous environment is given by 
\be
\tau^{(n)}_{i,j} = r_{i,j}^{(n)} \tau_{i,j}.
\ee      

For ${\bf u},{\bf v} \in\Z^2_+$ and $n\in \mathbb{N}$ denote the last passage time from ${\bf u}$ to ${\bf v}$ by
\be\label{microLPT}
G^{(n)}_{\bf u, \bf v}=\max_{\pi\in\Pi_{\bf u, \bf v}}\bigg\{\sum_{(i,j)\in\pi}r_{i,j}^{(n)}\tau_{i,j}\bigg\}
=\max_{\pi\in\Pi_{\bf u, \bf v}}\bigg\{\sum_{(i,j)\in\pi}\tau_{i,j}^{(n)}\bigg\}.
\ee

In order to prove a limit theorem of type \eqref{thm:classicalLPP} in the absence of stationarity, a number of assumptions on the speed function $c(x,y)$ need to be imposed. In general, we want to allow $c(x,y)$ to be discontinuous so we impose some regularity conditions on the type of discontinuity curves it is allowed to have.

\begin{assumption}[Discontinuity curves of $c(x,y)$] 
\label{ass:c}
	Function $c(x,y)$ is discontinuous on a (potentially) countable set of curves $H_c = \{ h_i\}_{i \in \mathcal I}$ satisfying the following properties 
		\begin{enumerate}
			\item $h_i$ is either a linear segment or strictly monotone.  
			\item All $h_i$ can be viewed as  graphs			
			\[
			h_i: [ z_i, w_i] = {\rm Dom}(h_i) \to \R,
			\] 
			\item If $h_i$ is strictly increasing, then 
			\begin{enumerate}
				\item $h_i$ is $C^1((z_i, w_i), \R)$. At the boundary points $z_i, w_i$ the derivative may take the value $\pm \infty, 0$.
				\item The equation $h_i'(s) = 0$ has finitely many solutions in $[ z_i, w_i]$. 
			\end{enumerate}
				\item If $h_i$ is strictly decreasing, then $h_i$ is continuous.					
				\item There are finitely many curves $h_i$ in any compact set $K \subseteq \bR^2_+$, satisfying $\it (1)-(4)$. Equivalently, accumulation points of different curves $\{h_j\}_{j}$ are not allowed. 
		\end{enumerate} 
\end{assumption}			

There are two types of points on these discontinuity curves,
\begin{enumerate}
	\item(Interior points) These are points ${\bf w}$ that belong to a single discontinuity curve $h_i$. For any point ${\bf w}$ of this form, we can find an $\e > 0$ so that $h_i$ partitions the ball $ B( {\bf w}, \e) $ into three disjoint sets,  $U_{\e, \bf w }$ (above $h_i$),  $L_{\e, \bf w }$ (below $h_i$) and 
$(h_i \cap B( {\bf w}, \e))$. 
	
	\item(Intersection/terminal points) These are points ${\bf w}$ that either belong on more than one discontinuity curve or they are terminal for $h_i$. There are finitely many of these points in any compact set.			
\end{enumerate}

The discontinuity curves $\{h_i\}_{i \in \mathcal I}$ separate $\R^2_+$ into open regions in which $c(x,y)$ is assumed continuous. The number of regions is finite in any compact set of $\R^2_+$.	Denote the set of regions by $\mathcal Q$. 

We also need to record some assumptions on the behaviour of $c(x,y)$ in $\mathcal Q$.

\begin{assumption}[Regularity properties of $c(x,y)$]		
\label{ass:c2}
\phantom{x}
\begin{enumerate}		
		\item $c(x,y)$ is continuous on any $Q \in \mathcal Q$,  lower-semicontinuous on $\R^2$, that further satisfies the following  stability assumption:  
		
		For every $i \in \mathcal I$ and interior point ${ \bf w } \in h_i$, there exists $\e = \e(i, {\bf w})>0$ so that for all ${\bf y} \in B( {\bf w}, \e) \cap h_i $ 
		there exists an open set $Q_{i, \bf w} \in \{ L_{\e, \bf w }, U_{\e, \bf w }\}$, so that 
		for any sequence ${\bf z}_n \in Q_{i, \bf w} \cap B( {\bf w}, \e)$ with ${\bf z}_n \to {\bf y}$,
		\be \label{eq:calder2}
		 \lim_{{{\bf z}_n} \to {\bf y}} c({\bf z}_n) = c( {\bf y}). 
		\ee
		\item For any compact set $K \subset \R^2_+$, there exist two constants $r_{\text{\rm low}}^{(K)} > 0$ and $r_{\text{\rm high}}^{(K)} < \infty$, so that 
		\[
			r_{\text{\rm low}}^{(K)} \le c(x,y) \le r_{\text{\rm high}}^{(K)} , \quad \forall (x, y) \in K.
		\]
	\end{enumerate}
\end{assumption}

Under Assumptions \ref{ass:c} and \ref{ass:c2}, the following law of large numbers holds.
\begin{theorem}[\cite{ciech2021},Theorem 2.6]
\label{thm:1}
Recall \eqref{microLPT}. Let $c(x,y)$ be a macroscopic speed function which satisfies Assumptions \ref{ass:c} and \ref{ass:c2}, and let ${\bf v} =(x,y) \in \R^2_+$. Define ${\bf v}_n = (\ce{nx}, \ce{ny})$. Then we have the scaling limit
\begin{equation}\label{eq:6}
\lim_{n\to\infty}n^{-1}G^{(n)}_{{\bf 0}, {\bf v}_n}= \Gamma_c({\bf v}) = \Gamma_c(x,y)\quad \P-\text{a.s.}
\end{equation}
where $\Gamma_c(x,y)$ is upper-semicontinuous on $\R_+^2$.
\end{theorem}

Function $\Gamma_c(x,y)$ has a variational formulation that mimics the microscopic last passage time. Namely
for a fixed $(x,y)$ in $\R^{2}_+$ and a speed function $c(\cdot, \cdot)$,
$\Gamma_c(x,y)$ is defined via the variational formula 
\begin{equation}\label{macroLPT}
\Gamma_c(x,y)=\sup_{\mathbf{x}(\cdot)\in\mathcal{H}(x,y)}\bigg\{\int_0^1\frac{\gamma(\mathbf{x}'(s))}{c(x_1(s),x_2(s))}ds\bigg\},
\end{equation}
where $\gamma(x,y)$ is from \eqref{eq:limshap}, $\mathbf{x}(s)=(x_1(s),x_2(s))$ denotes a path in $\mathbb{R}^2$ and 
\begin{align}
\mathcal{H}(x,y) &=\{{\bf x}\in C([0,1],\R^2_+):\mathbf{x} \text{ is piecewise }C^1,\mathbf{x}(0)=(0,0), \mathbf{x}(1)=(x,y), \notag \\
& \phantom{xxxxxxxxxxxxxxxxxxxxxxxx}\mathbf{x}'(s)\in \R^2_+ \text{ wherever the derivative is defined}\}. \label{eq:pathsH}
\end{align}

In this paper, we make the following assumption, to simplify a lot of the arguments.
\begin{assumption}\label{ass:c0}
The discontinuity curves $h_i$ of $c(x,y)$ cannot be (or contain) a vertical or a horizontal line segment. 
\end{assumption}

Under Assumption \ref{ass:c0}, the following holds, we can further restrict the set of admissible paths in the supremum above. 

 Let $\mathcal C$ denote the set of continuity points of the speed function $c(x,y)$,
\[
\mathcal C = \{ (x,y) \in \R^2:  c(x,y) \text{ is continuous at } (x,y) \}.
\]
By Assumptions \ref{ass:c}, \ref{ass:c2}, $\mathcal C$ has full measure in any subset of $\R^2$ and it is dense in $\R^2$. In particular, $\mathcal C$ has a full measure as a subset  $\mathcal C \subseteq {\rm Leb}(c)$.

\begin{lemma}[Changing the path space] 
\label{lem:betterpaths}
Assume \ref{ass:c}, \ref{ass:c2} and \ref{ass:c0}. The set of paths $\mathcal H(x,y)$ used in equation \eqref{macroLPT} can be changed to 
 \begin{align*}
\mathcal{H}^{\rm cont}(x,y) &=\{{\bf x}\in C([0,1],\R^2_+):\mathbf{x} \text{ is piecewise }C^1,\mathbf{x}(0)=(0,0), \mathbf{x}(1)=(x,y), \\
& \phantom{xxxxxxxxx} {\bf x }(s) \in \mathcal C \text{ for a.e. }  s\in [0,1], \,  \mathbf{x}'(s)\in \R^2_+ \text{ wherever the derivative is defined}\},
\end{align*}
without changing the value of the supremum in \eqref{macroLPT}.
\end{lemma}

\begin{theorem}[\cite{ciech2021}, Theorem 2.4] \label{thm:2.5}
Assume \ref{ass:c}, \ref{ass:c2} and \ref{ass:c0}.The shape function $\Gamma_c(x,y)$ in Theorem \ref{thm:1} is continuous.
\end{theorem}

\subsection{The space-time discontinuous TASEP} 

As mentioned earlier, here our goal is to prove hydrodynamic limits for the version of TASEP which is dual to this general inhomogeneous LPP model, and highlight its connection to the theory of Hamilton--Jacobi equation and scalar conservation laws with discontinuous fluxes.  

For this reason we prove the limiting profiles as solutions on all of $\mathbb R$ rather than on a compact set (e.g. the torus) because the probabilistic techniques - at least in this particular example - can offer a new perspective and contribute to the existing PDE theory.

Here we will assume that $c(x,y)$ grows according to the following assumption.

\begin{assumption}[Logarithmic growth conditions]
\label{assumption: growth rates} Fix any vector $\bf a \in \R^2$ and let $(x,y) \in \R^2$.  Define the rectangle  
\be\label{eq:rectangle} R^{\bf a}_{x,y} = \{ {\bf u} \in \bR^2: \| {\bf u} -{\bf a} \|_{\infty} \le  \| (x,y) -{\bf a} \|_{\infty} \},
\ee
and let
\[
M^{\bf a}_{x,y} = \sup_{ {\bf u} \in  R^{\bf a}_{x,y}  } c(\bf u).
\] 
We impose the following growth conditions for the speed function $c(x,y)$: for any fixed $ \bf a \in \bR^2$ 
     \begin{equation}\label{eq:2.7}
     \limsup_{|x|\vee |y| \to \infty}\frac{ M^{\bf a}_{x,y} }{\log(|x|\vee|y|)} \longrightarrow 0.
 \end{equation}
\end{assumption}

The reason for the sub-logarithmic growth is technical, and we remark later at the places
where it is needed, the graphical construction of the TASEP process. 
%Indeed, with the exception of this technicality, all the remaining arguments in the article can be pushed to the case in which $c(x,y)$ has sublinear growth. 

To define the TASEP, let, for any $n \in \N$ and $t > 0$ the configuration process be $\eta^n(t) = \{ \eta^n_x(t) \}_{x \in \Z}$ on $\{ 0, 1 \}^{\Z}$ given by  
\[
\eta^n_x(t) = \mathds1\{ \text{a particle exists at location $x$ at time $t$ }\}.
\]
At time $t = 0$, we start with an arbitrary  (random or deterministic) particle configuration $\{\eta_x(0)\}_{x \in \Z}$ and we can find a  non-negative measurable function $\rho_0$ with support in $\R$ and with $0 \le \rho_0 \le 1$, so that for any real $a < b$ the following initial law of large numbers holds. 
\be\label{eq:etaass} 
\lim_{n\to \infty} \frac{1}{n}\sum_{i = \fl{na}}^{\fl{nb}} \eta^n_i(0) = \int_{a}^b \rho_0(x) \, dx.
\ee
This LLN can be either weak or strong, but depending on which choice we make, the results in the theorems that follow will have the same regularity as this condition. 

Define the current function $v_0(b)$ by  
\be \label{eq:in-cur}
v_0(b) =  \int_{0}^b \rho_0(x) \, dx + v_0(0),
\ee
where $v_0(0)$ is a constant chosen arbitrarily. 
Because $\rho_0 \in [0,1]$, we have that $\displaystyle \int_{a}^{b}  \rho_0(x) \, dx =  v_0(b) - v_0(a) \le b-a$ for any $b > a$ and as such, $v_0$ is Lipschitz-1 continuous, with an almost everywhere derivative equal to $\rho_0$. The function $v_0$ can be viewed as the law of large numbers for a random process called the {\em height function}. In order to construct it, we use the configuration process $\eta^n$.

To define the height function $z^n$ at time $t = 0$,  $z^n(0) = \{z^{n}_i(0)\}_{i \in \Z}$ arbitrarily set the value of $z^n_0(0)$ (it can be any finite integer without loss of generality), subject to the constraint  
\be \label{eq:z0lim}
\lim_{n\to \infty} \frac{1}{n}z^n_0(0)=v_0(0).
\ee 
Then we use the initial configuration of the TASEP particles  $\{\eta_{i}^n(0)\}_{i \in \Z}$ to define $z_i^n(0)$ for all $i \in \Z\setminus\{0\}$ as  
\begin{equation} \label{Server-exclusion relation0}
    \eta^{n}_{i-1}(0) = z^{n}_{i}(0) - z_{i-1}^{n}(0).
\end{equation}
Note that this implies there is a particle at position $i-1$ if and only if $z^{n}_{i}(0) = 1+ z_{i-1}^{n}(0)$. 
Moreover, since consecutive differences are bounded by 1, we have that for any $k, \ell \in \Z$ we have the inequality 
\be \label{eq:DiscLip}
|z^n_k(0) - z^n_\ell(0)| \le |k-\ell|
\ee

By telescoping \eqref{Server-exclusion relation0}, we see that 
\be \label{eq:zic}
\lim_{n\to \infty} \frac{1}{n}(z^n_{\fl{nx}}(0) -z^n_0(0))= \lim_{n\to \infty} \frac{1}{n}\sum_{i = 0}^{\fl{nx}} \eta_i(0) = v_0(x) -v_0(0).
\ee

Define a {\em growth site} at time $0$ to be $(i, z_i^n(0))$ such that $ z^{n}_{i}(0) - z_{i-1}^{n}(0) = 1$ and  $z^{n}_{i+1}(0) - z_{i}^{n}(0) = 0$. 
These are places where the particle at location $i-1$ will be able to move one unit to the right when it wants to jump. In words, growth sites at time 0 are particle-hole pairs, places where particles are allowed to perform their first jump.   
The set of growth sites at time $0$ is denoted by $\mathcal G(0)$. 

The evolution of all processes is coupled by a field of mutually independent Poisson processes of rate 1, denoted by $\{ N_{i,j}(t)\}_{(i,j) \in \Z^2}$ and with $\bf P$ denotes their joint law. 

What changes with $n \in \N$ is the jump rate $\lambda^{n}_{i,j}$ of the independent Poisson processes $\{N^{(n)}_{i,j}\}_{(i, j) \in \Z^2}$ that are defined as time-changed versions of the original field, namely  
\[
N^{(n)}_{i,j}(t) = N_{i,j}(\lambda^{n}_{i,j} t).
\]
The only assumption we are making is that there is a background macroscopic function $ c(x,y): \R^2 \to \R_{>0}$ satisfying Assumptions \ref{ass:c}, \ref{ass:c2}, \ref{ass:c0} and \ref{assumption: growth rates} so that
\be\label{eq:lambda-disc} 
\lambda^{n}_{i,j} = c ((i + j)n^{-1}, jn^{-1}).
\ee
The coupling with an appropriate last passage model utilising the speed function $c(x,y)$ will be made clear in the next section. For now, we note that we will need to distinguish between the LPP environment utilising $c(x,y)$ versus the Poisson rates that are using the function 
\be\label{eq:tildec}
\tilde c(x,y) = c(x+y,y). 
\ee

The function $c$ is the speed function in the LPP coordinates, while
$\tilde c$ is the same environment written as will be used by the particle system 
coordinates. Throughout the paper these two functions are related by the
linear shear
\[
        \tilde c(x,y)=c(x+y,y),
        \qquad
        c(u,y)=\tilde c(u-y,y).
\]
All assumptions inherited from the discontinuous LPP law of large numbers
are imposed on $c$. All hydrodynamic and PDE statements are expressed via $\tilde c$.

For each $n \in \mathbb{N}$ and $t \in \R_{>0}$, we construct a $\mathbb{Z}$-valued height process $z^{n}(t) = (z^{n}_{i}(t): i \in \mathbb{Z})$ that evolves in parallel with the TASEP particles $\eta^n$. As a process, the height function  is {\bf decreasing in time} coordinate-wise, while obeying the constraint 
\begin{equation} \label{Server process constraint}
    0 \leq z^{n}_{i}(t) - z^{n}_{i-1}(t) \leq 1,
\end{equation}
for all $t \in \R_{>0}$ and taking values in $\Z$.  

In the case where $z^{n}_{i}(t) - z^{n}_{i-1}(t) = 1$, then we say that a TASEP particle sits at location $i-1$, and otherwise the location is empty i.e. for all times, the TASEP particle configuration is $\{\eta_i(t)\}_{i \in \Z}$ is given by 
\be \label{eq:parev}
z^{n}_{i}(t) - z^{n}_{i-1}(t) = \eta_{i-1}(t). 
\ee

Starting from $t=0$, the process $z$ attempts to jump one step down at all locations $i$ of growth sites, independently of anything, at the jump time of the Poisson $N^{(n)}_{i, -z_i(0)}$. Fix a time horizon $T$. Then by Assumption \ref{assumption: growth rates}, we can find arbitrarily large (in modulus) integers $i$ such that
\begin{equation}
    {\bf P}\{ N^{(n)}_{i, -z_i^n(0)}(T)= 0 \text{ for infinitely many } i \} = 1.
\end{equation}
To see this, use Assumption \ref{assumption: growth rates} to argue that we have $M_{i, -z_i^n(0)}^{\bf 0} < \delta \log(|i|\vee| -z_i^n(0)|+1) $ for any $\delta>0$ and all  $i$ large enough. Then by \eqref{eq:DiscLip},
\begin{equation*}
    {\bf P}\{N^{(n)}_{i, -z_i^n(0)}(T) = 0\}  \geq e^{- T M^{\bf 0}_{i,  -z_i^n(0)}} > e^{- \delta T \log(|i|\vee| z_i^n(0)|+1)} \ge (|i|+|z_0(0)|+1)^{- \delta T}.
\end{equation*}
This is the only place where the sublogarithmic growth condition in Assumption \ref{assumption: growth rates} is used. Now for $ \delta  < 1/T$ the sum of these probabilities over $i$ diverges, while the events are independent across the $i$ index. The statement then follows from the second Borel-Cantelli Lemma, for eventually all $n$. 
In particular, this implies that, for any time $T$, we can find arbitrarily large (in modulus) indices $L <0, U>0$ and random-sized boxes $B_{L,U} \subset \Z^2$, around $0$ so that 
\[B_{L,U} \subseteq [L,U]\times [-\max\{ U, |L| \} - |z^n_0(0)|, \max\{ U, |L| \} + |z^n_0(0)|].\]
In $B_{L,U}$, the evolution 
of the height process $z^n$ is not affected by Poisson $\{N^{(n)}_{i,j}\}_{(i, j) \in \Z^2}$ with $ (i,j) \notin B_{L,U}$.

\begin{figure}[ht]
% Gradient Info
  \tikzset {_dlyrgacgu/.code = {\pgfsetadditionalshadetransform{ \pgftransformshift{\pgfpoint{89.1 bp } { -128.7 bp }  }  \pgftransformscale{1.32 }  }}}
\pgfdeclareradialshading{_0a952aycu}{\pgfpoint{-72bp}{104bp}}{rgb(0bp)=(1,1,1);
rgb(2.5bp)=(1,1,1);
rgb(14.910714285714285bp)=(0.56,0.07,1);
rgb(400bp)=(0.56,0.07,1)}
\tikzset{_n0dvk4cfx/.code = {\pgfsetadditionalshadetransform{\pgftransformshift{\pgfpoint{89.1 bp } { -128.7 bp }  }  \pgftransformscale{1.32 } }}}
\pgfdeclareradialshading{_c9z6lgmt0} { \pgfpoint{-72bp} {104bp}} {color(0bp)=(transparent!0);
color(2.5bp)=(transparent!0);
color(14.910714285714285bp)=(transparent!37);
color(400bp)=(transparent!37)} 
\pgfdeclarefading{_rdz2r7nl8}{\tikz \fill[shading=_c9z6lgmt0,_n0dvk4cfx] (0,0) rectangle (50bp,50bp); } 

% Gradient Info
  
\tikzset {_u16zh8rhg/.code = {\pgfsetadditionalshadetransform{ \pgftransformshift{\pgfpoint{89.1 bp } { -128.7 bp }  }  \pgftransformscale{1.32 }  }}}
\pgfdeclareradialshading{_an1dxvxiw}{\pgfpoint{-72bp}{104bp}}{rgb(0bp)=(1,1,1);
rgb(2.5bp)=(1,1,1);
rgb(14.910714285714285bp)=(0.56,0.07,1);
rgb(400bp)=(0.56,0.07,1)}
\tikzset{_m96l3k2zt/.code = {\pgfsetadditionalshadetransform{\pgftransformshift{\pgfpoint{89.1 bp } { -128.7 bp }  }  \pgftransformscale{1.32 } }}}
\pgfdeclareradialshading{_9a6lqn3w7} { \pgfpoint{-72bp} {104bp}} {color(0bp)=(transparent!0);
color(2.5bp)=(transparent!0);
color(14.910714285714285bp)=(transparent!37);
color(400bp)=(transparent!37)} 
\pgfdeclarefading{_3jqkjquc1}{\tikz \fill[shading=_9a6lqn3w7,_m96l3k2zt] (0,0) rectangle (50bp,50bp); } 

% Gradient Info
  
\tikzset {_330i92zur/.code = {\pgfsetadditionalshadetransform{ \pgftransformshift{\pgfpoint{89.1 bp } { -128.7 bp }  }  \pgftransformscale{1.32 }  }}}
\pgfdeclareradialshading{_1lci1pl48}{\pgfpoint{-72bp}{104bp}}{rgb(0bp)=(1,1,1);
rgb(2.5bp)=(1,1,1);
rgb(14.910714285714285bp)=(0.56,0.07,1);
rgb(400bp)=(0.56,0.07,1)}
\tikzset{_oi9gq4szt/.code = {\pgfsetadditionalshadetransform{\pgftransformshift{\pgfpoint{89.1 bp } { -128.7 bp }  }  \pgftransformscale{1.32 } }}}
\pgfdeclareradialshading{_bhbzz3noh} { \pgfpoint{-72bp} {104bp}} {color(0bp)=(transparent!0);
color(2.5bp)=(transparent!0);
color(14.910714285714285bp)=(transparent!37);
color(400bp)=(transparent!37)} 
\pgfdeclarefading{_qfjfhwqkz}{\tikz \fill[shading=_bhbzz3noh,_oi9gq4szt] (0,0) rectangle (50bp,50bp); } 

% Gradient Info
  
\tikzset {_giqcsqxib/.code = {\pgfsetadditionalshadetransform{ \pgftransformshift{\pgfpoint{89.1 bp } { -128.7 bp }  }  \pgftransformscale{1.32 }  }}}
\pgfdeclareradialshading{_nmuduffyl}{\pgfpoint{-72bp}{104bp}}{rgb(0bp)=(1,1,1);
rgb(2.5bp)=(1,1,1);
rgb(14.910714285714285bp)=(0.56,0.07,1);
rgb(400bp)=(0.56,0.07,1)}
\tikzset{_9s1uu2kao/.code = {\pgfsetadditionalshadetransform{\pgftransformshift{\pgfpoint{89.1 bp } { -128.7 bp }  }  \pgftransformscale{1.32 } }}}
\pgfdeclareradialshading{_issr3qjft} { \pgfpoint{-72bp} {104bp}} {color(0bp)=(transparent!0);
color(2.5bp)=(transparent!0);
color(14.910714285714285bp)=(transparent!37);
color(400bp)=(transparent!37)} 
\pgfdeclarefading{_2l7obxm54}{\tikz \fill[shading=_issr3qjft,_9s1uu2kao] (0,0) rectangle (50bp,50bp); } 

% Gradient Info
  
\tikzset {_x5xws0sb9/.code = {\pgfsetadditionalshadetransform{ \pgftransformshift{\pgfpoint{89.1 bp } { -128.7 bp }  }  \pgftransformscale{1.32 }  }}}
\pgfdeclareradialshading{_e5h2195tt}{\pgfpoint{-72bp}{104bp}}{rgb(0bp)=(1,1,1);
rgb(2.5bp)=(1,1,1);
rgb(14.910714285714285bp)=(0.56,0.07,1);
rgb(400bp)=(0.56,0.07,1)}
\tikzset{_kb6sbn4o0/.code = {\pgfsetadditionalshadetransform{\pgftransformshift{\pgfpoint{89.1 bp } { -128.7 bp }  }  \pgftransformscale{1.32 } }}}
\pgfdeclareradialshading{_6j5ch2azu} { \pgfpoint{-72bp} {104bp}} {color(0bp)=(transparent!0);
color(2.5bp)=(transparent!0);
color(14.910714285714285bp)=(transparent!37);
color(400bp)=(transparent!37)} 
\pgfdeclarefading{_zvs913mvf}{\tikz \fill[shading=_6j5ch2azu,_kb6sbn4o0] (0,0) rectangle (50bp,50bp); } 

% Gradient Info
  
\tikzset {_s4wfgjn8z/.code = {\pgfsetadditionalshadetransform{ \pgftransformshift{\pgfpoint{89.1 bp } { -128.7 bp }  }  \pgftransformscale{1.32 }  }}}
\pgfdeclareradialshading{_9d2mgtb6j}{\pgfpoint{-72bp}{104bp}}{rgb(0bp)=(1,1,1);
rgb(2.5bp)=(1,1,1);
rgb(14.910714285714285bp)=(0.56,0.07,1);
rgb(400bp)=(0.56,0.07,1)}
\tikzset{_o7o322o5n/.code = {\pgfsetadditionalshadetransform{\pgftransformshift{\pgfpoint{89.1 bp } { -128.7 bp }  }  \pgftransformscale{1.32 } }}}
\pgfdeclareradialshading{_0nywkzaop} { \pgfpoint{-72bp} {104bp}} {color(0bp)=(transparent!0);
color(2.5bp)=(transparent!0);
color(14.910714285714285bp)=(transparent!37);
color(400bp)=(transparent!37)} 
\pgfdeclarefading{_w29ut9usx}{\tikz \fill[shading=_0nywkzaop,_o7o322o5n] (0,0) rectangle (50bp,50bp); } 

% Gradient Info
  
\tikzset {_n5s869br7/.code = {\pgfsetadditionalshadetransform{ \pgftransformshift{\pgfpoint{89.1 bp } { -128.7 bp }  }  \pgftransformscale{1.32 }  }}}
\pgfdeclareradialshading{_65q9ps644}{\pgfpoint{-72bp}{104bp}}{rgb(0bp)=(1,1,1);
rgb(2.5bp)=(1,1,1);
rgb(14.910714285714285bp)=(0.56,0.07,1);
rgb(400bp)=(0.56,0.07,1)}
\tikzset{_z0xcrlncm/.code = {\pgfsetadditionalshadetransform{\pgftransformshift{\pgfpoint{89.1 bp } { -128.7 bp }  }  \pgftransformscale{1.32 } }}}
\pgfdeclareradialshading{_p5qgb16lj} { \pgfpoint{-72bp} {104bp}} {color(0bp)=(transparent!0);
color(2.5bp)=(transparent!0);
color(14.910714285714285bp)=(transparent!37);
color(400bp)=(transparent!37)} 
\pgfdeclarefading{_a01hbpe3c}{\tikz \fill[shading=_p5qgb16lj,_z0xcrlncm] (0,0) rectangle (50bp,50bp); } 

% Gradient Info
  
\tikzset {_cd0nmm8wa/.code = {\pgfsetadditionalshadetransform{ \pgftransformshift{\pgfpoint{89.1 bp } { -128.7 bp }  }  \pgftransformscale{1.32 }  }}}
\pgfdeclareradialshading{_urlhw29dz}{\pgfpoint{-72bp}{104bp}}{rgb(0bp)=(1,1,1);
rgb(2.5bp)=(1,1,1);
rgb(14.910714285714285bp)=(0.56,0.07,1);
rgb(400bp)=(0.56,0.07,1)}
\tikzset{_11w77lqlh/.code = {\pgfsetadditionalshadetransform{\pgftransformshift{\pgfpoint{89.1 bp } { -128.7 bp }  }  \pgftransformscale{1.32 } }}}
\pgfdeclareradialshading{_d8fhb350k} { \pgfpoint{-72bp} {104bp}} {color(0bp)=(transparent!0);
color(2.5bp)=(transparent!0);
color(14.910714285714285bp)=(transparent!37);
color(400bp)=(transparent!37)} 
\pgfdeclarefading{_3abpagegl}{\tikz \fill[shading=_d8fhb350k,_11w77lqlh] (0,0) rectangle (50bp,50bp); } 

% Gradient Info
  
\tikzset {_kkjo516ck/.code = {\pgfsetadditionalshadetransform{ \pgftransformshift{\pgfpoint{89.1 bp } { -128.7 bp }  }  \pgftransformscale{1.32 }  }}}
\pgfdeclareradialshading{_drvzrahfu}{\pgfpoint{-72bp}{104bp}}{rgb(0bp)=(1,1,1);
rgb(2.5bp)=(1,1,1);
rgb(14.910714285714285bp)=(0.56,0.07,1);
rgb(400bp)=(0.56,0.07,1)}
\tikzset{_hteykxvvd/.code = {\pgfsetadditionalshadetransform{\pgftransformshift{\pgfpoint{89.1 bp } { -128.7 bp }  }  \pgftransformscale{1.32 } }}}
\pgfdeclareradialshading{_4nlf86idl} { \pgfpoint{-72bp} {104bp}} {color(0bp)=(transparent!0);
color(2.5bp)=(transparent!0);
color(14.910714285714285bp)=(transparent!37);
color(400bp)=(transparent!37)} 
\pgfdeclarefading{_w5az71zv7}{\tikz \fill[shading=_4nlf86idl,_hteykxvvd] (0,0) rectangle (50bp,50bp); } 

% Gradient Info
  
\tikzset {_zcqapn3su/.code = {\pgfsetadditionalshadetransform{ \pgftransformshift{\pgfpoint{89.1 bp } { -128.7 bp }  }  \pgftransformscale{1.32 }  }}}
\pgfdeclareradialshading{_907n669nz}{\pgfpoint{-72bp}{104bp}}{rgb(0bp)=(1,1,1);
rgb(2.5bp)=(1,1,1);
rgb(14.910714285714285bp)=(0.56,0.07,1);
rgb(400bp)=(0.56,0.07,1)}
\tikzset{_92j1bar5l/.code = {\pgfsetadditionalshadetransform{\pgftransformshift{\pgfpoint{89.1 bp } { -128.7 bp }  }  \pgftransformscale{1.32 } }}}
\pgfdeclareradialshading{_n5ni68por} { \pgfpoint{-72bp} {104bp}} {color(0bp)=(transparent!0);
color(2.5bp)=(transparent!0);
color(14.910714285714285bp)=(transparent!37);
color(400bp)=(transparent!37)} 
\pgfdeclarefading{_t43uh97cw}{\tikz \fill[shading=_n5ni68por,_92j1bar5l] (0,0) rectangle (50bp,50bp); } 

% Gradient Info
  
\tikzset {_dpr19xkn2/.code = {\pgfsetadditionalshadetransform{ \pgftransformshift{\pgfpoint{89.1 bp } { -128.7 bp }  }  \pgftransformscale{1.32 }  }}}
\pgfdeclareradialshading{_itvpsnfvz}{\pgfpoint{-72bp}{104bp}}{rgb(0bp)=(1,1,1);
rgb(2.5bp)=(1,1,1);
rgb(14.910714285714285bp)=(0.56,0.07,1);
rgb(400bp)=(0.56,0.07,1)}
\tikzset{_a8aan6uk9/.code = {\pgfsetadditionalshadetransform{\pgftransformshift{\pgfpoint{89.1 bp } { -128.7 bp }  }  \pgftransformscale{1.32 } }}}
\pgfdeclareradialshading{_um6cyef2t} { \pgfpoint{-72bp} {104bp}} {color(0bp)=(transparent!0);
color(2.5bp)=(transparent!0);
color(14.910714285714285bp)=(transparent!37);
color(400bp)=(transparent!37)} 
\pgfdeclarefading{_qu4qjuiu4}{\tikz \fill[shading=_um6cyef2t,_a8aan6uk9] (0,0) rectangle (50bp,50bp); } 
\tikzset{every picture/.style={line width=0.75pt}} %set default line width to 0.75pt        

\begin{tikzpicture}[x=0.75pt,y=0.75pt,yscale=-1,xscale=1]
%uncomment if require: \path (0,1878); %set diagram left start at 0, and has height of 1878

%Straight Lines [id:da4565396875010783] 
\draw    (41.07,1657.29) -- (121.72,1657.29) ;
%Straight Lines [id:da26256424611154305] 
\draw    (121.72,1636.86) -- (182.36,1636.86) ;
%Straight Lines [id:da4779534489675674] 
\draw    (182.36,1616.43) -- (251,1616.43) ;
%Straight Lines [id:da8334614155592197] 
\draw    (251,1596) -- (269.64,1596) ;
%Straight Lines [id:da48967152228015187] 
\draw    (270,1578) -- (292.32,1578) ;
%Straight Lines [id:da27270811585508126] 
\draw    (291,1560) -- (329.64,1560.43) ;
%Straight Lines [id:da5189644022665725] 
\draw    (330,1542) -- (353.36,1541.86) ;
%Straight Lines [id:da4017697770773805] 
\draw    (330,1542) -- (329.64,1560.43) ;
%Straight Lines [id:da3049338756673651] 
\draw    (291,1560) -- (291,1577.68) ;
%Straight Lines [id:da11840507240339215] 
\draw    (269.64,1578) -- (269.64,1596.43) ;
%Straight Lines [id:da12248350595801583] 
\draw    (251,1596) -- (251,1616.43) ;
%Straight Lines [id:da6211676233406802] 
\draw    (182.36,1616.43) -- (182.36,1636.86) ;
%Straight Lines [id:da6860562629608958] 
\draw    (121.72,1636.86) -- (121.72,1657.29) ;
%Straight Lines [id:da445227003458857] 
\draw    (353.36,1521.43) -- (422,1521.43) ;
%Straight Lines [id:da24150779500889896] 
\draw    (422,1501) -- (441.64,1501.43) ;
%Straight Lines [id:da5412672758922145] 
\draw    (441,1483) -- (501.64,1483) ;
%Straight Lines [id:da014104486199834843] 
\draw    (502,1465) -- (540.64,1465.43) ;
%Straight Lines [id:da47777036832939923] 
\draw    (541,1447) -- (560.64,1447.43) ;
%Straight Lines [id:da8266554324048416] 
\draw    (541,1447) -- (540.64,1465.43) ;
%Straight Lines [id:da40951338851801566] 
\draw    (502,1465) -- (501.64,1483.43) ;
%Straight Lines [id:da5027824272547646] 
\draw    (441,1483) -- (440.64,1501.43) ;
%Straight Lines [id:da3855795500675626] 
\draw    (422,1501) -- (422,1521.43) ;
%Straight Lines [id:da811104038639862] 
\draw    (353.36,1521.43) -- (353.36,1541.86) ;
%Straight Lines [id:da613844404915313] 
\draw    (37.64,1387) -- (567.64,1387) ;
%Straight Lines [id:da4964331764751436] 
\draw  [dash pattern={on 4.5pt off 4.5pt}]  (110,1389) -- (111.37,1651.42) ;
%Straight Lines [id:da8695440363903462] 
\draw  [dash pattern={on 4.5pt off 4.5pt}]  (172.5,1386) -- (172.5,1630) ;
%Straight Lines [id:da14312395691961544] 
\draw  [dash pattern={on 4.5pt off 4.5pt}]  (241,1387) -- (241.5,1609) ;
%Straight Lines [id:da6154503891991706] 
\draw  [dash pattern={on 4.5pt off 4.5pt}]  (262,1387) -- (262.5,1586) ;
%Straight Lines [id:da584080201029381] 
\draw  [dash pattern={on 4.5pt off 4.5pt}]  (281,1386) -- (281.5,1569) ;
%Straight Lines [id:da4364578065065364] 
\draw  [dash pattern={on 4.5pt off 4.5pt}]  (322,1386) -- (321.37,1558.42) ;
%Straight Lines [id:da258639231602919] 
\draw  [dash pattern={on 4.5pt off 4.5pt}]  (345,1386) -- (345.37,1538.42) ;
%Straight Lines [id:da34253816558936423] 
\draw  [dash pattern={on 4.5pt off 4.5pt}]  (411.5,1389) -- (411.5,1514) ;
%Straight Lines [id:da2536648998995712] 
\draw  [dash pattern={on 4.5pt off 4.5pt}]  (432.5,1387) -- (432.5,1494) ;
%Straight Lines [id:da8590792784683857] 
\draw  [dash pattern={on 4.5pt off 4.5pt}]  (491.5,1387) -- (491.5,1477) ;
%Straight Lines [id:da057989503617009164] 
\draw  [dash pattern={on 4.5pt off 4.5pt}]  (531.5,1386) -- (531,1462) ;
%Shape: Ellipse [id:dp7782871036529543] 
\path  [shading=_0a952aycu,_dlyrgacgu,path fading= _rdz2r7nl8 ,fading transform={xshift=2}] (100,1388.13) .. controls (100,1382.19) and (104.82,1377.37) .. (110.75,1377.37) .. controls (116.69,1377.37) and (121.51,1382.19) .. (121.51,1388.13) .. controls (121.51,1394.07) and (116.69,1398.88) .. (110.75,1398.88) .. controls (104.82,1398.88) and (100,1394.07) .. (100,1388.13) -- cycle ; % for fading 
 \draw   (100,1388.13) .. controls (100,1382.19) and (104.82,1377.37) .. (110.75,1377.37) .. controls (116.69,1377.37) and (121.51,1382.19) .. (121.51,1388.13) .. controls (121.51,1394.07) and (116.69,1398.88) .. (110.75,1398.88) .. controls (104.82,1398.88) and (100,1394.07) .. (100,1388.13) -- cycle ; % for border 

%Shape: Ellipse [id:dp6919696246852687] 
\path  [shading=_an1dxvxiw,_u16zh8rhg,path fading= _3jqkjquc1 ,fading transform={xshift=2}] (161.11,1388.13) .. controls (161.11,1382.19) and (165.92,1377.37) .. (171.86,1377.37) .. controls (177.8,1377.37) and (182.62,1382.19) .. (182.62,1388.13) .. controls (182.62,1394.07) and (177.8,1398.88) .. (171.86,1398.88) .. controls (165.92,1398.88) and (161.11,1394.07) .. (161.11,1388.13) -- cycle ; % for fading 
 \draw   (161.11,1388.13) .. controls (161.11,1382.19) and (165.92,1377.37) .. (171.86,1377.37) .. controls (177.8,1377.37) and (182.62,1382.19) .. (182.62,1388.13) .. controls (182.62,1394.07) and (177.8,1398.88) .. (171.86,1398.88) .. controls (165.92,1398.88) and (161.11,1394.07) .. (161.11,1388.13) -- cycle ; % for border 

%Shape: Ellipse [id:dp5970600493679222] 
\path  [shading=_1lci1pl48,_330i92zur,path fading= _qfjfhwqkz ,fading transform={xshift=2}] (228.08,1388.13) .. controls (228.08,1382.19) and (232.89,1377.37) .. (238.83,1377.37) .. controls (244.77,1377.37) and (249.59,1382.19) .. (249.59,1388.13) .. controls (249.59,1394.07) and (244.77,1398.88) .. (238.83,1398.88) .. controls (232.89,1398.88) and (228.08,1394.07) .. (228.08,1388.13) -- cycle ; % for fading 
 \draw   (228.08,1388.13) .. controls (228.08,1382.19) and (232.89,1377.37) .. (238.83,1377.37) .. controls (244.77,1377.37) and (249.59,1382.19) .. (249.59,1388.13) .. controls (249.59,1394.07) and (244.77,1398.88) .. (238.83,1398.88) .. controls (232.89,1398.88) and (228.08,1394.07) .. (228.08,1388.13) -- cycle ; % for border 

%Shape: Ellipse [id:dp6664917359985112] 
\path  [shading=_nmuduffyl,_giqcsqxib,path fading= _2l7obxm54 ,fading transform={xshift=2}] (249.59,1388.13) .. controls (249.59,1382.19) and (254.4,1377.37) .. (260.34,1377.37) .. controls (266.28,1377.37) and (271.1,1382.19) .. (271.1,1388.13) .. controls (271.1,1394.07) and (266.28,1398.88) .. (260.34,1398.88) .. controls (254.4,1398.88) and (249.59,1394.07) .. (249.59,1388.13) -- cycle ; % for fading 
 \draw   (249.59,1388.13) .. controls (249.59,1382.19) and (254.4,1377.37) .. (260.34,1377.37) .. controls (266.28,1377.37) and (271.1,1382.19) .. (271.1,1388.13) .. controls (271.1,1394.07) and (266.28,1398.88) .. (260.34,1398.88) .. controls (254.4,1398.88) and (249.59,1394.07) .. (249.59,1388.13) -- cycle ; % for border 

%Shape: Ellipse [id:dp9181708107873152] 
\path  [shading=_e5h2195tt,_x5xws0sb9,path fading= _zvs913mvf ,fading transform={xshift=2}] (270.14,1388.13) .. controls (270.14,1382.19) and (274.96,1377.37) .. (280.9,1377.37) .. controls (286.84,1377.37) and (291.65,1382.19) .. (291.65,1388.13) .. controls (291.65,1394.07) and (286.84,1398.88) .. (280.9,1398.88) .. controls (274.96,1398.88) and (270.14,1394.07) .. (270.14,1388.13) -- cycle ; % for fading 
 \draw   (270.14,1388.13) .. controls (270.14,1382.19) and (274.96,1377.37) .. (280.9,1377.37) .. controls (286.84,1377.37) and (291.65,1382.19) .. (291.65,1388.13) .. controls (291.65,1394.07) and (286.84,1398.88) .. (280.9,1398.88) .. controls (274.96,1398.88) and (270.14,1394.07) .. (270.14,1388.13) -- cycle ; % for border 

%Shape: Ellipse [id:dp8727053385582145] 
\path  [shading=_9d2mgtb6j,_s4wfgjn8z,path fading= _w29ut9usx ,fading transform={xshift=2}] (310.23,1388.13) .. controls (310.23,1382.19) and (315.04,1377.37) .. (320.98,1377.37) .. controls (326.92,1377.37) and (331.74,1382.19) .. (331.74,1388.13) .. controls (331.74,1394.07) and (326.92,1398.88) .. (320.98,1398.88) .. controls (315.04,1398.88) and (310.23,1394.07) .. (310.23,1388.13) -- cycle ; % for fading 
 \draw   (310.23,1388.13) .. controls (310.23,1382.19) and (315.04,1377.37) .. (320.98,1377.37) .. controls (326.92,1377.37) and (331.74,1382.19) .. (331.74,1388.13) .. controls (331.74,1394.07) and (326.92,1398.88) .. (320.98,1398.88) .. controls (315.04,1398.88) and (310.23,1394.07) .. (310.23,1388.13) -- cycle ; % for border 

%Shape: Ellipse [id:dp8532596321725717] 
\path  [shading=_65q9ps644,_n5s869br7,path fading= _a01hbpe3c ,fading transform={xshift=2}] (331.74,1388.08) .. controls (331.76,1382.14) and (336.6,1377.35) .. (342.54,1377.37) .. controls (348.48,1377.4) and (353.27,1382.23) .. (353.25,1388.17) .. controls (353.22,1394.11) and (348.38,1398.91) .. (342.44,1398.88) .. controls (336.5,1398.86) and (331.71,1394.02) .. (331.74,1388.08) -- cycle ; % for fading 
 \draw   (331.74,1388.08) .. controls (331.76,1382.14) and (336.6,1377.35) .. (342.54,1377.37) .. controls (348.48,1377.4) and (353.27,1382.23) .. (353.25,1388.17) .. controls (353.22,1394.11) and (348.38,1398.91) .. (342.44,1398.88) .. controls (336.5,1398.86) and (331.71,1394.02) .. (331.74,1388.08) -- cycle ; % for border 

%Shape: Ellipse [id:dp8689445377260001] 
\path  [shading=_urlhw29dz,_cd0nmm8wa,path fading= _3abpagegl ,fading transform={xshift=2}] (400.24,1388.13) .. controls (400.24,1382.19) and (405.05,1377.37) .. (410.99,1377.37) .. controls (416.93,1377.37) and (421.75,1382.19) .. (421.75,1388.13) .. controls (421.75,1394.07) and (416.93,1398.88) .. (410.99,1398.88) .. controls (405.05,1398.88) and (400.24,1394.07) .. (400.24,1388.13) -- cycle ; % for fading 
 \draw   (400.24,1388.13) .. controls (400.24,1382.19) and (405.05,1377.37) .. (410.99,1377.37) .. controls (416.93,1377.37) and (421.75,1382.19) .. (421.75,1388.13) .. controls (421.75,1394.07) and (416.93,1398.88) .. (410.99,1398.88) .. controls (405.05,1398.88) and (400.24,1394.07) .. (400.24,1388.13) -- cycle ; % for border 

%Shape: Ellipse [id:dp2538743240137955] 
\path  [shading=_drvzrahfu,_kkjo516ck,path fading= _w5az71zv7 ,fading transform={xshift=2}] (421.75,1388.13) .. controls (421.75,1382.19) and (426.56,1377.37) .. (432.5,1377.37) .. controls (438.44,1377.37) and (443.25,1382.19) .. (443.25,1388.13) .. controls (443.25,1394.07) and (438.44,1398.88) .. (432.5,1398.88) .. controls (426.56,1398.88) and (421.75,1394.07) .. (421.75,1388.13) -- cycle ; % for fading 
 \draw   (421.75,1388.13) .. controls (421.75,1382.19) and (426.56,1377.37) .. (432.5,1377.37) .. controls (438.44,1377.37) and (443.25,1382.19) .. (443.25,1388.13) .. controls (443.25,1394.07) and (438.44,1398.88) .. (432.5,1398.88) .. controls (426.56,1398.88) and (421.75,1394.07) .. (421.75,1388.13) -- cycle ; % for border 

%Shape: Ellipse [id:dp18118476010566376] 
\path  [shading=_907n669nz,_zcqapn3su,path fading= _t43uh97cw ,fading transform={xshift=2}] (480.75,1388.13) .. controls (480.75,1382.19) and (485.56,1377.37) .. (491.5,1377.37) .. controls (497.44,1377.37) and (502.25,1382.19) .. (502.25,1388.13) .. controls (502.25,1394.07) and (497.44,1398.88) .. (491.5,1398.88) .. controls (485.56,1398.88) and (480.75,1394.07) .. (480.75,1388.13) -- cycle ; % for fading 
 \draw   (480.75,1388.13) .. controls (480.75,1382.19) and (485.56,1377.37) .. (491.5,1377.37) .. controls (497.44,1377.37) and (502.25,1382.19) .. (502.25,1388.13) .. controls (502.25,1394.07) and (497.44,1398.88) .. (491.5,1398.88) .. controls (485.56,1398.88) and (480.75,1394.07) .. (480.75,1388.13) -- cycle ; % for border 

%Shape: Ellipse [id:dp65438436488949] 
\path  [shading=_itvpsnfvz,_dpr19xkn2,path fading= _qu4qjuiu4 ,fading transform={xshift=2}] (520.75,1388.13) .. controls (520.75,1382.19) and (525.56,1377.37) .. (531.5,1377.37) .. controls (537.44,1377.37) and (542.25,1382.19) .. (542.25,1388.13) .. controls (542.25,1394.07) and (537.44,1398.88) .. (531.5,1398.88) .. controls (525.56,1398.88) and (520.75,1394.07) .. (520.75,1388.13) -- cycle ; % for fading 
 \draw   (520.75,1388.13) .. controls (520.75,1382.19) and (525.56,1377.37) .. (531.5,1377.37) .. controls (537.44,1377.37) and (542.25,1382.19) .. (542.25,1388.13) .. controls (542.25,1394.07) and (537.44,1398.88) .. (531.5,1398.88) .. controls (525.56,1398.88) and (520.75,1394.07) .. (520.75,1388.13) -- cycle ; % for border 

%Curve Lines [id:da5624968402908924] 
\draw    (347.54,1374.55) .. controls (352.92,1360.95) and (362.28,1363.99) .. (370.08,1374.74) ;
\draw [shift={(371.18,1376.31)}, rotate = 236.31] [fill={rgb, 255:red, 0; green, 0; blue, 0 }  ][line width=0.08]  [draw opacity=0] (12,-3) -- (0,0) -- (12,3) -- cycle    ;
%Shape: Square [id:dp49279572065221755] 
\draw  [color={rgb, 255:red, 0; green, 0; blue, 0 }  ,draw opacity=1 ][fill={rgb, 255:red, 245; green, 166; blue, 35 }  ,fill opacity=0.42 ] (353.36,1521.43) -- (373.93,1521.43) -- (373.93,1542) -- (353.36,1542) -- cycle ;
%Straight Lines [id:da748708493315405] 
\draw [color={rgb, 255:red, 208; green, 2; blue, 27 }  ,draw opacity=1 ][line width=3]    (353.43,1542.5) -- (373.93,1542.5) ;
%Straight Lines [id:da4537753219841275] 
\draw [color={rgb, 255:red, 208; green, 2; blue, 27 }  ,draw opacity=1 ][line width=3]    (373.93,1521.43) -- (373.93,1542.5) ;
%Straight Lines [id:da6900455283245166] 
\draw  [dash pattern={on 4.5pt off 4.5pt}]  (568,1542.5) -- (363.68,1542.5) ;
%Straight Lines [id:da7939076531563335] 
\draw  [dash pattern={on 4.5pt off 4.5pt}]  (571,1521.43) -- (373.93,1521.43) ;
%Shape: Square [id:dp3861904780129314] 
\draw  [color={rgb, 255:red, 0; green, 0; blue, 0 }  ,draw opacity=1 ][fill={rgb, 255:red, 80; green, 227; blue, 194 }  ,fill opacity=0.52 ] (421.07,1501.43) -- (441.64,1501.43) -- (441.64,1522) -- (421.07,1522) -- cycle ;
%Shape: Square [id:dp24231016032424302] 
\draw  [color={rgb, 255:red, 0; green, 0; blue, 0 }  ,draw opacity=1 ][fill={rgb, 255:red, 80; green, 227; blue, 194 }  ,fill opacity=0.52 ] (502,1465) -- (520,1465) -- (520,1483) -- (502,1483) -- cycle ;
%Shape: Square [id:dp34367071696478824] 
\draw  [color={rgb, 255:red, 0; green, 0; blue, 0 }  ,draw opacity=1 ][fill={rgb, 255:red, 80; green, 227; blue, 194 }  ,fill opacity=0.52 ] (541,1447) -- (559,1447) -- (559,1465) -- (541,1465) -- cycle ;
%Shape: Square [id:dp7946090479199688] 
\draw  [color={rgb, 255:red, 0; green, 0; blue, 0 }  ,draw opacity=1 ][fill={rgb, 255:red, 80; green, 227; blue, 194 }  ,fill opacity=0.52 ] (291,1560) -- (309,1560) -- (309,1578) -- (291,1578) -- cycle ;
%Shape: Square [id:dp10838546457296283] 
\draw  [color={rgb, 255:red, 0; green, 0; blue, 0 }  ,draw opacity=1 ][fill={rgb, 255:red, 80; green, 227; blue, 194 }  ,fill opacity=0.52 ] (182.36,1616.56) -- (202.66,1616.56) -- (202.66,1636.86) -- (182.36,1636.86) -- cycle ;
%Shape: Square [id:dp7913556603930977] 
\draw  [color={rgb, 255:red, 0; green, 0; blue, 0 }  ,draw opacity=1 ][fill={rgb, 255:red, 80; green, 227; blue, 194 }  ,fill opacity=0.52 ] (121.72,1636.99) -- (142.02,1636.99) -- (142.02,1657.29) -- (121.72,1657.29) -- cycle ;

% Text Node
\draw (41,1625) node [anchor=north west][inner sep=0.75pt]    {$z^{n}(t)$};
% Text Node
\draw (574,1530) node [anchor=north west][inner sep=0.75pt]    {$z_{i}^{n}( t)$};
% Text Node
\draw (574,1505) node [anchor=north west][inner sep=0.75pt]    {$z_{i}^{n}( t-)$};
% Text Node
\draw (363,1389) node [anchor=north west][inner sep=0.75pt]    {$i$};
\end{tikzpicture}
\caption{Graphical representation of the TASEP and evolution of the height function $z^n(t)$ for a fixed $n$. When a particle has space to jump, the height function will be in a configuration that a top-left corner can be flipped to a low-right corner. All these locations are the growth sites, and they are marked with a shaded square below the level of $z^n$. In the diagram, the particle at location $i-1$ will jump at $i$ at time $t$ and the new height at location $k$ will lower by 1.}
\label{fig:zevolution}
\end{figure}
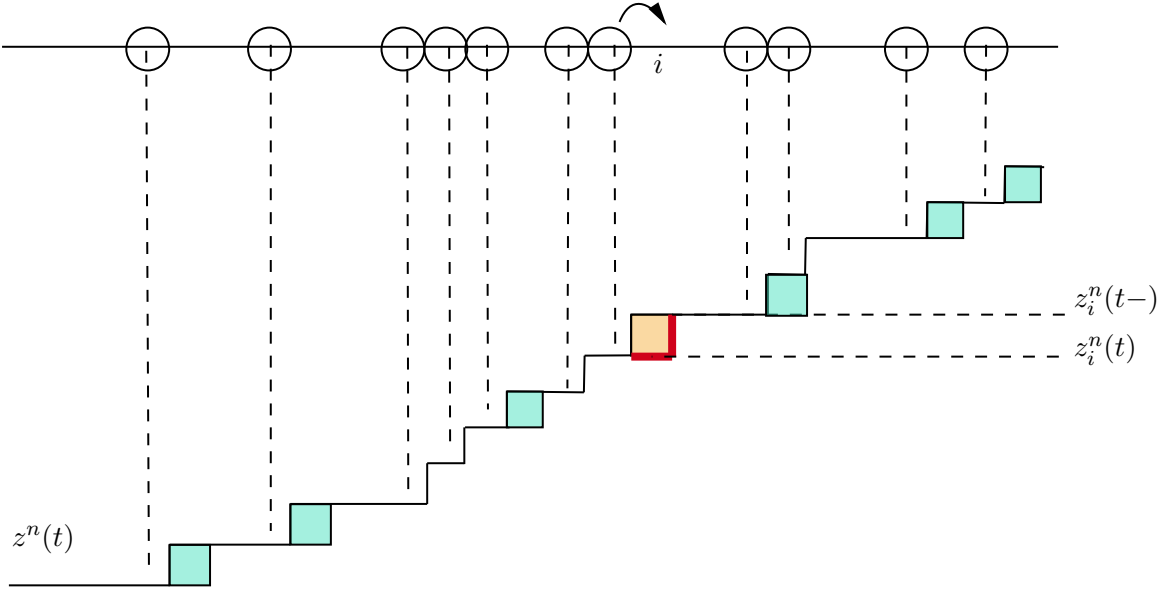
Now let us define the evolution for all times $t>0$.
Order the Poisson event times of the processes  $\{N^{(n)}_{i,j}\}_{(i, j) \in[L, U] \times [L, U]}$ 
up to time $T$, say 
$0 \le t_1 < t_2 < \ldots < t_N \le T$.  Follow through these event times, until a Poisson event is taking place at a position where $z_i$ is allowed to jump, say $i$, and assume this event time is $t_{k}$. There is no change to $z^n$ from time $t_{k-1}$ until $t_{k}^-$ and at time $t_{k}$ we have 
\[
z_{i}^{n}(t_{k}) = z_{i}^{n}(t_{k}^-)-1,
\] 
while all other values remain unchanged. This corresponds to a particle jumping from location $i-1$ to location $i$ in the TASEP at time $t_{k}$,  giving 
\[
\eta_{i-1}^{n}(t_{k}^-)(1 - \eta_{i}^{n}(t_{k}^-)) = 1 = \eta_{i}^{n}(t_{k})(1 - \eta_{i-1}^{n}(t_{k})).
\]
This particle jumped at the jump time of $N^{(n)}_{i, -z_i^{n}(t_{k}^-)}$. Now we can restart the process by the strong Markov property, and keep updating the $z$ process, until the end of our time horizon. 
A decrease in $z^{n}_{i}$ is equivalent to an exclusion particle jump from site $i-1$ to $i$. Dynamically speaking, we say that when $z^{n}_{i}$ is free to move, it is decreased by $1$ at the next event time of $N^{(n)}_{i,j}$, which will be exponentially distributed as $\tau^{(n)}_{i,j}$, with rate $\lambda^n_{i,j}$ and where the index $j = - z^{n}_{i}(t)$.

After this, update the growth sites in the interval we are looking at and restart from this height function using the strong Markov property. Proceed in a similar way until the time horizon is reached. We refer to this construction as the {\bf graphical construction}, as introduced by Harris in \cite{Harris1972}. 

As the TASEP uses rates sampled from the speed function $\tilde c$, it is convenient to rephrase Assumptions \ref{ass:c}, \ref{ass:c2}, \ref{ass:c0} and \ref{assumption: growth rates} for this speed function. 

\begin{figure}[ht]
\tikzset{every picture/.style={line width=0.75pt}} %set default line width to 0.75pt        

\begin{tikzpicture}[x=0.75pt,y=0.75pt,yscale=-1,xscale=1]
%uncomment if require: \path (0,861); %set diagram left start at 0, and has height of 861

%Shape: Arc [id:dp26095699260119143] 
\draw  [draw opacity=0][fill={rgb, 255:red, 184; green, 233; blue, 134 }  ,fill opacity=1 ][dash pattern={on 0.84pt off 2.51pt}] (482.9,258.07) .. controls (471.84,241.99) and (475.33,220) .. (491.01,208.42) .. controls (507.01,196.6) and (529.65,200.1) .. (541.6,216.23) -- (512.63,237.68) -- cycle ; \draw  [dash pattern={on 0.84pt off 2.51pt}] (482.9,258.07) .. controls (471.84,241.99) and (475.33,220) .. (491.01,208.42) .. controls (507.01,196.6) and (529.65,200.1) .. (541.6,216.23) ;  
%Shape: Rectangle [id:dp7096563738950605] 
\draw   (86,57) -- (604,57) -- (604,302) -- (86,302) -- cycle ;
%Curve Lines [id:da6095858575244465] 
\draw    (86,119) .. controls (126,89) and (172,285) .. (261,302) ;
%Shape: Polygon [id:ds6432442040833322] 
\draw   (554,58) -- (604,88) -- (554,168) -- (524,148) -- (504,88) -- cycle ;
%Curve Lines [id:da5938254678987595] 
\draw    (181,233) .. controls (221,203) and (298,88) .. (419,247) ;
%Curve Lines [id:da43203364509960474] 
\draw    (87,81) .. controls (127,51) and (378,49) .. (375,199) ;
%Curve Lines [id:da6597921678791138] 
\draw    (504,88) .. controls (544,58) and (378,49) .. (375,199) ;
%Straight Lines [id:da15394307735165913] 
\draw    (573,194) -- (451,279) ;
%Shape: Circle [id:dp3634000822186473] 
\draw  [color={rgb, 255:red, 189; green, 16; blue, 224 }  ,draw opacity=1 ][fill={rgb, 255:red, 255; green, 255; blue, 255 }  ,fill opacity=0.34 ][line width=3.75]  (90.75,116.37) .. controls (90.75,112.3) and (94.05,109) .. (98.13,109) .. controls (102.2,109) and (105.5,112.3) .. (105.5,116.37) .. controls (105.5,120.45) and (102.2,123.75) .. (98.13,123.75) .. controls (94.05,123.75) and (90.75,120.45) .. (90.75,116.37) -- cycle ;
%Shape: Circle [id:dp8399013615751825] 
\draw  [color={rgb, 255:red, 65; green, 117; blue, 5 }  ,draw opacity=1 ][fill={rgb, 255:red, 255; green, 255; blue, 255 }  ,fill opacity=0.34 ][line width=3.75]  (443.5,279) .. controls (443.5,274.86) and (446.86,271.5) .. (451,271.5) .. controls (455.14,271.5) and (458.5,274.86) .. (458.5,279) .. controls (458.5,283.14) and (455.14,286.5) .. (451,286.5) .. controls (446.86,286.5) and (443.5,283.14) .. (443.5,279) -- cycle ;
%Shape: Circle [id:dp4142788924413092] 
\draw  [color={rgb, 255:red, 65; green, 117; blue, 5 }  ,draw opacity=1 ][fill={rgb, 255:red, 255; green, 255; blue, 255 }  ,fill opacity=0.34 ][line width=3.75]  (411.5,247) .. controls (411.5,242.86) and (414.86,239.5) .. (419,239.5) .. controls (423.14,239.5) and (426.5,242.86) .. (426.5,247) .. controls (426.5,251.14) and (423.14,254.5) .. (419,254.5) .. controls (414.86,254.5) and (411.5,251.14) .. (411.5,247) -- cycle ;
%Shape: Circle [id:dp4705600687903383] 
\draw  [color={rgb, 255:red, 189; green, 16; blue, 224 }  ,draw opacity=1 ][fill={rgb, 255:red, 255; green, 255; blue, 255 }  ,fill opacity=0.34 ][line width=3.75]  (289,166.5) .. controls (289,162.36) and (292.36,159) .. (296.5,159) .. controls (300.64,159) and (304,162.36) .. (304,166.5) .. controls (304,170.64) and (300.64,174) .. (296.5,174) .. controls (292.36,174) and (289,170.64) .. (289,166.5) -- cycle ;
%Shape: Circle [id:dp5380446844859424] 
\draw  [color={rgb, 255:red, 208; green, 2; blue, 27 }  ,draw opacity=1 ][fill={rgb, 255:red, 255; green, 255; blue, 255 }  ,fill opacity=0.34 ][line width=3.75]  (367.5,199) .. controls (367.5,194.86) and (370.86,191.5) .. (375,191.5) .. controls (379.14,191.5) and (382.5,194.86) .. (382.5,199) .. controls (382.5,203.14) and (379.14,206.5) .. (375,206.5) .. controls (370.86,206.5) and (367.5,203.14) .. (367.5,199) -- cycle ;
%Shape: Circle [id:dp22443149398840245] 
\draw  [color={rgb, 255:red, 65; green, 117; blue, 5 }  ,draw opacity=1 ][fill={rgb, 255:red, 255; green, 255; blue, 255 }  ,fill opacity=0.34 ][line width=3.75]  (173.5,233) .. controls (173.5,228.86) and (176.86,225.5) .. (181,225.5) .. controls (185.14,225.5) and (188.5,228.86) .. (188.5,233) .. controls (188.5,237.14) and (185.14,240.5) .. (181,240.5) .. controls (176.86,240.5) and (173.5,237.14) .. (173.5,233) -- cycle ;
%Shape: Circle [id:dp22342977372095196] 
\draw  [color={rgb, 255:red, 189; green, 16; blue, 224 }  ,draw opacity=1 ][fill={rgb, 255:red, 255; green, 255; blue, 255 }  ,fill opacity=0.34 ][line width=3.75]  (173.5,66) .. controls (173.5,61.86) and (176.86,58.5) .. (181,58.5) .. controls (185.14,58.5) and (188.5,61.86) .. (188.5,66) .. controls (188.5,70.14) and (185.14,73.5) .. (181,73.5) .. controls (176.86,73.5) and (173.5,70.14) .. (173.5,66) -- cycle ;
%Shape: Circle [id:dp6514667747071226] 
\draw  [color={rgb, 255:red, 65; green, 117; blue, 5 }  ,draw opacity=1 ][fill={rgb, 255:red, 255; green, 255; blue, 255 }  ,fill opacity=0.34 ][line width=3.75]  (565.5,194) .. controls (565.5,189.86) and (568.86,186.5) .. (573,186.5) .. controls (577.14,186.5) and (580.5,189.86) .. (580.5,194) .. controls (580.5,198.14) and (577.14,201.5) .. (573,201.5) .. controls (568.86,201.5) and (565.5,198.14) .. (565.5,194) -- cycle ;
%Shape: Circle [id:dp6382413706728722] 
\draw  [color={rgb, 255:red, 189; green, 16; blue, 224 }  ,draw opacity=1 ][fill={rgb, 255:red, 255; green, 255; blue, 255 }  ,fill opacity=0.34 ][line width=3.75]  (474.5,69) .. controls (474.5,64.86) and (477.86,61.5) .. (482,61.5) .. controls (486.14,61.5) and (489.5,64.86) .. (489.5,69) .. controls (489.5,73.14) and (486.14,76.5) .. (482,76.5) .. controls (477.86,76.5) and (474.5,73.14) .. (474.5,69) -- cycle ;
%Shape: Circle [id:dp07715578363867481] 
\draw  [color={rgb, 255:red, 189; green, 16; blue, 224 }  ,draw opacity=1 ][fill={rgb, 255:red, 255; green, 255; blue, 255 }  ,fill opacity=0.34 ][line width=3.75]  (546.5,58) .. controls (546.5,53.86) and (549.86,50.5) .. (554,50.5) .. controls (558.14,50.5) and (561.5,53.86) .. (561.5,58) .. controls (561.5,62.14) and (558.14,65.5) .. (554,65.5) .. controls (549.86,65.5) and (546.5,62.14) .. (546.5,58) -- cycle ;
%Shape: Circle [id:dp9140741776177164] 
\draw  [color={rgb, 255:red, 189; green, 16; blue, 224 }  ,draw opacity=1 ][fill={rgb, 255:red, 255; green, 255; blue, 255 }  ,fill opacity=0.34 ][line width=3.75]  (500.96,90.19) .. controls (499.75,88.51) and (500.13,86.17) .. (501.81,84.96) .. controls (503.49,83.75) and (505.83,84.13) .. (507.04,85.81) .. controls (508.25,87.49) and (507.87,89.83) .. (506.19,91.04) .. controls (504.51,92.25) and (502.17,91.87) .. (500.96,90.19) -- cycle ;
%Shape: Circle [id:dp8195578524149979] 
\draw  [color={rgb, 255:red, 189; green, 16; blue, 224 }  ,draw opacity=1 ][fill={rgb, 255:red, 255; green, 255; blue, 255 }  ,fill opacity=0.34 ][line width=3.75]  (520.96,150.19) .. controls (519.75,148.51) and (520.13,146.17) .. (521.81,144.96) .. controls (523.49,143.75) and (525.83,144.13) .. (527.04,145.81) .. controls (528.25,147.49) and (527.87,149.83) .. (526.19,151.04) .. controls (524.51,152.25) and (522.17,151.87) .. (520.96,150.19) -- cycle ;
%Shape: Circle [id:dp9341858855665393] 
\draw  [color={rgb, 255:red, 189; green, 16; blue, 224 }  ,draw opacity=1 ][fill={rgb, 255:red, 255; green, 255; blue, 255 }  ,fill opacity=0.34 ][line width=3.75]  (549.73,171.08) .. controls (548.03,168.72) and (548.56,165.43) .. (550.92,163.73) .. controls (553.28,162.03) and (556.57,162.56) .. (558.27,164.92) .. controls (559.97,167.28) and (559.44,170.57) .. (557.08,172.27) .. controls (554.72,173.97) and (551.43,173.44) .. (549.73,171.08) -- cycle ;
%Shape: Circle [id:dp5292694888334302] 
\draw  [color={rgb, 255:red, 189; green, 16; blue, 224 }  ,draw opacity=1 ][fill={rgb, 255:red, 255; green, 255; blue, 255 }  ,fill opacity=0.34 ][line width=3.75]  (599.73,91.08) .. controls (598.03,88.72) and (598.56,85.43) .. (600.92,83.73) .. controls (603.28,82.03) and (606.57,82.56) .. (608.27,84.92) .. controls (609.97,87.28) and (609.44,90.57) .. (607.08,92.27) .. controls (604.72,93.97) and (601.43,93.44) .. (599.73,91.08) -- cycle ;
%Shape: Circle [id:dp10016280964093249] 
\draw  [color={rgb, 255:red, 208; green, 2; blue, 27 }  ,draw opacity=1 ][fill={rgb, 255:red, 255; green, 255; blue, 255 }  ,fill opacity=0.34 ][line width=3.75]  (507.96,75.19) .. controls (506.75,73.51) and (507.13,71.17) .. (508.81,69.96) .. controls (510.49,68.75) and (512.83,69.13) .. (514.04,70.81) .. controls (515.25,72.49) and (514.87,74.83) .. (513.19,76.04) .. controls (511.51,77.25) and (509.17,76.87) .. (507.96,75.19) -- cycle ;
%Shape: Circle [id:dp05151667135365445] 
\draw  [color={rgb, 255:red, 65; green, 117; blue, 5 }  ,draw opacity=1 ][fill={rgb, 255:red, 255; green, 255; blue, 255 }  ,fill opacity=0.34 ][line width=3.75]  (104.25,321) .. controls (104.25,316.86) and (107.61,313.5) .. (111.75,313.5) .. controls (115.89,313.5) and (119.25,316.86) .. (119.25,321) .. controls (119.25,325.14) and (115.89,328.5) .. (111.75,328.5) .. controls (107.61,328.5) and (104.25,325.14) .. (104.25,321) -- cycle ;
%Shape: Circle [id:dp20814210233312924] 
\draw  [color={rgb, 255:red, 189; green, 16; blue, 224 }  ,draw opacity=1 ][fill={rgb, 255:red, 255; green, 255; blue, 255 }  ,fill opacity=0.34 ][line width=3.75]  (104.25,347.5) .. controls (104.25,343.36) and (107.61,340) .. (111.75,340) .. controls (115.89,340) and (119.25,343.36) .. (119.25,347.5) .. controls (119.25,351.64) and (115.89,355) .. (111.75,355) .. controls (107.61,355) and (104.25,351.64) .. (104.25,347.5) -- cycle ;
%Shape: Circle [id:dp7229368689119451] 
\draw  [color={rgb, 255:red, 208; green, 2; blue, 27 }  ,draw opacity=1 ][fill={rgb, 255:red, 255; green, 255; blue, 255 }  ,fill opacity=0.34 ][line width=3.75]  (104.25,375) .. controls (104.25,370.86) and (107.61,367.5) .. (111.75,367.5) .. controls (115.89,367.5) and (119.25,370.86) .. (119.25,375) .. controls (119.25,379.14) and (115.89,382.5) .. (111.75,382.5) .. controls (107.61,382.5) and (104.25,379.14) .. (104.25,375) -- cycle ;
%Shape: Circle [id:dp3683168351929065] 
\draw  [dash pattern={on 4.5pt off 4.5pt}] (476.75,236.5) .. controls (476.75,217.03) and (492.53,201.25) .. (512,201.25) .. controls (531.47,201.25) and (547.25,217.03) .. (547.25,236.5) .. controls (547.25,255.97) and (531.47,271.75) .. (512,271.75) .. controls (492.53,271.75) and (476.75,255.97) .. (476.75,236.5) -- cycle ;
%Shape: Circle [id:dp3643674658324356] 
\draw  [fill={rgb, 255:red, 245; green, 166; blue, 35 }  ,fill opacity=1 ] (505.75,236.5) .. controls (505.75,233.05) and (508.55,230.25) .. (512,230.25) .. controls (515.45,230.25) and (518.25,233.05) .. (518.25,236.5) .. controls (518.25,239.95) and (515.45,242.75) .. (512,242.75) .. controls (508.55,242.75) and (505.75,239.95) .. (505.75,236.5) -- cycle ;

% Text Node
\draw (120,152) node [anchor=north west][inner sep=0.75pt]    {$h_{1}$};
% Text Node
\draw (227,260) node [anchor=north west][inner sep=0.75pt]    {$h_{2}$};
% Text Node
\draw (218,66) node [anchor=north west][inner sep=0.75pt]    {$h_{4}$};
% Text Node
\draw (132,65) node [anchor=north west][inner sep=0.75pt]    {$h_{3}$};
% Text Node
\draw (138,309) node [anchor=north west][inner sep=0.75pt]   [align=left] {Terminal /intersection points};
% Text Node
\draw (140,336) node [anchor=north west][inner sep=0.75pt]   [align=left] {Points that separate monotone pieces of curves};
% Text Node
\draw (141,363) node [anchor=north west][inner sep=0.75pt]   [align=left] {Points where the derivative of one discontinuity curve becomes infinite};
% Text Node
\draw (530,260.75) node [anchor=north west][inner sep=0.75pt]    {$L_{\varepsilon ,\ \mathrm{\bf w}}$};
% Text Node
\draw (453,185.75) node [anchor=north west][inner sep=0.75pt]    {$U_{\varepsilon ,\ \mathrm{\bf w}}$};
% Text Node
\draw (500,245) node [anchor=north west][inner sep=0.75pt]    {$\mathrm{\bf w}$};
\end{tikzpicture}
\caption{In this figure you can see a schematic of the decomposition of the discontinuity curves allowed by Assumption \ref{ass:pde-speed}. For example, while $h_3$ and $h_4$ can be viewed as one curve $h$ that is a graph, the monotonicity assumption would break since there is a local maximum for $h$. Thus, that local maximum becomes a terminal point, and $h$ is split into two monotone pieces. Straight line discontinuities are allowed, as shown in the top right corner, as long their slope is not $0$ or 
$-1$. \\
In the lower right corner, you can see an interior point $\bf w$ and the two sets $U_{\e, \bf w}$ and $L_{\e, \bf w}$. The shaded upper set tells us that those are the lower values for the function $\tilde c(x,y)$ and as such the value of $\tilde c(x,y)$ on the discontinuity curve is the continuous extension of the value in the interior of $U_{\e, \bf w}$, as needed by Assumption \ref{ass:pde-speed}(A-2) (equiv. \ref{ass:c2} for $c$).}
\label{fig:DiscCurves}
\end{figure}
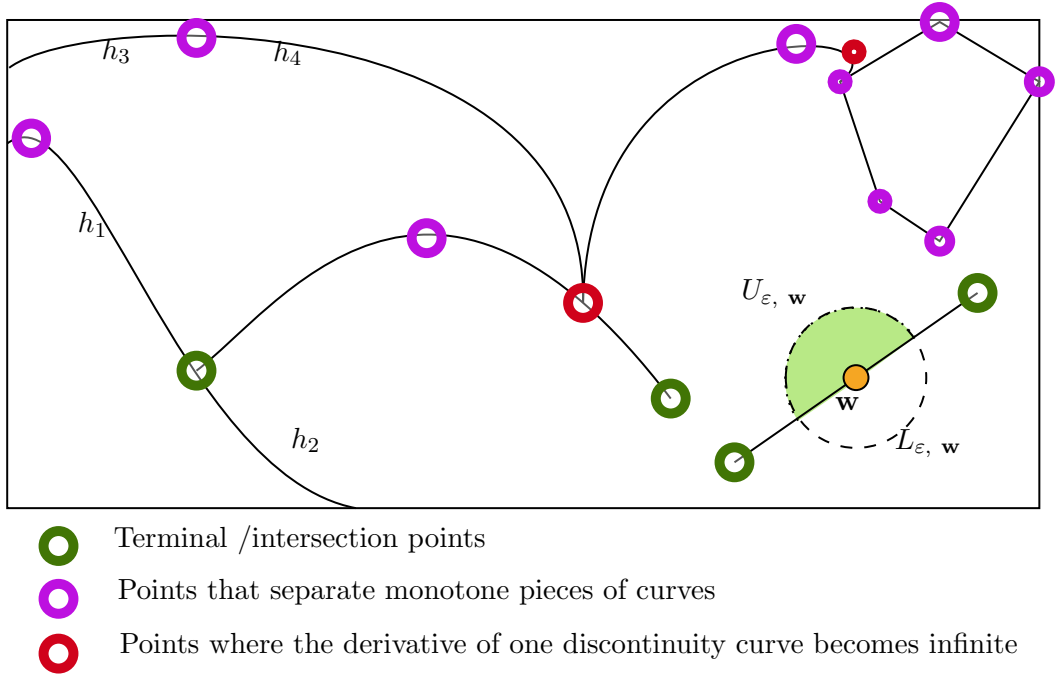

We now record the corresponding assumptions on $\tilde c$. It is straightforward to check that any $\tilde c$ that satisfies Assumption \ref{ass:pde-speed}, then the equivalent $c$ will satisfy \ref{ass:c}, \ref{ass:c2}, \ref{ass:c0} and \ref{assumption: growth rates}. For a graphical representation of these assumptions, see Figure \ref{fig:DiscCurves}.

\begin{assumption}[Assumptions for $\tilde c$]
\label{ass:pde-speed}
We assume the following regularity conditions for the function $\tilde c$:
\begin{enumerate} 
\item[\textnormal{(A-1)}]-Discontinuity curves:
		\begin{enumerate}
			\item There are finitely many discontinuity curves $\{ h_i\}_{i \in I}$ in any compact set $K$. 
			\item All $h_i$ can be viewed as graphs			
			\[
			h_i: [ z_i, w_i] = {\rm Dom}(h_i) \to \R,
			\] 
            or are vertical line segments.
			\item Each discontinuity curve $h_i$ is strictly monotone. In particular this excludes horizontal discontinuity segments. 
            \item $h_i$ is $C^1((z_i, w_i), \R)$. At the boundary points $z_i, w_i$ the derivative may take the value $\pm \infty, 0$.
            \item No discontinuity curves satisfy $h'(s) =-1$ for $s$ in some interval. 
		\end{enumerate} 
\medskip

\item[\textnormal{(A-2)}]-Regularity:
The function $\tilde c$ is continuous away from its discontinuity
curves, locally bounded above and below away from zero, and lower
semicontinuous. At interior points of each discontinuity curve, the value of
$\tilde c$ is obtained as the one-sided continuous extension from the
lower-speed side. At any terminal point, there will be a continuity region of $\tilde c$ adjacent to it with the value of $\tilde c$ at the point to equal the limits of values from that continuity region.
\medskip

\item[\textnormal{(A-3)}]-Growth:
The speed function has sub-logarithmic growth at infinity. 
\end{enumerate}
\end{assumption}

Our first theorem is a law of large numbers for the height function and the particle density at time $t$. 

\begin{theorem}(Hydrodynamic limit)\label{thm:LLN}
 For any $n \in \N$, let $z^n$ evolve using a field of independent Poisson processes $\{N^{(n)}_{i,j}\}_{(i,j) \in \Z^2}$ with rates $\{\lambda^n_{i,j}\}_{(i,j) \in \Z^2}$ given by 
\[
        \lambda^n_{i,j}=\tilde c(i/n,j/n)
        = c((i+j)/n,j/n),
\]
and let $\tilde c$ satisfy Assumption \ref{ass:pde-speed}.
Assume \eqref{eq:etaass}, \eqref{eq:in-cur} and \eqref{eq:z0lim} where the limits are assumed either ${\bf P}-$ a.s.~or in ${\bf P}-$probability. 

Then for any $x \in \bR$ and $t \in \bR_+$ we have the existence of a function $v(x,t)$ such that the limit 
\be \label{eq:vdef}
\lim_{n \to \infty} \frac{1}{n} z^{n}_{\fl{nx}}(nt) = v(x,t) 
\ee
holds in the same mode of convergence as for $t=0$, i.e. ${\bf P}-$ a.s.~or in ${\bf P}-$probability. 

Moreover, $v(x,t )$ satisfies the following properties
\begin{enumerate}
\item $v(x,t)$ is Lipschitz-1 in the $x$-variable.
\item  $v(x,t)$ is locally Lipschitz, jointly in $(x,t)$. 
\item $v(x,t)$ is differentiable in both $x$ and $t$ for a.e.\ $(x,t)$. 
\end{enumerate}
The a.e. $x$-derivative $\rho(x,t)$ is the density profile for the particle system, as given by the LLN
\be\label{eq:etaasst} 
\lim_{n\to \infty} \frac{1}{n}\sum_{i = \fl{na}}^{\fl{nb}} \eta^n_i(nt) = \int_{a}^b \rho(x,t) \, dx.
\ee
\end{theorem}

\subsection{Discontinuous Hamilton--Jacobi equations and scalar conservation laws}

In the classical setting of Rost \cite{rost1981}, the homogeneous TASEP limiting density $\rho(x,t)$ is the unique entropy solution to the scalar conservation law 
\[
\begin{cases}
\rho_t + (\rho(1-\rho))_x =0 \\
\medskip

\rho(x,0) = \rho_0(x),
\end{cases}
\]
while the limiting height function $v(x,t)$ can be viewed as the current of the TASEP process and is the unique viscosity solution of the Hamilton--Jacobi equation
\[
\begin{cases}
v_t + v_x(1-v_x) =0 \\
\medskip

v(x,0) = v_0(x).
\end{cases}
\] 

It is natural to expect that the LLNs $v(x,t), \rho(x,t)$ of Theorem \ref{thm:LLN} are also solutions of some appropriate PDEs. We will show that irrespective of the space-time discontinuities in the Poisson rates, the limits $v(x,t)$ and $\rho(x,t)$ from Theorem \ref{thm:LLN} are strong solutions to appropriate Hamilton--Jacobi or scalar conservation laws, but with discontinuous coefficients and with a more elaborate flux function.

To bridge the LLN with a PDE theory, we will need to provide a Lax-Oleinik variational characterization of the current $v(x,t)$.

We begin from the homogeneous TASEP flux function 
\begin{equation} \label{little f defn}
    f(\rho) = 
\begin{cases}
\rho(1- \rho), & 0 \le \rho \le 1\\
-\infty, & \text{otherwise},
\end{cases}
\end{equation}
and its negative Legendre dual 
\begin{equation} \label{Legendre dual relation}
   \psi(y) = - f^*(y) = - \inf_{\rho \in \R} \{ y \rho - f(\rho) \} = 
\begin{cases}
-y, & y < -1,\\
 \frac{1}{4}(1-y)^2, & y \in [-1,1],\\ 
0, & y >1. 
\end{cases} 
\end{equation}
This has a relationship with the last passage shape function $\gamma(x,y)$ from \eqref{eq:limshap}
\begin{equation} \label{g condition}
    \gamma(x+\psi(x),\psi(x)) = 1, \quad \text{for } x \in [-1,1].
\end{equation}
In other words, the curve $\ell(x)= (x+\psi(x),\psi(x)), x \in [-1,1]$ is a level curve of the limiting function $\gamma(x,y)$.

Using the Legendre dual $ \psi$ we can obtain a variational characterisation of $v(x,t)$ from Theorem \ref{thm:LLN}. We define, for a given initial condition $v_0(x)$, 
\begin{align}
\mathcal{H}_{x,t}^{v_0} &= \Big\{ {\bf w}(\cdot) = (w_{1}(\cdot), w_{2}(\cdot)): [0,t] \to \bR^2: \notag \\
& \phantom{xxxxxxxx}{\bf w}(\cdot)\text{ continuous piecewise $C^{1}$ paths }, \notag \\
&\phantom{xxxxxxxx} w_{1}(t) = x, w_{2}(0) = -v_0(w_1(0)), \notag \\
&\phantom{xxxxxxxx}w_{2}'(s) = \tilde c(w_1(s), w_2(s))\psi\left(\frac{w_1'(s)}{\tilde c(w_1(s), w_2(s))}\right) \text{ for a.e. } s \in[0,t]\Big\}. \label{eq:pathsw}
\end{align}
It will be verified during the proof of the next proposition that this set is non-empty; in fact the path set contains re-parametrised versions of the paths in $\mathcal H(x,y)$ from  \eqref{eq:pathsH}.
 
\begin{proposition}[Lax-Oleinik formula] \label{thm:lang}
Let $\tilde c (x,y)$ be given by \eqref{eq:tildec}. Then, the law of large numbers limit \eqref{eq:vdef} $v(x,t)$ can be written as follows:
\begin{equation} \label{v definition}
     v(x,t) = \sup_{{\bf w}(\cdot) \in \mathcal{H}^{v_0}_{x,t}} \left\{v_{0}(w_{1}(0)) - \int^{t}_{0}\tilde c({\bf w}(s))\psi\left( \frac{w_{1}'(s)}{\tilde c({\bf w}(s))} \right) ds \right\}.
\end{equation}
\end{proposition}

Note that since we have the derivative of $w_2$, we also have that for all admissible paths,
\be \label{eq:w2value}
w_2(t) = w_2(0) + \int_{0}^t\tilde c({\bf w}(s))\psi\left( \frac{w_{1}'(s)}{\tilde c({\bf w}(s))} \right) ds 
\ee
so equation 
\eqref{v definition} can be written as
\be \label{eq:vasw2}
 v(x,t) = -\inf_{{\bf w}(\cdot) \in \mathcal{H}^{v_0}_{x,t}} \left\{ w_2(t) \right \}.
\ee

An appropriate initial value problem for the current $v(x,t)$ in the discontinuous TASEP setting would correspond to
\begin{equation} \label{eq:genHJ}
    \begin{cases}
   % v_t(x,t) + H(\tilde c(x,-v(x,t)), v_x(x,t)) = 0, \hspace{15.pt} &(x,t) \in \bR \times (0,T], \\
     v_t + H( x, v, v_x) = 0, \hspace{15.pt} &(x,t) \in \bR \times (0,T], \\
    v(x,0) = v_0(x), \hspace{15.pt} &(x,t) \in \bR \times \{0\}.
    \end{cases}
\end{equation}
The function $H$ here is the Hamiltonian for this particle system. For the TASEP, it is given by
\begin{equation} \label{eq:discoflux}
    H(x, v, v_x) = \begin{cases}
        \tilde c(x,-v(x,t)) v_x(x,t)(1-v_x(x,t)), \hspace{15.pt} &v_x(x,t) \in [0,1] \\
        -\infty, \hspace{15.pt} &v_x(x,t) \notin [0,1].
    \end{cases}
\end{equation}

Note that in our case $\rho(x,t) = v_x(x,t) \in [0,1]$ so \eqref{eq:discoflux} is always well-defined. The remaining results discuss notions of solutions for the initial value problem \eqref{eq:genHJ} with Hamiltonian \ref{eq:discoflux}.

The main task is to determine what type of solution $v(x,t)$ is. Depending on what regularity properties one needs on $v(x,t)$, different regularity properties must be imposed on the speed function $\tilde c(x,y)$, just due to the analytical difficulties afforded by the discontinuity curves. For stronger statements we will need to substitute Assumption \ref{ass:pde-speed}-A-1(e), with the much stronger  

\begin{assumption}\label{ass:c0strong}

Let the discontinuity curves $\{\tilde h_i\}_{i \in \N}$ of $\tilde c$ satisfy 
    \[
    \tilde h_j'(s)\notin[-1,0]
    \qquad\text{for a.e. }s\in I_j.
    \]
\end{assumption}

\begin{remark} The equivalent assumption for the function $c(x,y)$ is that there are no decreasing discontinuity curves. \qed
\end{remark}

The theorems below change flavour depending on the allowable discontinuous $\tilde c$. Stronger control on the discontinuity curves implies stronger results. 

First, $v(x,t)$ will be a strong solution to \eqref{eq:genHJ} at all points of continuity of $\tilde c(x, -v(x,t))$ and differentiability of $v$. For this Assumption \ref{ass:pde-speed} suffices.
These continuity points however, depending on the structure of the discontinuity curves or the initial conditions, may not be dense in $\R^2$, so Assumption \ref{ass:c0strong} needs to be imposed to guarantee that, giving an a.e.~classical solution to the PDE. 

Under Assumption \ref{ass:pde-speed}, $v(x,t)$ will be a viscosity solution for the discontinuous Hamilton--Jacobi equation in an appropriate sense. To obtain uniqueness of the viscosity solution and uniqueness $v_x= \rho$ as a solution to the scalar conservation law we further need to eliminate temporal discontinuities from $\tilde c$ and assume an integrated current entropy condition.

\begin{theorem}[Strong solution to Hamilton--Jacobi equation] \label{thm: existence} 
Assume \ref{ass:pde-speed}. The limiting current $v(x,t)$, given by equation \eqref{eq:vdef}, is a solution to the Cauchy problem \eqref{eq:genHJ} with flux given by \eqref{eq:discoflux}, at all points $(x,t)$ that are differentiability points of $v(x,t)$ and continuity points of $\tilde c(x,-v(x,t))$. 

Furthermore, under Assumption \ref{ass:c0strong}
almost every $(x,t) \in \bR \times \bR_+$ is a differentiability point of $v(x,t)$ and continuity point of $\tilde c(x,-v(x,t))$.
\end{theorem}

It turns out that Assumption \ref{ass:c0strong} is unnecessarily strong if one only cares about a special class of weak solutions to \eqref{eq:no-t-HJ} known as viscosity solutions. 
There are various notions of viscosity in the literature of discontinuous Hamilton--Jacobi equations. For the purposes of this article, the following definition records the viscosity notion selected by the
variational formula when no extra assumptions are imposed on $\tilde c$. It differs from the standard Ishii discontinuous-viscosity
definition, as the upper and lower envelopes are switched.

Let
\[
\tilde c_*({\bf z}):=\liminf_{{\bf z}'\to {\bf z}}\tilde c({\bf z}'),
\qquad
\tilde c^*({\bf z}):=\limsup_{{\bf z}'\to {\bf z}}\tilde c({\bf z}')
\]
and define, for $0\le p\le 1$,
\[
H_{\rm low}(x,r,p)=\tilde c_*(x,-r)p(1-p),
\qquad
H_{\rm up}(x,r,p)=\tilde c^*(x,-r)p(1-p).
\]
Since $\tilde c$ is lower semicontinuous, $\tilde c_*=\widetilde c$.

Then an {\bf envelope selected viscosity solution} to  \eqref{eq:genHJ} is defined as 
\begin{definition}[Envelop-selected viscosity solution to \eqref{eq:genHJ}]
%\label{def:selected-viscosity}
\label{defn: classical viscosity solution}
\rm 
Let $v:\mathbb R\times[0,\infty)\to\mathbb R$ be continuous. We say that
$v$ is an envelope-selected viscosity solution of
\[
        v_t+H_{\rm low}(x,v,v_x)=0
\]
if $v(x,0)=v_0(x)$, and for every $\phi\in C^1$:

\begin{enumerate}
\item(\rm Selected viscosity subsolution) if $v-\phi$ has a local maximum at $(x_0,t_0)$, with $t_0>0$, then
\[
        \phi_t(x_0,t_0)
        +
       H_{\rm low}(x_0,v(x_0,t_0),\phi_x(x_0,t_0))
        \le 0;
\]
\item(\rm Selected viscosity supersolution) if $v-\phi$ has a local minimum at $(x_0,t_0)$, with $t_0>0$, then
\[
        \phi_t(x_0,t_0)
        +
       H_{\rm up}(x_0,v(x_0,t_0),\phi_x(x_0,t_0))
        \ge 0.
\]
\end{enumerate}
\end{definition}

\begin{remark} 
In the usual Ishii discontinuous-viscosity definition, the
subsolution inequality is tested with the upper semicontinuous envelope of the
Hamiltonian, while the supersolution inequality is tested with the lower
semicontinuous envelope. Definition \ref{defn: classical viscosity solution} uses the opposite convention: the
subsolution inequality is written with $H_{\rm low}$, while the supersolution
inequality is written with $H_{\rm up}$.

This convention is the one naturally produced by the Lax--Oleinik formula in
the present model. More precisely, the variational argument yields one-sided
viscosity inequalities controlled by the lower and upper envelopes of the
discontinuous coefficient $\tilde c$. We therefore use the term
\emph{envelope-selected} to emphasize that the solution concept is tied to the
microscopic variational representation, rather than to the standard Ishii
choice of envelopes. 

Indeed, under extra assumptions on the structure of the discontinuity curves and on $\tilde c$ (such as monotonicity), we expect, based on our proofs, that this notion of viscosity can be strengthened to more classical ones. There are two competing forces: The first is that the near-optimal paths of \eqref{v definition} want to stay in regions were $\tilde c$ is small. The second is that the values of $\tilde c$ may not allow for an admissible path approximations from the smaller side, while the value on the discontinuity does. controlling these two with extra assumptions will allow for more strict definitions of viscosity solutions. 

The idea to relax the viscosity definition is analogous in spirit to
the flux-limited viscosity framework \cite{ImbertMonneau2017}, although our Hamiltonian and the space--height dependence of the coefficient place the model outside the standard network setting. \qed
\end{remark}

\begin{theorem}\label{thm:visco} Under Assumption \ref{ass:pde-speed},
the limiting current $v = v(x,t)$ from \eqref{eq:vdef} is an envelope selected viscosity solution to \eqref{eq:genHJ}, in the sense of Definition \ref{defn: classical viscosity solution}. 
\end{theorem}

\begin{remark}This theorem guarantees the existence of an a.e.~ differentiable function $v$ that solves the PDE in the selected viscosity sense, but we do not have uniqueness under the full general conditions on the discontinuities of function $\tilde c
(x,y)$. However, notions of uniqueness of weak solutions for discontinuous PDEs similar to \eqref{eq:genHJ} have been well studied when $\tilde c$ is an actual coefficient and has no dependence on $v$ (i.e. in our case this happens when $\tilde c(x,y) = \tilde c(x)$ only) and under various regularity conditions of $H$ and $\tilde c$. \qed
% (see for example \cite{ChenHu2008, adimurthiXXX, georgiou2010}).
\end{remark}

 Indeed, the situation where we have time-homogeneity and only spatial discontinuities
 \be\label{eq:uni-cond}
 \tilde c(x,y) = \tilde c(x)
 \ee
 is more amenable for uniqueness analysis. For the remaining results, we will assume that $\tilde c$ is lower semicontinuous and satisfies  \eqref{eq:uni-cond}, which gives a speed function $\tilde c(x)$ that needs to satisfy Assumption \ref{ass:pde-speed}-A-1. We denote by  
\be \label{eq:no-t-dependance}
G(x, \rho) = \tilde c(x)f(\rho), 
\ee
where $f$ is given by \eqref{little f defn}.
We know that for every fixed $(x,t)$, $G(x, \rho)$ is concave on the interval $[0,1]$.
The resulting Hamilton--Jacobi equation (which is a special case of \eqref{eq:genHJ}) is 
\be \label{eq:no-t-HJ}
\begin{cases}
v_t + G(x, v_x) =0 \\
\medskip

v(x,0) = v_0(x).
\end{cases}
\ee

When we have \eqref{eq:no-t-HJ}, then by formally differentiating the PDE
 with respect to $x$, the scalar conservation law for $\rho$ is given by the following initial value problem
\begin{equation} \label{SCL}
    \begin{cases}
    \rho_t + (\tilde c(x)\rho(1-\rho))_x = 0, \hspace{15.pt} &(x,t) \in \bR \times (0,T], \\
    \rho(x,0) = \rho_0(x), \hspace{15.pt} &(x,t) \in \bR \times \{0\}.
    \end{cases}
\end{equation}

\begin{remark}
    For the equation to make weak sense, the only assumption needed is that $\tilde c(x)\rho(x,t)(1-\rho(x,t)) \in L^1_{\rm loc}(\R)$ which follows since $\tilde c(x) \in L^\infty_{\rm loc}(\R)$ by Assumption \ref{ass:c2}. \qed
\end{remark} 

While Theorem \ref{thm: existence} guarantees the existence of a strong solution, the main difficulty comes when trying to determine uniqueness in a good class of functions. For this, it is often necessary to impose additional  entropy criteria on the solutions to recover uniqueness.

Aside from uniqueness of $v(x,t)$ in some nice class of functions, we are also interested in showing that our density $\rho(x,t)$ is a unique weak solution to \eqref{SCL}

A locally bounded measurable function $\lambda(x,t)$ is a \textit{weak solution} to \eqref{SCL} provided that for every compactly supported and continuously differentiable test function $\phi$ on $\bR \times [0,\infty)$,
\be \label{weak soln integral form}
\int^{\infty}_0 \int_{\bR} \left(\lambda(x,t) \phi_t(x,t) + G(x, \lambda(x,t)) \phi_x(x,t)\right) dx dt + \int_{\bR} \rho_0(x) \phi(x,0) dx = 0.
\ee

In the spatially inhomogeneous case, since $v$ is locally Lipschitz, 
$v_x=\rho\in[0,1]$ a.e., and $v_t+G(x,v_x)=0$ a.e., differentiating this identity 
in the sense of distributions shows that $\rho=v_x$ is a weak solution of \eqref{SCL}
with initial datum $\rho_0=v_0'$.

\begin{theorem}[Uniqueness criteria for viscosity and entropy solutions]
\label{thm:uniqueness}

Assume that the speed function $\tilde c(x,y)$ satisfies Assumption \ref{ass:pde-speed} and
is dependent on $x$ only, that is to say that $\tilde c(x,y) = \tilde c(x)$. 

(HJ) Then, the limiting current $v = v(x,t)$ is the unique solution in the class of non-decreasing Lipschitz-1 functions in the space variable, in the sense of Definition \ref{defn: classical viscosity solution}.

    (SCL) Moreover, suppose that $\lambda(x,t)$ is a non-negative, locally bounded measurable function on $\bR \times \bR_+$ lying inside the interval $[0,1]$, that is also a weak solution to \eqref{SCL} in the sense of \eqref{weak soln integral form}. Let $\rho(x,t) = v_x(x,t)$ where $v$ solves \eqref{eq:no-t-HJ}. 
    
    Then, for any fixed $t > 0$, we have that for almost every $x \in \bR$,
    \be \label{maximum principle}
    \int^t_0 G(x, \lambda(x,s)) ds \leq \int^t_0 G(x, \rho(x,s)) ds.
    \ee
    Furthermore, if equality in \eqref{maximum principle} holds for almost all $(x,t)$, it implies that $\lambda(x,t) = \rho(x,t)$ for almost all $(x,t)$. 
    \end{theorem}

\begin{remark}
As a consequence of the above theorem, the TASEP density $\rho$ can be thought of as the unique weak solution to \eqref{SCL} that maximises current over time and its space integral is the unique viscosity solution in the class of non-decreasing Lipschitz-1 functions in the space variables. \qed
\end{remark}

\section{Law of large numbers for the height functions} 
\label{sec:Law of large numbers}

In this section, we prove a law of large numbers for the processes $z^{n}$. The method is a technical generalisation of Sepp{\"a}l{\"a}inen's coupling (see e.g. \cite{Seppalainen1999KExclusion}, \cite{Seppalainen2001SlowBond}) using the graphical construction of TASEP that we defined in Section \ref{sec:results}.  

From Theorem \ref{thm:1}, we have the LLN in equation \eqref{eq:6}, but the implicit assumption is that when we compute last passage times, we always start at the origin. Here we need the limit as a function of the starting point. To this end, define 
$$ \Gamma_c((0,0), (x,y)) : = \Gamma_c(x,y)$$
and similarly define $\Gamma_c$ 
	from any starting point $(a,b)$ to any terminal point $(x,y)$, $(x,y) \ge (a,b)$ by 
\begin{equation}\label{macroLPT2}
\Gamma_c((a,b),(x,y))=\sup_{\mathbf{x}(\cdot)\in\mathcal{H}((a,b), (x,y))}\bigg\{\int_0^1\frac{\gamma(\mathbf{x}'(s))}{c(x_1(s),x_2(s))}ds\bigg\},
\end{equation}
where 
\begin{align*}
\mathcal{H}((a,b),(x,y)) &=\{{\bf x}\in C([0,1],\R^2_+):\mathbf{x} \text{ is piecewise }C^1,\mathbf{x}(0)=(a,b), \mathbf{x}(1)=(x,y), \\
& \phantom{xxxxxxxxxxxxxxxxxx}\mathbf{x}'(s)\in \R^2_+ \text{ wherever the derivative is defined}\}.
\end{align*}

For any vector $(a,b) \in \R^2$, let the shift operator $T_{a,b}$ acting on the speed function $c(x,y)$ be defined by 
\[
(T_{a,b}c)(x,y) = c(x+a, y+b).
\]
The following proposition extends Theorem 2.6 of \cite{ciech2021}, to converging sequences to starting and ending points. Moreover, it shows  that changing the starting point corresponds to shifting the speed function and starting from $0$.

We record the generalised LLN for the space-time discontinuous last passage model that we will need to use. This is not included in \cite{ciech2021}. Their proof is postponed to Appendix \ref{app:LLN}.

\begin{proposition} \label{lem:SPEP}
Let $c(x,y)$ satisfy Assumptions \ref{ass:c}, \ref{ass:c2}, \ref{ass:c0} and consider the discretisation 
$c^n(i,j) = c(in^{-1}, jn^{-1})$. Consider further (not necessarily deterministic) sequences of vectors $(a_n, b_n)$ and $(s_n, t_n)$ which satisfy 
\[
n^{-1}(a_n, b_n) \to (a,b) \text{ and }  n^{-1}(s_n, t_n) \to (s,t),
\] 
with $(s_n,t_n)-(a_n, b_n) \in \R^2_+$.
Then 
\[
\lim_{n\to \infty} \frac{1}{n}G^{(n)}_{(a_n, b_n), (s_n, t_n)} = \Gamma_c((a,b),(s,t)) = \Gamma_{T_{a,b}c}(s-a,t-b).
\] 
\end{proposition}

\subsection{Step initial conditions for the height process and a coupling to the last passage time model.} 

In this subsection we define a family of auxiliary height functions. These functions are subject to the same evolution rules as the original height function $-z^n$ (see \eqref{Server process constraint} and \eqref{eq:parev}) and their evolutions will be coupled through a common realisation of the Poisson processes $\{N_{i,j}\}_{(i,j) \in \Z^2}$. Moreover, each one of these auxiliary processes will be coupled to a last passage time process.  

For any $k \in \Z$ define the shift vector
\[
{\bf k} = (k, -z^n_k(0)).
\]
Each member of the family of height functions $\xi^{n,\bf k} = \{ \xi^{n,\bf k}_{i}\}_{ i \in \mathbb{Z}}$ corresponds to a distinct exclusion process, as given by equation \eqref{eq:parev}. 

For each $\bf k$,  $n \in \mathbb{N}$ let the function $\xi^{n,\bf k}$ evolve {\em upward} from the initial condition

\begin{equation} \label{aux height process initial conditions}
    \xi^{n,\bf k}_{i}(0) = 
        \begin{cases}
            0, \hspace{27.pt} i \geq 0, \\
            -i, \hspace{20.pt} i < 0.
        \end{cases}
\end{equation}
The index $k$ indicates which Poisson processes are used for the evolution of $ \xi^{n,\bf k}$, given a background height function $z^n$.
 Indeed, $ \xi^{n,\bf k}_{i}$ attempts to increase by 1 at the event times of $N^{(n)}_{i+k,j'}(t)$ with $j'= \xi^{n,\bf k}_{i}(t) - z^{n}_{k}(0)$,
  and it succeeds as long as this does not violate the constraints
 \begin{equation} \label{auxiliary height process constraints}
    \xi^{n,\bf k}_{i}(t) \leq \xi^{n,\bf k}_{i-1}(t) \hspace{10pt} \text{and} \hspace{10pt} \xi^{n,\bf k}_{i}(t) \leq \xi^{n,\bf k}_{i+1}(t) + 1.
\end{equation}

This way, the auxiliary height functions $ \xi^{n,\bf k}(\cdot)$ are coupled via the Poisson processes to the height function $z^{n}(\cdot)$ whose evolution we are interested in. As usual, we use this construction to reach a microscopic variational formula for $z^{n}(\cdot)$ in terms of these auxiliary height functions.

\begin{lemma}[The envelope property]
\label{lemma: Envelope property}

 For each $n \in \mathbb{N}$, $i \in \mathbb{Z}$ and $t \in \mathbb{R}_{+}$,

\begin{equation} \label{Envelope Property}
    z^{n}_{i}(t) = \sup_{k \in \mathbb{Z}} \left\{ z^{n}_{k}(0) - \xi^{n,\bf k}_{i-k}(t) \right\}. 
\end{equation}
\end{lemma}

The proof of Lemma \ref{lemma: Envelope property} can be found in Appendix \ref{app:B}, as it is only a technical extension and bookkeeping of the known homogeneous version; see for example Lemma 4.1 in \cite{Seppalainen2001SlowBond}. Note that the only dependence of $\xi^n$ on $z^n$ is the index shift described by $\bf k$.

As usual, \eqref{Envelope Property} will play the role of the microscopic version of the Hopf-Lax formula (see \cite{Evans2010}, Chapter 3) and it will connect the particle system to the relevant PDE. Notice that the initial conditions are separated from the temporal evolution, while the temporal evolution is now restricted to a very special class of initial height functions. 

To define the coupling with the LPP, define the time taken for the auxiliary height function $ \xi^{n,\bf k}_{i}$ to reach level $j \geq -i$ by
\begin{equation} \label{eq:LLPP}
    L^{n, \bf k}_{i,j} = \inf \left\{ t \geq 0: \xi^{n,\bf k}_{i}(t) \geq j \right\}.
\end{equation}
This is naturally defined in the lattice wedge 
$$\mathcal{W} = \{ (i,j): i\in \Z,  j \ge -i \mathds1\{ i \le 0\} \}$$ 
and has boundary conditions prescribed by
\begin{equation}
    L^{n, \bf k}_{i,j} = 0 \text{ for } (i,j) \in \partial \mathcal{W}.
\end{equation}

By looking at constraints \eqref{Server process constraint} and the dynamics of $ \xi^{n,\bf k}_{i}$, we can define $L^{n, \bf k}_{i,j}$ using a recurrence relation:

\begin{equation*}
    L^{n, \bf k}_{i,j} = \max \{ L^{n, \bf k}_{i+1,j-1}, L^{n, \bf k}_{i-1,j}\} +\tilde \omega^{n,\bf k}_{i,j}
\end{equation*}
where $\tilde \omega^{n, \bf k}_{i,j}$ is an exponential waiting time representing the time  for which $ \xi^{n,\bf k}_{i}$ waits to jump from $j-1$ to $j$ once it is able to do so. From the evolution dynamics of the auxiliary height functions, we have that $\tilde \omega^{n,\bf k}_{i,j}$ has a rate $\tilde c^{n}(i+k,j - z^{n}_{k}(0))$.

We want to make use of results from \cite{ciech2021} which are stated for the lattice $\Z^2_+$ (in coordinates $u,v$), while we need them for $\mathcal W$ (in coordinates $(i,j)$), satisfying  
\begin{equation*}
    u = i+j, \hspace{15.pt} v = j.
\end{equation*}
From this transformation, the boundary conditions are now described by 
\begin{equation}
    G^{n,\bf k}_{u,v} = 0 \text{ for } (u,v) \in \partial\Z^2_+.
\end{equation}
The previous recurrence relation becomes
\begin{equation}
    G^{n,\bf k}_{u,v} = \max \{G^{n,\bf k}_{u,v-1}, G^{n, \bf k}_{u-1,v}\} +  \omega^{n, \bf k}_{u,v}
\end{equation}
where $ \omega^{n,\bf k}_{u,v}$ are exponential random variables, with rate $\tilde c^{n}(u-v+k, v - z^{n}_{k}(0))$. Note that 
\be \label{eq:tildec=cshift}
\tilde c (x,y) = c(x+y,y), \quad c(s,t) = \tilde c(s-t,t).
\ee
In the sequel, the superscript $\sim$ will be used when we are referring to the speed function of the particle system, while the absence of it will refer to the speed function in the LPP.

From the boundary conditions, the recurrence relation, and the fact that  
\[
\tilde \omega_{i,j}^{n, \bf k}  \stackrel{d}{=}   \omega_{i+j,j}^{n, \bf k} 
\]
we have
\begin{equation} \label{eq:LtoG}
    \{L^{n, \bf k}_{i,j}: (i,j) \in \mathcal{W} \}\overset{d}{=} \{ G^{n, \bf k}_{i+j,j}: (i+j,j) \in \Z^2_+ \}.
\end{equation}
Note that $G^{n, \bf k}_{i+j,j}$ is the last passage time in the corner growth model 
starting from $(k -z^n_0(k), -z^n_0(k))$.

\subsection{Limiting behaviour of $\xi$.}

A limit for $n^{-1}z^n$ will be obtained by scaling the right-hand side of \eqref{Envelope Property} by $n$ and letting it tend to infinity. For the method to work, we need a scaling limit for the initial conditions, guaranteed by the assumptions once the relevant maximisers are localized, and a scaling limit for $\xi^{n, \bf k}$. We begin with the latter. To do this, we need to define a couple of functions. First, for any $q \in \mathbb{R}$, we define the $q$-shifted last passage times

\begin{align} \label{Gamma q defn}
     \Gamma_{c}^{q}(x,y) 
     &= \Gamma_c((q-v_0(q), -v_0(q)),(x+q-v_0(q), y-v_0(q)) = \Gamma_{T_{q -v_0(q), -v_0(q)}c}(x,y).
\end{align}
From this, we can define functions 
\begin{equation} \label{g^q definition}
    g^{q}(r,t) = \inf \{ y: (r+y,y) \in \mathbb{R}_{+}^{2}, \Gamma_{c}^{q}(r+y,y) \geq t \}.
\end{equation}
These functions $ g^{q}(r,t)$ correspond to the limiting level curves of the height functions $\xi$, evolving with rates $\tilde c$ using \eqref{eq:tildec=cshift}.

\begin{remark}
The first argument of $g^q$ is a displacement measured from the starting
point $q$, not an absolute spatial coordinate. Thus $g^q(r,t)$ denotes the
macroscopic height reached at relative position $r$ by the auxiliary
step-initial process whose origin is attached to the initial site $q$.
Consequently, when we need to denote the physical endpoint to be $x$, the corresponding argument
of $g^q$ will be $x-q$. \qed

\end{remark}

For $x\in\R$, define the lower wedge boundary
\begin{equation}\label{eq:beta-def}
\beta(x):=(-x)_+=\max\{0,-x\}.
\end{equation}
Then the definition \eqref{g^q definition} can be written as
\begin{equation}\label{eq:gq-beta-form}
 g^q(x,t)=\inf\{y\ge \beta(x):\ \Gamma_c^q(x+y,y)\ge t\}.
\end{equation}
In particular,
\begin{equation}\label{eq:levelcurve}
\Gamma_c^q\bigl(x+g_q(x,t),g_q(x,t)\bigr)\ge t.
\end{equation}
The inequality in \eqref{eq:levelcurve} may be strict when the infimum is attained on the active wedge boundary $y=\beta(x)$. However, if the level curve has lifted from the wedge boundary, then equality holds, as shown in the next lemma.

\begin{lemma}\label{lem:interior-level-curve}
Assume \ref{ass:c}, \ref{ass:c2} and \ref{ass:c0}. Fix $q,x\in\mathbb R$ and $t>0$. If
$
 g^q(x,t)>\beta(x),
$
then
\begin{equation}\label{eq:levelcurve-interior}
\Gamma_c^q\bigl(x+g^q(x,t),g^q(x,t)\bigr)=t.
\end{equation}
\end{lemma}

\begin{proof}
By \eqref{eq:gq-beta-form} we already know that
\[
\Gamma_c^q\bigl(x+g^q(x,t),g^q(x,t)\bigr)\ge t.
\]
Assume by contradiction that the inequality is strict. Since $g^q(x,t)>\beta(x)$, we can choose $\e>0$ small enough so that $g_q(x,t)-\varepsilon\ge \beta(x)$. By continuity of $\Gamma_c^q$ on $\R_+^2$,
\[
\Gamma_c^q\bigl(x+g^q(x,t)-\varepsilon,g^q(x,t)-\varepsilon\bigr)>t
\]
for all sufficiently small $\varepsilon>0$, contradicting the infimum property in \eqref{eq:gq-beta-form}. Hence equality must hold.
\end{proof}

It will be convenient to have a version of $\Gamma$ directly in the $\tilde c$ environment, suitable for the TASEP. To this end, define the macroscopic wedge 
\[
\mathbb W = \{ (x,y) \in \R^2: x \ge -y, y \ge 0\}, \quad \partial \mathbb W = \{(x,-x): x \le 0\}\cup\{(x,0): x >0\}.
\]
We also define the wedge ordering: 
\be \label{eq:wedgeorder}
(x,y) \le_{\mathbb W} (w,z) \Longleftrightarrow x+y \le w+z,\, \text{and } y\le z.
\ee

For any $(x,y) \in \mathbb W$ define the set of macroscopic paths from $(0,0)$ to $(x,y)$, as long as $(0,0) \le_{\mathbb W} (x,y)$,
\begin{align*}
\widetilde{\mathcal{H}}(x,y) &=\{{\bf w}=(w_1, w_2)\in C([0,1],\R^2_+):\mathbf{w} \text{ is piecewise }C^1,\mathbf{w}(0)=(0,0), \mathbf{w}(1)=(x,y), \\
& \phantom{xxxxxxxxxxxxxxxxxxxxxxxx} \mathbf w '(s) \in \mathbb W \text{ wherever the derivative is defined}\}.
\end{align*}
The limiting shape function in $\mathbb W$ for the passage times $L$ in \eqref{eq:LLPP} is then given by 
\begin{equation}\label{macroLPTtilde}
\widetilde \Gamma ^q_{\tilde c}(x,y):=\sup_{\mathbf{w}(\cdot)\in\mathcal{\widetilde H}(x,y)}\bigg\{\int_0^1\frac{ \gamma({(w_1 +w_2)}'(s), w_2'(s))}{\tilde c(w_1(s)+q,w_2(s)-v_0(q))}ds\bigg\},
\end{equation}
where $\gamma(x,y)=(\sqrt{x}+\sqrt{y})^2$ is the last-passage constant in a homogeneous rate 1 environment.

Note that a path ${\bf x} = (x_1, x_2) \in \mathcal H (x+y,y)$ corresponds to a path ${\bf w} \in \widetilde{\mathcal H}(x,y)$ by $w_1 = x_1 - x_2, w_2 = x_2$ or equivalently, $x_1 = w_1+w_2, x_2 = w_2$. Now, if one transforms paths from $w$ to $x$, \eqref{macroLPTtilde} becomes \eqref{Gamma q defn} at point $(x+y,y)$.

The fact that \eqref{macroLPTtilde} is the limiting shape function of $L$ passage times then follows immediately from \eqref{eq:LtoG}.
Moreover, we can also write 
\begin{equation} \label{g^q definition til}
    g^{q}(x,t) = \inf \{ y: (x,y) \in  \mathbb W, \widetilde \Gamma_{\tilde c}^{q}(x,y) \geq t \}.
\end{equation}

The next lemma verifies that $g^q(x,t)$ is jointly continuous in all parameters. 
\begin{lemma}\label{lem:gjointc}
Assume \ref{ass:c}, \ref{ass:c2} and \ref{ass:c0}.
The function $g^{q}(x,t)$ is jointly continuous in $(x,q,t)$.
\end{lemma}

\begin{proof}
Define
\[
H(q,x,y):=\Gamma_c^q(x+y,y),
\qquad
D:=\{(q,x,y)\in\mathbb R^3:\ y\ge \beta(x)\}.
\]
By Theorem \ref{thm:2.5}, $\Gamma_c$ is continuous, and since $v_0$ is Lipschitz-$1$, the shift appearing in \eqref{Gamma q defn} depends continuously on $q$. Hence $H$ is jointly continuous on $D$.

For each fixed $(q,x)$, the map $y\mapsto H(q,x,y)$ is continuous and strictly increasing on $[\beta(x),\infty)$. Indeed, if $y_2>y_1\ge \beta(x)$, then both coordinates of
\[
(x+y_2,y_2)
\]
are strictly larger than those of $(x+y_1,y_1)$, and strict monotonicity of $\Gamma_c^q$ in both coordinates gives
\[
H(q,x,y_2)>H(q,x,y_1).
\]
By \eqref{eq:gq-beta-form}, $g^q(x,t)$ is the leftmost point of the superlevel set
\[
\{y\ge \beta(x): H(q,x,y)\ge t\}.
\]
This immediately gives the monotonicities: if $x_1>x_2$, then $g^q(x_1,t)\le g^q(x_2,t)$, and if $t_1>t_2$, then $g^q(x,t_1) \ge g^q(x,t_2)$.

We now prove joint continuity. Let
\[
(q_n,x_n,t_n)\to (q_0,x_0,t_0),
\qquad
y_n:=g^{q_n}(x_n,t_n),
\qquad
y_0:=g^{q_0}(x_0,t_0).
\]
We show that $y_n\to y_0$.

First, the sequence $\{y_n\}_n$ is locally bounded. Indeed, choose $M>y_0+1$ so large that
$
H(q_0,x_0,M)>t_0+1.
$
By continuity of $H$, for all large $n$ we have $H(q_n,x_n,M)>t_n$, and hence by the definition of $y_n$,
\[
\beta(x_n)\le y_n\le M.
\]
Therefore every subsequence of $\{y_n\}_n$ has a convergent further subsequence. Let $y_{n_k}\to y$.

We claim that $y=y_0$. Suppose first that $y<y_0$. Since
\[
H(q_{n_k},x_{n_k},y_{n_k})\ge t_{n_k}
\]
for all $k$, continuity of $H$ yields
\[
H(q_0,x_0,y)\ge t_0.
\]
But $y<y_0$ contradicts the definition of $y_0$ as the infimum of those $y\ge \beta(x_0)$ with $H(q_0,x_0,y)\ge t_0$.

Now suppose that $y>y_0$. Set $\varepsilon=(y-y_0)/2>0$. For all large $k$ we have
\[
y_{n_k}-\varepsilon\ge \beta(x_{n_k}).
\]
By minimality of $y_{n_k}$, this implies $
H(q_{n_k},x_{n_k},y_{n_k}-\varepsilon)<t_{n_k}.
$
Passing to the limit gives
\be \label{eq:coniq}
H(q_0,x_0,y-\varepsilon)\le t_0.
\ee
However, $y-\varepsilon>y_0$, and since $y\mapsto H(q_0,x_0,y)$ is strictly increasing on $[\beta(x_0),\infty)$, we must have
$H(q_0,x_0,y-\varepsilon)>t_0$ if $y_0>\beta(x_0)$, while if $y_0=\beta(x_0)$ then strict increase gives
$
H(q_0,x_0,y-\varepsilon)>H(q_0,x_0,\beta(x_0))\ge t_0.
$

Both cases contradict \eqref{eq:coniq}. Therefore $y=y_0$.

Since every convergent subsequence of $\{y_n\}_n$ has the same limit $y_0$, we conclude that $y_n\to y_0$. This proves joint continuity of $g^q(x,t)$ in $(x,q,t)$.
\end{proof}

As in the homogeneous environment case, level curves can be used to derive large scale limits for the auxiliary height functions $ \xi^{n,\bf k}$ which we show in the following proposition. In other words, the height functions $\xi$ act as microscopic level curves for the last passage time. 

In the sequel we will utilise the vectors 
\[
\fl{\bf{nq}} := (\fl{nq}, -z^n_{\fl{nq}}(0)), \quad q \in \R.  
\]
Note that by the initial LLN for $z^n_{\fl{nx}}(0)$ \eqref{eq:zic} we have 
\be \label{eq:shiftlim}
\lim_{n\to \infty} n^{-1}\fl{\bf{nq}}= (q, -v_0(q)) = \lim_{n\to \infty}n^{-1}(\fl{nq}, -\fl{n v_0(q)}).
\ee
\begin{proposition}
\label{prop: Limit for xi}
For all $x,q \in \mathbb{R}$, and $t \in \mathbb{R}_{+}$,

\begin{equation} \label{auxilliary function limit}
    \lim_{n \to \infty}n^{-1}\xi^{n,\bf \lfloor nq \rfloor}_{\lfloor nx \rfloor}(nt) = g^{q}(x,t) \hspace{15.pt} {\bf P-}\text{a.s.}
\end{equation}
\end{proposition}

\begin{proof}
Use relationship \eqref{eq:LtoG} to map to passage times on $\Z_+^2$ and use Proposition \ref{lem:SPEP} together with equation \eqref{eq:shiftlim} to obtain

\begin{equation} \label{convergence of L to Gamma}
    \lim_{n \to \infty}n^{-1}L^{n,\lfloor {\bf nq} \rfloor}_{\fl{nx}, \fl{ny }} = \Gamma_{c}^{q}(x+y,y).
\end{equation}

Now we want to use this to derive equation \eqref{auxilliary function limit}. Fix $x,q \in \mathbb{R}, t > 0$ and let $A$ be the event (of probability $1$) on which
\begin{equation*}
    n^{-1}L^{n,\lfloor \bf nq \rfloor}_{\fl{ nx' }, \fl{ny}} \to \Gamma^{q}_{c}(x'+y,y), \quad {\bf P}-\text{a.s.}
\end{equation*}
for all $(x',y)$ such that $x'$ is either $x$ or rational, and $y$ can be written as  $y= g^{q}(x',t) \pm \varepsilon$ for rational $\e$. We will be working on $A$ in order to find an upper and lower bound for \eqref{auxilliary function limit}. 

For the lower bound, take $x' > x$ and assume first $g^q(x,t)>0$. Let $\e_0 > 0$ small enough and let rational $x'$ close enough to $x$ and so that $g^{q}(x',t) \geq g^{q}(x,t) - \varepsilon_0>0$ which can be done by Lemma \ref{lem:gjointc}. 

Then we can find a rational $\e < \e_0$ small enough such that $y=g^{q}(x',t) - \varepsilon > 0$. 

Then, by the monotonicity and continuity of $\Gamma^{q}_{c}(x' + y,y)$ in the $y$ variable, there exists a $\delta>0$ such that 
\begin{equation*}
    \Gamma^{q}_{c}(x' + y,y) = \Gamma^{q}_{c}(x'+ g^q(x',t) - \varepsilon, g^{q}(x',t) - \varepsilon) \leq t - \delta.
\end{equation*}
Thus for all $n$ large enough $L^{n,\lfloor \bf nq \rfloor}_{\lfloor nx' \rfloor, \lfloor ny \rfloor} \leq nt - n\delta/2$, which implies with \eqref{eq:LLPP} that for large enough $n$
\begin{align*}
    \xi^{n, \lfloor \bf nq \rfloor}_{\lfloor nx \rfloor}(nt) &\geq \xi^{n, \lfloor \bf nq \rfloor}_{\lfloor nx' \rfloor}(nt) \geq \lfloor ny \rfloor \geq ng^{q}(x',t) - n\varepsilon - 1 \geq ng^{q}(x,t) - 2n\varepsilon_0 - 1.
\end{align*}
First divide both sides by $n$ and take the $\liminf$ as $n\to \infty$. Then, let $\varepsilon_0 \to 0$ on the RHS to finally obtain
\begin{equation*}
    \varliminf_{n \to \infty} n^{-1}\xi^{n,\lfloor \bf nq \rfloor}_{\lfloor nx \rfloor}(nt) \geq g^{q}(x,t)
\end{equation*}
on the full probability event $A$ and when $g^q(x,t)>0$. 

If  $g^q(x,t)=0$ the inequality in the last display is automatically true. 

Finally, the reverse inequality which gives the upper bound follows in a similar way, by approximating $x$ from the left.
\end{proof}

\begin{lemma}\label{lem:jointxicontrol}
For any $\e>0$, $t>0$ and $a, x \in \R$. There exists a $\delta = \delta(x, a, \e)>0$ such that 
\[
\varliminf_{n\to \infty} \frac{1}{n} \min_{ \fl{na} \le k \le \ce{na+n\delta}} \xi_{\fl{nx} - k}^{n,\bf k}(nt) > g^{a}(x-a, t) - \e.
\] 
\end{lemma}

\begin{proof}%[Proof of Lemma \ref{lem:jointxicontrol}] 
Fix an $\e>0$ and assume $\delta < \e/4$. Moreover, fix an $\e_1> 0$ such that 
\be \label{eq:later}
g^{a - \e_1}(x-a+\e_1, t) - g^{a}(x-a, t) > -\e/4 \quad \text{ and } \quad \e_1 < \e/4.
\ee
The first inequality is allowed by Lemma \ref{lem:gjointc}. 

Recall that the process $ \xi^{n,\bf k}$ attempts to increase by 1 at location $i$ at the Poisson event times of $N^{(n)}_{i+k, \xi^{n,\bf k}_i - z_k^n(0)}$. Since $k \in  [\fl{na}, \ce{na+n\delta}]$ and $z_k$ is increasing in $k$, the smallest possible value for $ \xi^{n,\bf k}_i - z_k^n(0)$ is $- z_{\ce{na+n\delta}}^n(0)$ as $ \xi^{n,\bf k}_i \ge 0$. 
Then, by \eqref{eq:zic}, for any $\e_1 > 0$ we can find  $n$ large enough, so that   $$-z_{\ce{na+n\delta}}^n(0) > -\fl{n(v_0(a+\delta) + \e_1)}.$$  

% Similarly,  for any fixed $i$,  the leftmost first index of the Poisson processes involved in the evolution of $ \xi^{n,\bf k}$ at $i$ are $i + k \ge i+ \fl{na} > i + \fl{n(a - \e_1)}$. 

Set \[ K = \fl{n(a - \e_1)} \text{ and } Z =  \fl{n(v_0(a+\delta) + 2 \e_1)}. \]

We now define a new process $\xi^{n, (K, -Z)}:=\xi^{ n, K, Z}$, coupled with the same Poisson processes such that
\[
 \xi^{n, K, Z}_{i}(0) = 
\begin{cases}
0, & i \ge 0\\
-i, & i < 0
\end{cases} 
\] 
where it attempts to jump at location $i$ according to the events of the Poisson  process $N^{(n)}_{i+K, \xi^{n, K, Z}_{i}(t) - Z}$. In other words, this new 
$\xi$-process starts evolving at location $(K,-Z)$.

Now note that, by definition, at time $t = 0$ and for any $i \in \Z$ and $k \in  [\fl{na}, \ce{na+n\delta}]$, we have
\[
-Z +  \xi^{n, K, Z}_{i-K}(0) < -z_k^{n}(0) + \xi^{n, \bf k}_{i-k}(0).
\]
We will show that  for any time $t$ and $k \in  [\fl{na}, \ce{na+n\delta}]$ we have 
\be\label{eq:ok?}
-Z +  \xi^{n, K, Z}_{i-K}(t) \le -z_k^{n}(0) + \xi^{n, \bf k}_{i-k}(t), \quad k \in  [\fl{na}, \ce{na+n\delta}].
\ee
 If the inequality in \eqref{eq:ok?} holds, then we have for $i = \fl{nx}$ and time $nt$ that 
\begin{align*}
-Z +  \xi^{n, K, Z}_{\fl{nx} - K}(t) &\le  \min_{k \in  [\fl{na}, \ce{na+n\delta}]} \{ -z_k^{n}(0) + \xi^{n, \bf k}_{\fl{nx} - k}(t) \} \\
&\le - z_{\fl{na}}^n(0) +  \min_{k \in  [\fl{na}, \ce{na+n\delta}]} \{\xi^{n, \bf k}_{\fl{nx} - k}(t)\}.
\end{align*}
Divide by $n$, and let $n \to \infty$ to obtain, also using Proposition \ref{prop: Limit for xi}:
\[
-v_0(a+ \delta) - 2 \e_1 + g^{a - \e_1}(x - a+\e_1, t) \le -v_0(a) + \varliminf_{n\to \infty} \frac{1}{n}   \min_{k \in  [\fl{na}, \ce{na+n\delta}]} \{\xi^{n, \bf k}_{\fl{nx} - k}(t)\}.
\]
Since $v_0$ is Lipschitz-1 continuous,  $v_0(a+ \delta) - v_0(a) < \delta < \e /4$, while  $g^{a - \e_1}(x-a+\e_1, t) - g^{a}(x-a, t) > -\e/4$ by \eqref{eq:later}. Since $\e_1 < \e/4$ we get the bound. 

It now remains to prove \eqref{eq:ok?}. For this we mimic the proof of the envelope property in Lemma \ref{lemma: Envelope property}. We first order all possible Poisson events up to a time horizon $T$, in the random spatial box $B_{L, U}$ for which we know that the $z^n$- process and also the evolution of the process  $\xi^{n, K, Z}(t)$ is not affected by Poisson processes outside $B_{L, U}$. The existence of such a box $B_{L,U}$ follows just like in the graphical construction in Section \ref{sec:results}, but by looking at pairs of independent Poisson processes.

Equation \eqref{eq:ok?} holds at time 0, and we assume that it remains true up to and including the $m$-th Poisson event (here $m$ can be 0 so the initial induction step is also proven by the argument below), at time $t_m$. Therefore, until time $t_{m+1}^-$ the equation holds, giving for all $i \in \Z,$ and $k \in  [\fl{na}, \ce{na+n\delta}]$
\be \label{eq:ind-hyp}
-Z +  \xi^{n, K, Z}_{i-K}(t_{m+1}^-) \le -z_k^{n}(0) + \xi^{n, \bf k}_{i-k}(t_{m+1}^-)
\ee
Now, assume that, for some reason, at $t_{m+1}$ we can find a $k \in  [\fl{na}, \ce{na+n\delta}]$ and an $i$ such that 
\be\label{eq:contind0}
-Z +  \xi^{n, K, Z}_{i-K}(t_{m+1}) > -z_k^{n}(0) + \xi^{n, \bf k}_{i-k}(t_{m+1}).
\ee
Since the inequality is reversed and $Z, z^n_k(0)$ are constant in time, it must be that $t_{m+1}$ was a jump time for $\xi_{i-K}^{n,K, Z}$ with the jump having taken place. Process $\xi_i^{n,K, Z}$ can only increase by one at the Poisson event, so we must have 
\be \label{eq:touching}
-Z +  \xi^{n, K, Z}_{i-K}(t_{m+1}^-) = -z_k^{n}(0) + \xi^{n, \bf k}_{i-k}(t_{m+1}^-)
\ee
 and the Poisson event is for 
 \[
 N^{(n)}_{i, -Z +  \xi^{n, K, Z}_{i-K}(t_{m+1}^-)} = N^{(n)}_{i, -z_k^{n}(0) + \xi^{n, \bf k}_{i-k}(t_{m+1}^-)}.
\]
This means that both processes $ \xi^{n, K, Z}_{i-K}$ and  $\xi^{n, \bf k}_{i-k}$ attempt to jump at time $t_{m+1}$. 

Since the jump actually occurred, by \eqref{auxiliary height process constraints} we get  
%$\xi^{n, K, Z}_{i-K}(t_{m+1}) =  \xi^{n, K, Z}_{i-K}(t_{m+1}-) + 1$ it must be that 
\[ 
\xi^{n, K, Z}_{i-K-1}(t_{m+1}^-) = \xi^{n, K, Z}_{i-K}(t_{m+1}^-) + 1=\xi^{n, K, Z}_{i-K+1}(t_{m+1}^-) + 1,
\]
and, by \eqref{eq:touching}, 
\begin{align*}
\xi^{n, K, Z}_{i-K-1}(t_{m+1}^-) =  - z_k^n(0) + Z + \xi^{n,\bf k}_{i-k}(t_{m+1}^-) +1\ge - z_k^n(0) + Z + \xi^{n, \bf k}_{i-k-1}(t_{m+1}^-).
\end{align*}
In other words, by the inequality above, and the induction hypothesis \eqref{eq:ind-hyp} for index $i-1$ we have 
\be \label{eq:thisone}
-Z + \xi^{n, K, Z}_{(i-1)-K}(t_{m+1}^-) =  - z^n_k(0) + \xi^{n, \bf k}_{(i-1)^-k}(t_{m+1}^-).
\ee
Moreover, 
\begin{align*}
\xi^{n, K, Z}_{i+1-K}(t_{m+1}^-)= \xi^{n, K, Z}_{i-K}(t_{m+1}^-) = -z^n_k(0) + Z + \xi^{n, \bf k}_{i-k}(t_{m+1}^-)\ge -z^n_k(0) + Z + \xi^{n, \bf k}_{i+1-k}(t_{m+1}^-)
\end{align*}
and since the induction hypothesis \eqref{eq:ind-hyp} needs to be true for index $i+1$ we have 
\be \label{eq:andthisone}
-Z+ \xi^{n, K, Z}_{i+1-K}(t_{m+1}^-)= -z^n_k(0) + \xi^{n, \bf k}_{i+1-k}(t_{m+1}^-).
\ee
Equations \eqref{eq:touching},\eqref{eq:thisone} and \eqref{eq:andthisone} give that $\xi_{i-k}^{n, \bf k}$ can also jump at $t_{m+1}$. This contradicts our assumption \eqref{eq:contind0}. Therefore the inequality is preserved for the event at $t_{m+1}$.
\end{proof}

\subsection{Hydrodynamic limit for the height function.}

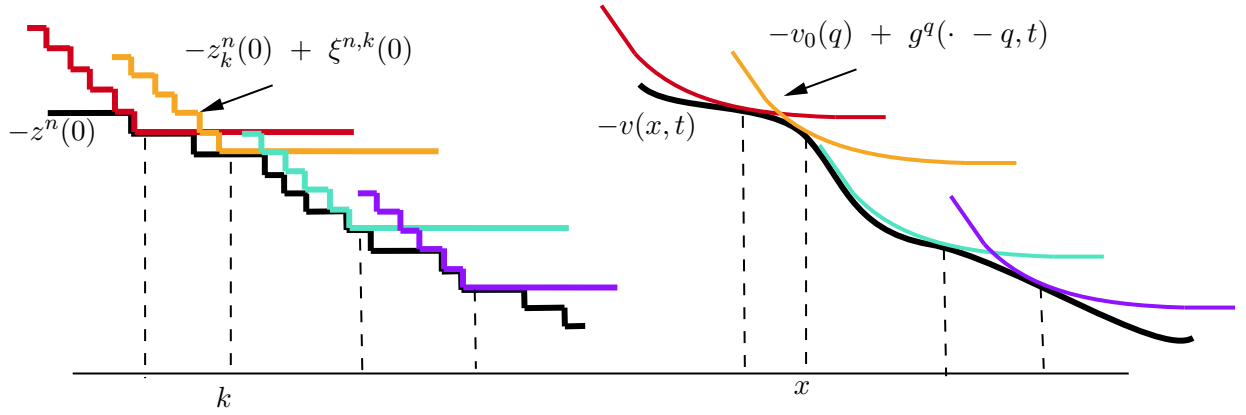
\begin{figure}[ht]

\tikzset{every picture/.style={line width=0.75pt}} %set default line width to 0.75pt        

\begin{tikzpicture}[x=0.75pt,y=0.75pt,yscale=-1,xscale=1]
%uncomment if require: \path (0,2283); %set diagram left start at 0, and has height of 2283

%Straight Lines [id:da6673877640702358] 
\draw [color={rgb, 255:red, 0; green, 0; blue, 0 }  ,draw opacity=1 ][line width=2.25]    (37.13,2049.95) -- (78.96,2049.86) ;
%Straight Lines [id:da6138225464756363] 
\draw [color={rgb, 255:red, 0; green, 0; blue, 0 }  ,draw opacity=1 ][line width=2.25]    (78.98,2060.31) -- (110.43,2060.25) ;
%Straight Lines [id:da5217841123960628] 
\draw [color={rgb, 255:red, 0; green, 0; blue, 0 }  ,draw opacity=1 ][line width=2.25]    (110.46,2070.7) -- (146.06,2070.63) ;
%Straight Lines [id:da19416843975554554] 
\draw [color={rgb, 255:red, 0; green, 0; blue, 0 }  ,draw opacity=1 ][line width=2.25]    (146.08,2081.08) -- (155.75,2081.06) ;
%Straight Lines [id:da04936494562414051] 
\draw [color={rgb, 255:red, 0; green, 0; blue, 0 }  ,draw opacity=1 ][line width=2.25]    (155.96,2090.27) -- (167.53,2090.25) ;
%Straight Lines [id:da049784146147514474] 
\draw [color={rgb, 255:red, 0; green, 0; blue, 0 }  ,draw opacity=1 ][line width=2.25]    (166.87,2099.46) -- (186.91,2099.2) ;
%Straight Lines [id:da515178694317613] 
\draw [color={rgb, 255:red, 0; green, 0; blue, 0 }  ,draw opacity=1 ][line width=2.25]    (187.12,2108.63) -- (199.23,2108.67) ;
%Straight Lines [id:da409937545164495] 
\draw [color={rgb, 255:red, 0; green, 0; blue, 0 }  ,draw opacity=1 ][line width=2.25]    (187.12,2108.63) -- (186.91,2099.2) ;
%Straight Lines [id:da7685700242957959] 
\draw [color={rgb, 255:red, 0; green, 0; blue, 0 }  ,draw opacity=1 ][line width=2.25]    (166.87,2099.46) -- (166.85,2090.41) ;
%Straight Lines [id:da5708052634052408] 
\draw [color={rgb, 255:red, 0; green, 0; blue, 0 }  ,draw opacity=1 ][line width=2.25]    (155.77,2090.27) -- (155.75,2080.84) ;
%Straight Lines [id:da27558848014428694] 
\draw [color={rgb, 255:red, 0; green, 0; blue, 0 }  ,draw opacity=1 ][line width=2.25]    (146.08,2081.08) -- (146.06,2070.63) ;
%Straight Lines [id:da1448806212411189] 
\draw [color={rgb, 255:red, 0; green, 0; blue, 0 }  ,draw opacity=1 ][line width=2.25]    (110.46,2070.7) -- (110.43,2060.25) ;
%Straight Lines [id:da17377128814852405] 
\draw [color={rgb, 255:red, 0; green, 0; blue, 0 }  ,draw opacity=1 ][line width=2.25]    (78.98,2060.31) -- (78.96,2049.86) ;
%Straight Lines [id:da29258955394457165] 
\draw [color={rgb, 255:red, 0; green, 0; blue, 0 }  ,draw opacity=1 ][line width=2.25]    (199.25,2119.13) -- (234.86,2119.05) ;
%Straight Lines [id:da934525277962692] 
\draw [color={rgb, 255:red, 0; green, 0; blue, 0 }  ,draw opacity=1 ][line width=2.25]    (234.88,2129.51) -- (245.07,2129.27) ;
%Straight Lines [id:da5673167358856972] 
\draw [color={rgb, 255:red, 0; green, 0; blue, 0 }  ,draw opacity=1 ][line width=2.25]    (244.75,2138.7) -- (276.21,2138.63) ;
%Straight Lines [id:da18107982274624523] 
\draw [color={rgb, 255:red, 0; green, 0; blue, 0 }  ,draw opacity=1 ][line width=2.25]    (276.41,2147.84) -- (296.46,2147.58) ;
%Straight Lines [id:da5789930119935943] 
\draw [color={rgb, 255:red, 0; green, 0; blue, 0 }  ,draw opacity=1 ][line width=2.25]    (296.66,2157.01) -- (306.85,2156.77) ;
%Straight Lines [id:da36424274926788547] 
\draw [color={rgb, 255:red, 0; green, 0; blue, 0 }  ,draw opacity=1 ][line width=2.25]    (296.66,2157.01) -- (296.46,2147.58) ;
%Straight Lines [id:da18752063706123345] 
\draw [color={rgb, 255:red, 0; green, 0; blue, 0 }  ,draw opacity=1 ][line width=2.25]    (276.41,2147.84) -- (276.21,2138.41) ;
%Straight Lines [id:da08556943708454268] 
\draw [color={rgb, 255:red, 0; green, 0; blue, 0 }  ,draw opacity=1 ][line width=2.25]    (244.75,2138.7) -- (244.55,2129.27) ;
%Straight Lines [id:da8488488697651487] 
\draw [color={rgb, 255:red, 0; green, 0; blue, 0 }  ,draw opacity=1 ][line width=2.25]    (234.88,2129.51) -- (234.86,2119.05) ;
%Straight Lines [id:da049223085598394944] 
\draw [color={rgb, 255:red, 0; green, 0; blue, 0 }  ,draw opacity=1 ][line width=2.25]    (199.25,2119.13) -- (199.23,2108.67) ;
%Straight Lines [id:da2493773106760897] 
\draw [color={rgb, 255:red, 208; green, 2; blue, 27 }  ,draw opacity=1 ][line width=2.25]    (48.57,2028.01) -- (58.24,2027.99) ;
%Straight Lines [id:da36135744564503225] 
\draw [color={rgb, 255:red, 208; green, 2; blue, 27 }  ,draw opacity=1 ][line width=2.25]    (58.44,2038.17) -- (70.02,2038.14) ;
%Straight Lines [id:da08808191801738252] 
\draw [color={rgb, 255:red, 208; green, 2; blue, 27 }  ,draw opacity=1 ][line width=2.25]    (58.26,2038.17) -- (58.24,2027.75) ;
%Straight Lines [id:da8232696758035489] 
\draw [color={rgb, 255:red, 208; green, 2; blue, 27 }  ,draw opacity=1 ][line width=2.25]    (48.57,2028.01) -- (48.55,2016.46) ;
%Straight Lines [id:da12190096958899066] 
\draw [color={rgb, 255:red, 208; green, 2; blue, 27 }  ,draw opacity=1 ][line width=2.25]    (70.87,2049.29) -- (80.54,2049.27) ;
%Straight Lines [id:da1179634997733443] 
\draw [color={rgb, 255:red, 208; green, 2; blue, 27 }  ,draw opacity=1 ][line width=2.25]    (80.75,2059.45) -- (190.74,2059.49) ;
%Straight Lines [id:da4553311418197439] 
\draw [color={rgb, 255:red, 208; green, 2; blue, 27 }  ,draw opacity=1 ][line width=2.25]    (80.56,2059.45) -- (80.54,2049.03) ;
%Straight Lines [id:da9715681313652025] 
\draw [color={rgb, 255:red, 208; green, 2; blue, 27 }  ,draw opacity=1 ][line width=2.25]    (70.87,2049.29) -- (70.85,2037.75) ;
%Straight Lines [id:da20391211271810672] 
\draw [color={rgb, 255:red, 208; green, 2; blue, 27 }  ,draw opacity=1 ][line width=2.25]    (26.78,2007.26) -- (36.45,2007.24) ;
%Straight Lines [id:da40636788425178694] 
\draw [color={rgb, 255:red, 208; green, 2; blue, 27 }  ,draw opacity=1 ][line width=2.25]    (36.66,2017.42) -- (48.24,2017.39) ;
%Straight Lines [id:da16690580496015983] 
\draw [color={rgb, 255:red, 208; green, 2; blue, 27 }  ,draw opacity=1 ][line width=2.25]    (36.47,2017.42) -- (36.45,2007) ;

%Straight Lines [id:da01975290456446388] 
\draw [color={rgb, 255:red, 245; green, 166; blue, 35 }  ,draw opacity=1 ][line width=2.25]    (91.1,2040.63) -- (100.77,2040.61) ;
%Straight Lines [id:da5653276355447281] 
\draw [color={rgb, 255:red, 245; green, 166; blue, 35 }  ,draw opacity=1 ][line width=2.25]    (100.98,2049.82) -- (112.55,2049.79) ;
%Straight Lines [id:da6857084645353757] 
\draw [color={rgb, 255:red, 245; green, 166; blue, 35 }  ,draw opacity=1 ][line width=2.25]    (100.79,2049.82) -- (100.77,2040.39) ;
%Straight Lines [id:da8633366336283944] 
\draw [color={rgb, 255:red, 245; green, 166; blue, 35 }  ,draw opacity=1 ][line width=2.25]    (91.1,2040.63) -- (91.08,2030.18) ;
%Straight Lines [id:da28565761478179086] 
\draw [color={rgb, 255:red, 245; green, 166; blue, 35 }  ,draw opacity=1 ][line width=2.25]    (113.4,2059.89) -- (123.07,2059.87) ;
%Straight Lines [id:da946381547128662] 
\draw [color={rgb, 255:red, 245; green, 166; blue, 35 }  ,draw opacity=1 ][line width=2.25]    (123.28,2069.08) -- (233.27,2069.12) ;
%Straight Lines [id:da8544166703155145] 
\draw [color={rgb, 255:red, 245; green, 166; blue, 35 }  ,draw opacity=1 ][line width=2.25]    (123.09,2069.08) -- (123.07,2059.65) ;
%Straight Lines [id:da81832822462964] 
\draw [color={rgb, 255:red, 245; green, 166; blue, 35 }  ,draw opacity=1 ][line width=2.25]    (113.4,2059.89) -- (113.38,2049.44) ;
%Straight Lines [id:da45053531868489405] 
\draw [color={rgb, 255:red, 245; green, 166; blue, 35 }  ,draw opacity=1 ][line width=2.25]    (69.32,2021.85) -- (78.99,2021.83) ;
%Straight Lines [id:da26053535498108693] 
\draw [color={rgb, 255:red, 245; green, 166; blue, 35 }  ,draw opacity=1 ][line width=2.25]    (79.19,2031.04) -- (90.77,2031.01) ;
%Straight Lines [id:da11451367769753773] 
\draw [color={rgb, 255:red, 245; green, 166; blue, 35 }  ,draw opacity=1 ][line width=2.25]    (79.01,2031.04) -- (78.99,2021.61) ;

%Straight Lines [id:da2875549324579152] 
\draw [color={rgb, 255:red, 80; green, 227; blue, 194 }  ,draw opacity=1 ][line width=2.25]    (156.46,2079.15) -- (166.12,2079.13) ;
%Straight Lines [id:da49728822253553895] 
\draw [color={rgb, 255:red, 80; green, 227; blue, 194 }  ,draw opacity=1 ][line width=2.25]    (166.33,2088.34) -- (177.91,2088.32) ;
%Straight Lines [id:da8133035545592245] 
\draw [color={rgb, 255:red, 80; green, 227; blue, 194 }  ,draw opacity=1 ][line width=2.25]    (166.14,2088.34) -- (166.12,2078.91) ;
%Straight Lines [id:da2838879804066037] 
\draw [color={rgb, 255:red, 80; green, 227; blue, 194 }  ,draw opacity=1 ][line width=2.25]    (156.46,2079.15) -- (156.43,2068.7) ;
%Straight Lines [id:da3868517304597887] 
\draw [color={rgb, 255:red, 80; green, 227; blue, 194 }  ,draw opacity=1 ][line width=2.25]    (178.76,2098.42) -- (188.43,2098.4) ;
%Straight Lines [id:da080992522823922] 
\draw [color={rgb, 255:red, 80; green, 227; blue, 194 }  ,draw opacity=1 ][line width=2.25]    (188.63,2107.61) -- (298.62,2107.65) ;
%Straight Lines [id:da18063047191664983] 
\draw [color={rgb, 255:red, 80; green, 227; blue, 194 }  ,draw opacity=1 ][line width=2.25]    (188.45,2107.61) -- (188.43,2098.18) ;
%Straight Lines [id:da26722754499631274] 
\draw [color={rgb, 255:red, 80; green, 227; blue, 194 }  ,draw opacity=1 ][line width=2.25]    (178.76,2098.42) -- (178.74,2087.96) ;
%Straight Lines [id:da8647540607375217] 
\draw [color={rgb, 255:red, 80; green, 227; blue, 194 }  ,draw opacity=1 ][line width=2.25]    (134.67,2060.37) -- (144.34,2060.35) ;
%Straight Lines [id:da7264720528316383] 
\draw [color={rgb, 255:red, 80; green, 227; blue, 194 }  ,draw opacity=1 ][line width=2.25]    (144.55,2069.56) -- (156.12,2069.54) ;
%Straight Lines [id:da8958138112428368] 
\draw [color={rgb, 255:red, 80; green, 227; blue, 194 }  ,draw opacity=1 ][line width=2.25]    (144.36,2069.56) -- (144.34,2060.13) ;

%Straight Lines [id:da7581084601779554] 
\draw [color={rgb, 255:red, 144; green, 19; blue, 254 }  ,draw opacity=1 ][line width=2.25]    (214.13,2109.01) -- (223.62,2108.99) ;
%Straight Lines [id:da5178347897930962] 
\draw [color={rgb, 255:red, 144; green, 19; blue, 254 }  ,draw opacity=1 ][line width=2.25]    (223.82,2118.2) -- (235.18,2118.18) ;
%Straight Lines [id:da09423801872532733] 
\draw [color={rgb, 255:red, 144; green, 19; blue, 254 }  ,draw opacity=1 ][line width=2.25]    (223.64,2118.2) -- (223.62,2108.77) ;
%Straight Lines [id:da4009600903378797] 
\draw [color={rgb, 255:red, 144; green, 19; blue, 254 }  ,draw opacity=1 ][line width=2.25]    (214.13,2109.01) -- (214.11,2098.56) ;
%Straight Lines [id:da8292851415245306] 
\draw [color={rgb, 255:red, 144; green, 19; blue, 254 }  ,draw opacity=1 ][line width=2.25]    (236.01,2128.27) -- (245.5,2128.25) ;
%Straight Lines [id:da61312013230129] 
\draw [color={rgb, 255:red, 144; green, 19; blue, 254 }  ,draw opacity=1 ][line width=2.25]    (245.7,2137.46) -- (323,2137.5) ;
%Straight Lines [id:da5703953026222024] 
\draw [color={rgb, 255:red, 144; green, 19; blue, 254 }  ,draw opacity=1 ][line width=2.25]    (245.52,2137.46) -- (245.5,2128.03) ;
%Straight Lines [id:da7373432363902237] 
\draw [color={rgb, 255:red, 144; green, 19; blue, 254 }  ,draw opacity=1 ][line width=2.25]    (236.01,2128.27) -- (235.99,2117.82) ;
%Straight Lines [id:da22650228744776968] 
\draw [color={rgb, 255:red, 144; green, 19; blue, 254 }  ,draw opacity=1 ][line width=2.25]    (192.76,2090.23) -- (202.25,2090.21) ;
%Straight Lines [id:da12302626756418389] 
\draw [color={rgb, 255:red, 144; green, 19; blue, 254 }  ,draw opacity=1 ][line width=2.25]    (202.45,2099.42) -- (213.81,2099.4) ;
%Straight Lines [id:da057447325630865964] 
\draw [color={rgb, 255:red, 144; green, 19; blue, 254 }  ,draw opacity=1 ][line width=2.25]    (202.27,2099.42) -- (202.25,2089.99) ;
%Straight Lines [id:da334150239978104] 
\draw  [dash pattern={on 4.5pt off 4.5pt}]  (85.96,2062.38) -- (86,2183) ;
%Straight Lines [id:da6908167539043386] 
\draw  [dash pattern={on 4.5pt off 4.5pt}]  (129.01,2074.42) -- (129,2180.64) ;
%Straight Lines [id:da981544017252813] 
\draw  [dash pattern={on 4.5pt off 4.5pt}]  (193.85,2112.94) -- (195,2179.86) ;
%Straight Lines [id:da37790806109866204] 
\draw  [dash pattern={on 4.5pt off 4.5pt}]  (251.42,2138.47) -- (252,2178.29) ;
%Curve Lines [id:da38973255296963305] 
\draw [line width=2.25]    (334.13,2035.61) .. controls (348,2050.21) and (399.58,2042.9) .. (417.71,2061.87) .. controls (435.85,2080.83) and (439.66,2108.86) .. (479.34,2115.77) .. controls (519.02,2122.68) and (594.89,2173.21) .. (611.89,2162.77) ;
%Curve Lines [id:da9873348135352068] 
\draw [color={rgb, 255:red, 208; green, 2; blue, 27 }  ,draw opacity=1 ][line width=1.5]    (331.89,2020.19) .. controls (352.99,2042.99) and (384.23,2050.59) .. (432.35,2051.98) ;
%Straight Lines [id:da5439657814675882] 
\draw [color={rgb, 255:red, 208; green, 2; blue, 27 }  ,draw opacity=1 ][line width=1.5]    (315,1996) -- (331.89,2020.19) ;
%Straight Lines [id:da1877449296420447] 
\draw [color={rgb, 255:red, 208; green, 2; blue, 27 }  ,draw opacity=1 ][line width=1.5]    (432.35,2051.98) -- (457.68,2051.98) ;

%Curve Lines [id:da881242870518263] 
\draw [color={rgb, 255:red, 245; green, 166; blue, 35 }  ,draw opacity=1 ][line width=1.5]    (397.45,2043.21) .. controls (418.56,2066.02) and (449.8,2073.62) .. (497.92,2075) ;
%Straight Lines [id:da9987211643388759] 
\draw [color={rgb, 255:red, 245; green, 166; blue, 35 }  ,draw opacity=1 ][line width=1.5]    (380.57,2019.02) -- (397.45,2043.21) ;
%Straight Lines [id:da9316334957993592] 
\draw [color={rgb, 255:red, 245; green, 166; blue, 35 }  ,draw opacity=1 ][line width=1.5]    (497.92,2075) -- (523.25,2075) ;

%Curve Lines [id:da9342458294204039] 
\draw [color={rgb, 255:red, 80; green, 227; blue, 194 }  ,draw opacity=1 ][line width=1.5]    (441.35,2090.2) .. controls (462.46,2113.01) and (493.7,2120.61) .. (541.82,2121.99) ;
%Straight Lines [id:da42329716767478875] 
\draw [color={rgb, 255:red, 80; green, 227; blue, 194 }  ,draw opacity=1 ][line width=1.5]    (424.47,2066.02) -- (441.35,2090.2) ;
%Straight Lines [id:da5114293196574284] 
\draw [color={rgb, 255:red, 80; green, 227; blue, 194 }  ,draw opacity=1 ][line width=1.5]    (541.82,2121.99) -- (567.15,2121.99) ;

%Curve Lines [id:da7217606350480529] 
\draw [color={rgb, 255:red, 144; green, 19; blue, 254 }  ,draw opacity=1 ][line width=1.5]    (507.2,2115.77) .. controls (528.31,2138.58) and (559.55,2146.18) .. (607.67,2147.56) ;
%Straight Lines [id:da877314637603216] 
\draw [color={rgb, 255:red, 144; green, 19; blue, 254 }  ,draw opacity=1 ][line width=1.5]    (490.32,2091.59) -- (507.2,2115.77) ;
%Straight Lines [id:da11908765838667779] 
\draw [color={rgb, 255:red, 144; green, 19; blue, 254 }  ,draw opacity=1 ][line width=1.5]    (607.67,2147.56) -- (633,2147.56) ;

%Straight Lines [id:da32004685997806037] 
\draw  [dash pattern={on 4.5pt off 4.5pt}]  (385.96,2051.38) -- (387,2180.64) ;
%Straight Lines [id:da014673536325762737] 
\draw  [dash pattern={on 4.5pt off 4.5pt}]  (417.71,2061.87) -- (417.71,2178.29) ;
%Straight Lines [id:da17986874035749778] 
\draw  [dash pattern={on 4.5pt off 4.5pt}]  (486.85,2117.66) -- (488,2182) ;
%Straight Lines [id:da709041419865948] 
\draw  [dash pattern={on 4.5pt off 4.5pt}]  (535.42,2136.9) -- (537,2178.29) ;
%Straight Lines [id:da996787204810329] 
\draw    (49.64,2180.64) -- (579.64,2180.64) ;
%Straight Lines [id:da8297904206298349] 
\draw    (150,2036) -- (115.26,2048.75) ;
\draw [shift={(113.38,2049.44)}, rotate = 339.85] [fill={rgb, 255:red, 0; green, 0; blue, 0 }  ][line width=0.08]  [draw opacity=0] (12,-3) -- (0,0) -- (12,3) -- cycle    ;
%Straight Lines [id:da6034536468244966] 
\draw    (443,2028) -- (408.26,2040.75) ;
\draw [shift={(406.38,2041.44)}, rotate = 339.85] [fill={rgb, 255:red, 0; green, 0; blue, 0 }  ][line width=0.08]  [draw opacity=0] (12,-3) -- (0,0) -- (12,3) -- cycle    ;

% Text Node
\draw (410,2183) node [anchor=north west][inner sep=0.75pt]    {$x$};
% Text Node
\draw (397,2002) node [anchor=north west][inner sep=0.75pt]    {$-v_{0}( q) \ +\ g^{q}( \cdot \ -q,t)$};
% Text Node
\draw (104,2008) node [anchor=north west][inner sep=0.75pt]    {$-z_{k}^{n}( 0) \ +\ \xi ^{n,k}( 0)$};
% Text Node
\draw (120,2186) node [anchor=north west][inner sep=0.75pt]    {$k$};
% Text Node
\draw (311,2048.02) node [anchor=north west][inner sep=0.75pt]    {$-v( x,t)$};
% Text Node
\draw (15.36,2050.91) node [anchor=north west][inner sep=0.75pt]    {$-z^{n}( 0)$};

\end{tikzpicture}
\caption{{\bf Left:} Representation of the envelope property, showing  relationship between $-z^n(0)$ and the auxiliary processes $\xi^n$. Each $\xi$ is placed so that it is above $-z^n$ but touching at the starting point $(k,-z^n_k(0))$. As all coupled processes evolve, all $\xi^n$ remain above, and at least one remains tangent to $-z^n_k(t)$ for any $k$ at each spatial point.\\
{\bf Right:} The same relationship remains in the limit; the scaled height function converges to $-v(x,t)$ and the scaled $\xi$ converge to level curves of passage times. This gives rise to the variational formula \eqref{v tilde definition}}
\label{fig:envelope property}
\end{figure}

Having found the scaling limits for the auxiliary height functions $\xi^{n,k}$, we substitute these into equation \eqref{Envelope Property} from Lemma \ref{lemma: Envelope property} to find scaling limits for $z^{n}$. 

We define the candidate limit to be
\begin{equation} \label{v tilde definition}
    v(x,t) = \sup_{q \in \mathbb{R}} \left\{v_{0}(q) - g^{q}(x-q,t) \right\}. 
\end{equation}

\begin{proposition}
\label{prop: Limit for z}
For all $x \in \bR$ and $t \in \bR_{+}$,
\begin{equation}
    \lim_{n \to \infty}n^{-1}z^{n}_{\fl{nx}}(nt) = v(x,t).
\end{equation}
\end{proposition}

\begin{remark} \label{remark: limit for z}
Because the height functions $z^{n}$ are linked to exclusion particle configurations $\eta^n$ through \eqref{Server-exclusion relation0}, we can also write Proposition \ref{prop: Limit for z} in terms of $\eta^n$. In fact, for any $b > a$ we have
\begin{equation*}
    \lim_{n \to \infty}n^{-1}\sum^{\fl{nb}}_{i = \fl{na}}\eta^n_i(nt) = \lim_{n \to \infty}n^{-1}\left(z^n_{\fl{nb}} - z^n_{\fl{na}} \right) = {v}(b,t) - {v}(a,t). \qed
\end{equation*}
\end{remark}

Before proving Proposition \ref{prop: Limit for z} fully we need one more lemma which gives us a useful way to truncate \eqref{Envelope Property}.  Fix $x,t$ and define for any $-\infty < r_{1} < r_{2} < \infty$ the truncated variational formula,

\begin{equation}
      \zeta^{n}_{x,t}(r_{1},r_{2}) = \sup_{ q \in [ r_{1},   r_{2} ]} \left\{ z^{n}_{\lfloor nq \rfloor}(0) - \xi^{n,\lfloor \bf nq \rfloor}_{\lfloor nx \rfloor - \lfloor nq \rfloor}(nt) \right\}.
\end{equation}

\begin{lemma}
\label{lemma: truncation for q}
 Fix $(x,t) \in \mathbb{R} \times \mathbb{R}_{+}$ and suppose that the speed function $c$ satisfies Assumption \ref{assumption: growth rates}. Then we can find deterministic functions $R_L(x,t)$ and $R_U(x,t)$ so that $-\infty < R_{L} < 0  < R_{U} < +\infty$ such that, for every $r_1 \in (-\infty, R_L]$, $r_2 \in [R_U, \infty)$ and for all $n$ large enough there exists a $C=C(x,t)$ such that 
\begin{equation}\label{eq:qtr}
    \mathbf{P} \left\{ z^{n}_{\lfloor nx \rfloor}(nt) \neq   \zeta^{n}_{x,t}(r_{1}, r_{2}) \right\} \le e^{-C(x, t) n}.
\end{equation}
\end{lemma}

\begin{proof}
Fix $(x,t) \in \bR \times \bR_{+}$ and $n \in \bN$. Suppose that there exists a $y \in \bR$ such that $\fl{nx} - \fl{ny} > 0$ and $\xi^{n,\fl{\bf ny}}_{\fl{nx}-\fl{ny}}(nt) = 0$, it then follows that for every $y' < y$, 
\begin{equation*}
    \xi^{n,\fl{\bf ny'}}_{\fl{nx}-\fl{ny'}}(nt) \geq \xi^{n,\fl{\bf ny}}_{\fl{nx}-\fl{ny}}(nt) = 0.
\end{equation*}
This and the monotonicity of $z^n$ imply that
\begin{equation*}
    z^{n}_{\fl{ny'}}(0) - \xi^{n,\fl{\bf ny'}}_{\fl{nx}-\fl{ny'}}(nt) \leq z^{n}_{\fl{ny}}(0) - \xi^{n,\fl{\bf ny}}_{\fl{nx}-\fl{ny}}(nt),
\end{equation*}
and hence that the supremum in \eqref{Envelope Property} does not depend on any $y' < y$. 

The same can be argued for the other side. Suppose that there exists a $v \in \bR$ such that $\fl{nx} - \fl{nv} < 0$ and $\xi^{n,\fl{\bf nv}}_{\fl{nx}-\fl{nv}}(nt) = \fl{nv} - \fl{nx}$. Take any $v' > v$ to write  
\begin{equation*}
    \xi^{n,\fl{\bf nv'}}_{\fl{nx}-\fl{nv'}}(nt) \geq \fl{nv'} - \fl{nx} \geq  \xi^{n,\fl{\bf nv}}_{\fl{nx}-\fl{nv}}(nt) + \fl{nv'} - \fl{nv} \ge \xi^{n,\fl{\bf nv}}_{\fl{nx}-\fl{nv}}(nt) + z_{\fl{nv'}}^n(0) - z^n_{\fl{nv}}(0).
\end{equation*}
The last inequality follows from the Lipschitz-1 property of $z^n$ in \eqref{eq:DiscLip}. Rearranging the terms
\begin{equation*}
    z^{n}_{\fl{nv'}}(0) - \xi^{n,\fl{\bf nv'}}_{\fl{nx}-\fl{nv'}}(nt) \leq z^{n}_{\fl{nv}}(0) - \xi^{n,\fl{\bf nv}}_{\fl{nx}-\fl{nv}}(nt),
\end{equation*}
and hence, that the supremum in \eqref{Envelope Property} does not depend on any $v' > v$.

Thus, to prove the lemma, it suffices to find $r_1(x,t), r_2(x,t)$ such that for all $n$ large enough 
\begin{equation}
    \mathbf{P} \left\{ \exists N < \fl{nr_1} :  \xi^{n,  \bf N }_{\lfloor nx \rfloor - N }(nt) > 0 \right\} \leq e^{-C(x,t) n},
\end{equation}
and
\begin{equation}
   \mathbf{P} \left\{  \exists N > \fl{nr_2} :  \xi^{n,  \bf N }_{\lfloor nx \rfloor - N}(nt) > N - \lfloor nx \rfloor \right\} \leq e^{-C(x,t)n}, 
\end{equation}
for some constant $C(x,t)$ and where $N \in \Z$. %Then the Lemma follows from the first Borel-Cantelli Lemma. 

Let us begin with the first statement. Initially $\xi^{n, \lfloor \bf nr_{1} \rfloor}_{\lfloor nx \rfloor - \lfloor nr_{1} \rfloor}(0) = 0$ for $\lfloor nr_{1} \rfloor < \lfloor nx \rfloor$ by construction. Now we bound the first jump time 
\begin{equation*}
    \tau_{1}^{n} = \inf \{s \in \mathbb{R}_{+} : \xi^{n, \lfloor \bf nr_{1} \rfloor}_{\lfloor nx \rfloor - \lfloor nr_{1} \rfloor}(s) > 0 \}.
\end{equation*}
That is the time when the particle located at $-1$ of process $\xi^{n, \lfloor \bf nr_{1} \rfloor}$ moves to location $\lfloor nx \rfloor - \lfloor nr_{1} \rfloor$. The movement of the first particle is unobstructed, so at each site $i$ the first particle waits for an exponentially distributed amount of time with rate $ c^{n}_{i, r_1} = \tilde c(n^{-1}(i + \fl{nr_1}) , -n^{-1}z^n_{\fl{nr_1}}(0))$ and then jumps.
We can then write
\begin{align*}
  {\bf P} \left\{\xi^{n,\lfloor \bf nr_{1} \rfloor}_{\lfloor nx \rfloor - \lfloor nr_{1} \rfloor}(nt) > 0 \right\} =  {\bf P}\{ \tau^{n}_{1} < nt \} =  {\bf P}\bigg\{ \sum_{k = 0}^{\fl{nx} - \fl{nr_1}-1} \frac{1}{ c^{n}_{k-1, r_1}} \omega_{k}< nt \bigg\},
\end{align*}
where the $\omega_k$ are i.i.d. Exp(1) random variables. Note that 
\[
 \sup_{i \in [0, \fl{nx} - \fl{nr_1}] } c^{n}_{i, r_1} \le  \sup_{ (u,v) \in R^{\bf 0}_{ |x| + 2|r_1|, 2|r_1|}} c(u,v) = M^{\bf 0}_{x, r_1},
\]
where the rectangle is defined in \eqref{eq:rectangle} and by Assumption \ref{assumption: growth rates} we may bound 
\begin{align}
 {\bf P}\bigg\{ \sum_{k = 0}^{\fl{nx} - \fl{nr_1}-1} \frac{1}{ c^{n}_{k-1, r_1}} \omega_{k}< nt \bigg\} &\le   {\bf P}\bigg\{ \sum_{k = 0}^{\fl{nx} - \fl{nr_1}-1} \omega_{k}< nt M^{\bf 0}_{x, r_1} \bigg\} \label{eq:pild}\\ 
 & \le   {\bf P}\bigg\{ \sum_{k = 0}^{\fl{nx} - \fl{nr_1}-1} \omega_{k}< nt \log(|x| + 2 |r_1|) \bigg\}, \textrm{ for $r_1$ large.} \notag
\end{align}
Fix an $\e > 0$ and implicitly define a function $R_L=R_L(x,t, \e) < 0 $ so that 
\[
t = \frac{(x - R_L)(1 - \e)}{\log(|x| + 2|R_L|)}.
\]
Then, by possibly reducing $R_L$ further, we have that 
for any  $r_1 <  R_L(x,t, \e)$  
\be\label{eq:x-tdep}
t < \frac{(x - r_1)(1 - \e)}{\log(|x| + 2|r_1|)},
\ee
and the right-hand side in \eqref{eq:x-tdep} is provably monotonically increasing in $|r_1|$. This makes the last probability in \eqref{eq:pild} a large deviation probability for i.i.d.~ sums of exponential random variable, which we can bound by $e^{ - n(x - r_1)C(\e)}$ by a Chernoff bound, for all $n(x-r_1)>n_0$ large enough.

When $n_0$ and $r_1$ are identified we can further improve the bound. Fix $n > n_0$ and let $N < nr_1$, $N \in \bZ_-$. Then we have 
\[
  {\bf P}\bigg\{ \sum_{k = 0}^{\fl{nx} - N} \omega_{k}< nt \log(|x| - 2 N n^{-1}) \bigg\} \le   {\bf P}\bigg\{ \sum_{k = 0}^{\fl{nx} - N} \omega_{k}< (nx - N)(1 - \e)  \bigg\} \le e^{-(nx - N)C(\e)},
\]
and in particular 
\[
 {\bf P}\big\{ \exists N < \fl{nr_1}: \xi^{n, \bf N}_{\lfloor nx \rfloor - N}(nt) > 0 \big\} \le \sum_{-N = \fl{nr_1}}^{\infty}  e^{-(nx - N)C(\e)} = C_{\e} e^{-(nx - nr_1)C(\e)}.
\]
Then by the Borel-Cantelli lemma we can find a random $n_1> n_0$ after which, for all $n>n_1$, we have 
\[
\xi^{n, \bf k}_{\lfloor nx \rfloor - k}(nt) = 0, \text{ for all } k < \fl{nr_1}.
\]

The second statement follows by a symmetric argument. Note that initially $\xi^{n, \lfloor \bf nr_{2} \rfloor}_{\lfloor nx \rfloor - \lfloor nr_{2} \rfloor}(0) = \lfloor nr_{2} \rfloor - \lfloor nx \rfloor$ for $\lfloor nr_{2} \rfloor > \lfloor nx \rfloor$, and the first jump time is given by
\begin{equation*}
    \tau_{2}^{n} = \inf \{s \in \mathbb{R}_{+} : \xi^{n, \lfloor \bf nr_{2} \rfloor}_{\lfloor nx \rfloor - \lfloor nr_{2} \rfloor}(s) > \lfloor nr_{2} \rfloor - \lfloor nx \rfloor\}.
\end{equation*}
We again want to bound this by considering the maximum speed of the particles. We can say that $\tau^{n}_{2}$ is bounded by a sequence of $\lfloor nr_{2} \rfloor - \lfloor nx \rfloor$ i.i.d.~exponential random variables of rate $M^{\bf 0}_{x, r_2}$, and the proof follows in a similar manner. 
\end{proof}

\begin{proof}[Proof of Proposition \ref{prop: Limit for z}]

We start by rewriting \eqref{Envelope Property} with the correct scaling
\begin{equation}
    n^{-1}z^{n}_{\lfloor nx \rfloor}(nt) = \sup_{q \in \mathbb{R}} \left\{ n^{-1}z^{n}_{\lfloor nq \rfloor}(0) - n^{-1}\xi^{n,\lfloor nq \rfloor}_{\lfloor nx \rfloor - \lfloor nq \rfloor}(nt) \right\}.
\end{equation} 
We want to take limits as $n\to \infty$ inside the supremum 
The candidate limit is the function 
\begin{equation}
   v(x,t) = \sup_{q \in \mathbb{R}} \left\{ v_{0}(q) - g^{q}(x-q,t)  \right\}, 
\end{equation}
as the first term in braces  is the limit of $ n^{-1}z^{n}_{\lfloor nq \rfloor}(0)$ from the assumptions. The limit of the second term can be argued as follows: Fix any $y \in \R$ and take $x' > x$. Then use $\lfloor n(x'-y) \rfloor \geq \lfloor nx \rfloor - \lfloor ny \rfloor \geq \lfloor n(x-y) \rfloor -1$ and the monotonicity of $\xi^{n, \lfloor \bf nq \rfloor}(t)$ to write 
 \begin{equation*}
    n^{-1}\xi^{n, \lfloor \bf nq \rfloor}_{\lfloor n(x'-y) \rfloor}(nt) \leq n^{-1}\xi^{n, \lfloor \bf nq \rfloor}_{\lfloor nx \rfloor - \lfloor ny \rfloor}(nt) \leq n^{-1}\xi^{n, \lfloor \bf nq \rfloor}_{\lfloor n(x-y) \rfloor}(nt) +n^{-1}.
\end{equation*}
Let $n\to \infty$ and apply Proposition \ref{prop: Limit for xi} to obtain limits for the left- and right-hand sides of the inequality above. Then let $y = q$ and after that use the continuity of $g^{q}(x,t)$ to take $x' \to x$ to obtain 
\begin{equation}
    \lim_{n \to \infty}n^{-1}\xi^{n, \lfloor \bf nq \rfloor}_{\lfloor nx \rfloor - \lfloor nq \rfloor}(nt) = g^{q}(x-q,t).
\end{equation}
We want to show that, for all $\varepsilon > 0$, we have the limit

\begin{equation} \label{v tilde limit}
    \lim_{n \to \infty}\mathbf{P} \left\{\left|n^{-1}z^{n}_{\lfloor nx \rfloor}(nt) - v(x,t) \right| > \varepsilon \right\} = 0.
\end{equation}
This gives us the convergence in probability. We will prove \eqref{v tilde limit} but the arguments go through for the almost sure convergence
\begin{equation*}
    \lim_{n \to \infty}n^{-1}z^{n}_{\lfloor nx \rfloor}(nt) = v(x,t).
\end{equation*}
in a similar way. 

For \eqref{v tilde limit} we bound separately  each of 
\[
\mathbf{P} \left\{n^{-1}z^{n}_{\lfloor nx \rfloor}(nt) - v(x,t) > \varepsilon \right\}, \quad \text{ and } \quad  \mathbf{P} \left\{n^{-1}z^{n}_{\lfloor nx \rfloor}(nt) - v(x,t)  < -\varepsilon \right\}.
\] 
For the left tail, given an $\varepsilon > 0$, pick $q = q(\e) \in \mathbb{R}$ such that
\begin{equation*}
    v_{0}(q) - g^{q}(x-q,t) \geq v(x,t) - \varepsilon/2.
\end{equation*}
In particular, we also have that for $n$ large enough and the same $q$, Lemma \ref{lemma: Envelope property} gives
\[
z^{n}_{\lfloor nx \rfloor}(nt) \ge z^{n}_{\lfloor nq \rfloor}(0) - \xi^{n, \lfloor \bf nq \rfloor}_{\lfloor nx \rfloor - \lfloor nq \rfloor}(nt).
\]
Then 
\begin{align*} 
%1 -& {\bf P} \left\{n^{-1}z^{n}_{\lfloor nx \rfloor}(nt) - \tilde{v}(x,t)  \le - \varepsilon \right\} = 
{\bf P}  \left\{n^{-1}z^{n}_{\lfloor nx \rfloor}(nt)< v(x,t)  - \varepsilon \right\} 
&\le {\bf P}  \left\{n^{-1}(z^{n}_{\lfloor nq \rfloor}(0) - \xi^{n, \lfloor \bf nq \rfloor}_{\lfloor nx \rfloor - \lfloor nq \rfloor}(nt))<  v_{0}(q) - g^{q}(x-q,t) - \e/2\right\}\\
&\le   {\bf P}  \left\{ n^{-1}z^{n}_{\lfloor nq \rfloor}(0) \leq v_{0}(q) -\e/4 \right\} \\
&\phantom{xxxxxx}+  {\bf P}  \left\{ n^{-1}\xi^{n, \lfloor \bf nq \rfloor}_{\lfloor nx \rfloor - \lfloor nq \rfloor}(nt) \geq g^{q}(x-q,t) + \varepsilon/4 \right\}.
\end{align*}
Both of the last two probabilities are near 0 for $n$ large enough (see \eqref{eq:zic} and \eqref{auxilliary function limit}), which gives us  
\[
\lim_{n\to \infty} {\bf P} \left\{n^{-1}z^{n}_{\lfloor nx \rfloor}(nt) - v(x,t)  \le - \varepsilon \right\} = 0.
\]

For the right tail we need to utilise Lemma \ref{lemma: truncation for q}. For $x$ and $t$ fixed, choose  $r_{1}, r_{2}$ as in the statement of the lemma. We will then substitute $z^{n}$ with the random variable $\zeta^n_{x,t}(r_1, r_2)$ that agrees with $z^n$ on a high probability event by equation \eqref{eq:qtr} and show
\begin{equation} \label{zeta lower bound}
    \lim_{n \to \infty}\mathbf{P} \left\{ n^{-1}   \zeta^{n}_{x,t}(r_{1},r_{2}) > v(x,t) + \varepsilon \right\} = 0.
\end{equation}

To do this, consider the partition $\mathscr P = \{a_i \}_{i \in [m]}$ of the interval $r_1 = a_0 < a_1 < \ldots < a_m = r_2$ so that 
Lemmas \ref{lem:gjointc} and \ref{lem:jointxicontrol} can be applied in each subinterval, with an error of at most $\e/3$ for each application below. %\note{say better}
Then we can write 
\begin{align*}
      \zeta^{n}_{x,t}(r_{1},r_{2}) &= \max_{0 \leq j < m} \max_{\fl{na_{j}} \leq k \leq \fl{na_{j+1}}} \left\{ z^{n}_{k}(0) - \xi^{n,\bf k}_{\fl{nx} - k}(nt) \right\} \\
    & \le \max_{0 \leq j < m} \max_{\fl{na_{j}} \leq k \leq \fl{na_{j+1}}} \left\{ z^{n}_{\fl{na_{j+1}}}(0) - \xi^{n,\bf k}_{\fl{nx} - \fl{na_j}}(nt) \right\} \\
     & \le  \max_{0 \leq j < m} \left\{ z^{n}_{\fl{na_{j+1}}}(0) - \xi^{n, \fl{ \bf na_j}}_{\fl{nx} - \fl{na_j}}(nt) \right\}  \\
    &\phantom{xxxxxxx}+  \max_{0 \leq j < m} \left\{ \xi^{n, \fl{\bf na_j}}_{\fl{nx} - \fl{na_j}}(nt) -  \min_{\fl{na_{j}} \leq k \leq \fl{na_{j+1}}} \xi^{n,\bf k}_{\fl{nx} - \fl{na_j}}(nt) \right\}.
\end{align*}
From Lemma \ref{lem:gjointc} and Lemma \ref{lem:jointxicontrol}
\[
\varliminf_{n\to \infty} \frac{1}{n}  \min_{\fl{na_{j}} \leq k \leq \fl{na_{j+1}}} \xi^{n,\bf k}_{\fl{nx} - \fl{na_j}}(nt) \ge \inf_{q \in [a_j, a_{j+1}]}  g^{q}(x-a_j, t) -\e > g^{a_j}(x-a_j, t) -2\e/3 .
\]
With this in place we have 
\begin{align*}
\varlimsup_{n\to \infty} \frac{1}{n} \zeta^{n}_{x,t}(r_{1},r_{2}) & \le  \max_{0 \leq j < m}  \left\{ v_{0}(a_{j+1}) - g^{a_{j}}(x-a_{j},t) \right\}+ 2\e/3\\
&\le  \max_{0 \leq j < m} \left\{ v_{0}(a_{j+1}) - g^{a_{j+1}}(x-a_{j+1},t) \right\} +\e\\
&\le  v(x,t) + \varepsilon.
\end{align*}
This proves \eqref{zeta lower bound}, and consequently \eqref{v tilde limit}.
\end{proof}

\begin{corollary} \label{cor: restriction to compact interval} 
For any $x \in \bR, t \in \R_+$, there exist deterministic $r_L(x,t) < 0 < r_{U}(x,t)$ so that almost surely
\begin{equation}\label{eq:goodv}
    v(x,t) =\lim_{n \to \infty}n^{-1}z^{n}_{\fl{nx}}(nt) = \max_{q \in [r_L, r_U]} \left\{v_{0}(q) - g^{q}(x-q,t) \right\} = \sup_{q \in \mathbb{R}} \left\{v_{0}(q) - g^{q}(x-q,t) \right\}.
\end{equation}
\end{corollary}

\begin{proof}
This follows from the fact that $\zeta$ agrees with $z^n$ for $q$ in a compact interval, so the supremum can be obtained by checking only the $q \in [r_L, r_U]$ and the fact that it is a maximum follows from the continuity of $v_0$ and $g^q$ in the $q$-argument.  
\end{proof}

\begin{lemma}\label{lem:maximiser-off-contact}
Fix $(x,t)\in\mathbb R\times\mathbb R_+$ and define
\[
F_{x,t}(q):=v_0(q)-g^q(x-q,t),
\qquad q\in [r_L(x,t),r_U(x,t)],
\]
where $r_L(x,t),r_U(x,t)$ are as in Corollary \ref{cor: restriction to compact interval}. Define the two contact sets
\[
A_0(x,t):=\{q\in [r_L,r_U]: q\le x,\ g^q(x-q,t)=0\},
\]
and
\[
A_1(x,t):=\{q\in [r_L,r_U]: q\ge x,\ g^q(x-q,t)=q-x\}.
\]
Then
\[
\max_{q\in[r_L,r_U]}F_{x,t}(q)
=
\max_{q\in [r_L,r_U]\setminus (\operatorname{int}A_0(x,t)\cup \operatorname{int}A_1(x,t))}F_{x,t}(q).
\]
In particular, there exists a maximiser $q^*$ such that either
\[
g^{q^*}(x-q^*,t)>\max\{0,q^*-x\},
\]
or $q^*$ is an endpoint of a connected component of $A_0(x,t)$ or of $A_1(x,t)$.
\end{lemma}

\begin{proof}
By Corollary \ref{cor: restriction to compact interval}, the supremum in the variational formula is attained on the compact interval $[r_L,r_U]$.
Let $I$ be a connected component of $A_0(x,t)$. For $q\in I$ we have
\[
F_{x,t}(q)=v_0(q).
\]
% Since
% $0\le \rho_0\le 1$, equation~\eqref{eq:in-cur} implies that for $a<b$,
% \[
% v_0(b)-v_0(a)=\int_a^b \rho_0(s)\,ds\ge 0.
% \]
% Hence 
We have that $v_0$ is nondecreasing, so the maximum of $F_{x,t}$ on $I$ is attained at the right endpoint of $I$.

Let $J$ be a connected component of $A_1(x,t)$. For $q\in J$ we have
\[
F_{x,t}(q)=v_0(q)-(q-x)=x+(v_0(q)-q).
\]
Since $v_0$ is Lipschitz-1, for $a<b$,
\[
v_0(b)-v_0(a)\le b-a,
\]
which gives
\[
(v_0(b)-b)-(v_0(a)-a)\le 0.
\]
Thus $q\mapsto v_0(q)-q$ is nonincreasing, so the maximum of $F_{x,t}$ on $J$ is attained at the left endpoint of $J$.

Therefore no maximiser is needed in the interior of either contact set, and the claimed reduction follows.
\end{proof}

\begin{corollary}\label{cor:maximiser-levelcurve}
Fix $(x,t)\in\mathbb R\times\mathbb R_+$. There exists a maximiser $q^*\in [r_L(x,t),r_U(x,t)]$ in Corollary \ref{cor: restriction to compact interval} such that
\begin{equation}\label{eq:maximiser-levelcurve}
\Gamma_c^{q^*}\bigl(x-q^*+g^{q^*}(x-q^*,t),g^{q^*}(x-q^*,t)\bigr)=t.
\end{equation}
\end{corollary}

\begin{proof}
Choose $q^*$ as in Lemma~\ref{lem:maximiser-off-contact}. If
\[
g^{q^*}(x-q^*,t)>\max\{0,q^*-x\},
\]
then \eqref{eq:maximiser-levelcurve} follows immediately from Lemma~\ref{lem:interior-level-curve} applied with $x$ replaced by $x-q^*$.

It remains to consider the case where $q^*$ is an endpoint of a connected component of $A_0(x,t)$ or $A_1(x,t)$. Then there exists a sequence $q_n\to q^*$ with $q_n\in [r_L,r_U]$ lying outside the corresponding contact component. For each $n$,
\[
g^{q_n}(x-q_n,t)>\max\{0,q_n-x\},
\]
so Lemma~\ref{lem:interior-level-curve} gives
\[
\Gamma_c^{q_n}\bigl(x-q_n+g^{q_n}(x-q_n,t),g^{q_n}(x-q_n,t)\bigr)=t.
\]
Now use joint continuity of $g^q(x,t)$ from Lemma~3.3 and continuity in $q$ of $\Gamma_c^q$ to pass to the limit $n\to\infty$. This yields \eqref{eq:maximiser-levelcurve}.
\end{proof}

\begin{proof}[Proof of Theorem \ref{thm:LLN}]
The first statement follows from Proposition \ref{prop: Limit for z}, which also gives a variational characterisation for $ v$.  Lipschitz-1 continuity of ${v}$ with respect to the $x$ variable follows from the exclusion rule which implies
\[
| z^n_{\fl{nx}}(nt) - z^n_{\fl{ny}}(nt) | \le |nx - ny|.
\]
In turn, this implies that $\tilde{v}$ is almost everywhere differentiable with respect to $x$ and that $\rho(x,t)$ is well defined almost everywhere. Since $\rho$ is well defined almost everywhere, the last statement of the theorem then follows directly from the first (see remark \ref{remark: limit for z}).

For the remaining Lipschitz and differentiability statements, it suffices to show local Lipschitz continuity with respect to the $t$ variable. For that, we will show that for any $t, t'$ in a fixed neighbourhood $U$ of some $s_0 >0$ we have that there exists a constant $M =M(x, s_0, U)$ such that
\[ 
| v(x,t) - v(x, t')| \le M|t - t'|.
\] 
Fix $s_0>0$ and a radius $r>0$. 

Let, without loss of generality,   $t > t' > 0$ so that $|t-s_0| < r$ and $|t' - s_0| < r$.  Let  $x \in \mathbb{R}$. From Corollary \ref{cor: restriction to compact interval} we can find a $q_{0} \in \mathbb{R}$ that satisfies \eqref{eq:maximiser-levelcurve} and 
\begin{equation}
    v(x,t') = v_{0}(q_0) - g^{q_0}(x-q_0,t').
\end{equation}
From this we can say
\[
   0\le  v(x,t') - v(x,t) = v_{0}(q_{0}) - g^{q_{0}}(x-q_{0},t') - \max_{q' \in \mathbb{R}}\{v_{0}(q') - g^{q'}(x-q',t)\},
\]
and we can set $q' = q_{0}$ in the supremum to obtain
\[
    v(x,t') - v(x,t) \leq  g^{q_{0}}(x-q_{0},t) - g^{q_{0}}(x-q_{0},t'). 
\]
Let  $g^{q_{0}}(x-q_{0},t) = y_{0}$, $g^{q_{0}}(x-q_{0},t') = y_{0}'$. 
We also have  $\Gamma^{q_0}_c(x-q_0+y_0, y_0) = t$,  $\Gamma^{q_0}_c(x-q_0+y_0', y_0') \ge t'$, so we can write
\begin{align*}
t - t' &\ge \Gamma^{q_0}_c(x-q_0+y_0, y_0) - \Gamma^{q_0}_c(x-q_0+y_0', y_0') \\
 & \ge  \Gamma^{q_0}_c((x-q_0+y_0', y_0'),(x-q_0+y_0, y_0)), \quad \text{by super-additivity}\\
 & \ge \int_0^1 \frac{ \gamma(\bf x'(s))}{ c(x_1(s) + q_0-v_0(q_0), x_2(s) - v_0(q_0)) } \, ds, \\
 &\phantom{xxxxxxx}\quad \text{where } {\bf x}(s) = (x- q_0 + y_0'+ (y_0- y_0')s, y_0'   + (y_0 - y_0')s)\\
 &\ge \frac{4}{c_{\max}}(y_0 - y_0') \ge  \frac{4}{c_{\max}}( v(x,t') - v(x,t) ).
\end{align*}
The constant  $c_{\max}$ is the maximum value of $c(x_1(s) + q_0 -v_0(q_0), x_2(s) - v_0(q_0))$ in the square with lower left corner $(x-q_0+y_0', y_0')$ and upper right corner $(x-q_0+y_0, y_0)$ so it is a function of $q_0, x,t,t'$. 

Define 
\[y_{\min} = g^{q_{0}}(x-q_{0},s-r), \quad y_{\max} = g^{q_{0}}(x-q_{0}, s+r). \]
Since $t > t'$ it follows that 
\[ y_{\max} \ge y_{0} \geq y_{0}' \geq y_{\min}.\]  
Therefore we can enlarge $c_{\max}$ further by assuming it is the supremum of the values in the square 
$$R_{(x-q_0+y_{\min}, y_{\min}), (x-q_0+y_{\max}, y_{\max})}. $$ This way, $c_{\max}$ only depends on the $r$-neighbourhood of $s$ rather than the $t,t'$ and it can be used for the definition of locally Lipschitz. 
\end{proof}

\section{Solutions to discontinuous Hamilton--Jacobi equations} 
\label{sec:visco}

In this section we will prove all the claimed PDE connections. It is more convenient to work with wedge last passage times and their a.s.~ limits

If
$(x_1,y_1)$ and $(x_2,y_2)$ are absolute wedge points, we write the absolute shape function $\widetilde\Gamma_{\tilde c}$ as 
\[
\widetilde\Gamma_{\tilde c}((x_1,y_1),(x_2,y_2)) := \Gamma_c((x_1+y_1, y_1), (x_2+y_2, y_2)),
\]
for the wedge passage time in the global $\tilde c$-environment from
$(x_1,y_1)$ to $(x_2,y_2)$. 

Often we use the $q$-{\bf rooted relative shape function} $\widetilde\Gamma_{\tilde c}^{\,q}(r,s)$, where in absolute wedge coordinates becomes
\[
\widetilde\Gamma_{\tilde c}^{\,q}(r,s)
=
\widetilde\Gamma_{\tilde c}
\bigl((q,-v_0(q)),(q+r,-v_0(q)+s)\bigr) = \Gamma^q_c(r+s,s).
\]
With notation $\widetilde\Gamma_{\tilde c}^{\,q}(r,s)$ the argument $(r,s)$ corresponds to displacement in the relative coordinate system rather than the terminal point. 

Moreover, in the case of the $q-$rooted relative shape function, it may be that the starting point is not the root $(q, -v_0(q))$, and notation
$
\widetilde\Gamma_{\widetilde c}^{\,q}((a,b),(r,s)) 
$ is used. This is again in relative coordinates, so it should read as 
\[
\widetilde\Gamma_{\tilde c}^{\,q}((a,b),(r,s)) = \widetilde\Gamma_{\tilde c}((q,-v_0(q))+(a,b), (q,-v_0(q))+(r,s)) = \widetilde\Gamma_{T_{(q,-v_0(q)}\tilde c}((a,b),(r,s)),
\]
and $(a,b) \le_{\mathbb W} (r,s)$.
The most common use of such notation is on lifted level-curve points: Their absolute coordinates will be
$
\bigl(x,g^q(x-q,t)-v_0(q)\bigr)
$
while their relative $q$-coordinates
$
\bigl(x-q,g^q(x-q,t)\bigr).
$

Consequently, if $x$ denotes an absolute spatial coordinate, then the relevant
relative wedge coordinate is $r=x-q$, and
\[
\widetilde\Gamma_{\tilde c}^{\,q}
\bigl(x-q,g^q(x-q,t)\bigr)\ge t.
\]
If \(g^q(x-q,t)>\beta(x-q)\), then equality holds.

\subsection{Strong solution at differentiability points}

In this section we establish further analytical properties for the level curves $g^q(x-q,t)$ that will then be used throughout in the proof for the existence of solutions for the PDE. 

\begin{lemma}\label{lem:recbelowlc}
Recall \eqref{g^q definition til}. For any point $(x, g^z(x-z,t)-v_0(z))$ on the level curve that satisfies \eqref{eq:maximiser-levelcurve} we have that the (potentially degenerate) parallelogram %\note{check parallelogram}
\[A_{z,x,t} = \{ (u,v): z - v_0(z) \le u+v \le x +g^z(x-z,t)-v_0(z), -v_0(z) \le v \le g^z(x-z,t)-v_0(z)  \} \]
touches the level curve only at point $(x, g^z(x-z,t)-v_0(z))$ and otherwise is below it.
\end{lemma}

\begin{proof} Assume by way of contradiction that there is a point $(u_0,r_0)$ in $A_{z,x,t}$ so that $r_0 = g^z(u_0-z, t) -v_0(z)$. Then we have that $\widetilde \Gamma_{\tilde c}^z(u_0-z,g^z(u_0-z, t))=t$, while (recall the wedge ordering \eqref{eq:wedgeorder}) $(u_0,r_0) \le_{\mathbb W} (x, g^z(x-z,t)-v_0(z))$.
But then, by super-additivity
\begin{align*}
t &= \widetilde \Gamma^z_{\tilde c}(x-z, g^z(x-z,t)) \ge \widetilde \Gamma^z_{\tilde c}(u_0-z, g^z(u_0-z,t)) + \widetilde \Gamma^z_{\tilde c}((u_0-z, g^z(u_0-z,t)),(x-z, g^z(x-z,t))\\
&= t + \widetilde \Gamma^z_{\tilde c}((u_0-z, g^z(u_0-z,t)),(x-z, g^z(x-z,t)).
\end{align*}
The last term is strictly positive if and only if $(u_0,r_0) <_{\mathbb W} (x, g^z(x-z,t)-v_0(z))$, which leads to the desired contradiction unless $(u_0,r_0)= (x, g^z(x-z,t)-v_0(z))$.

Therefore, the curve $g^{z}(\cdot-z, t)-v_0(z)$ only touches $A_{z,x,t}$  at the north east corner. 
\end{proof}

An immediate corollary is that for any fixed $z\in \R$ and any $(u_0, r_0) \in A^{\circ}_{z,x,t}$,
\be \label{eq:tinyhcurve}
\widetilde\Gamma^z_{\tilde c}(u_0-z, r_0+v_0(z)) < t.
\ee

\begin{lemma}\label{lem:psitaylor}Recall $\psi(x)$ from \eqref{Legendre dual relation}, 
and assume $\alpha \in (-1,1)$ and $0<\delta <\frac{(1+\alpha)^2}{2(1-\alpha)}$. Then we have the bound 
\be
\gamma(\alpha + \psi(\alpha + \delta), \psi(\alpha + \delta)) \le 1 - C_{\alpha}\delta +o(\delta), \text{ for } \alpha \in (-1,1),
\ee
with $0 < C_{\alpha} < \infty$.
\end{lemma}

\begin{proof}
First consider a Taylor expansion around $\alpha > -1$ for $\psi(\alpha + \delta)$. We have 
\be
\psi(\alpha + \delta) = \psi(\alpha) + \psi'(\alpha)\delta + o(\delta) = \psi(\alpha) - \frac{(1-\alpha)}{2}\delta +o(\delta)
\ee
By the inequality $\sqrt{x-\e} < \sqrt{x} - \frac{\e}{2\sqrt{x}}$ and the monotonicity of $\gamma$ we have 
\begin{align*}
\gamma(\alpha + \psi(\alpha + \delta), \psi(\alpha + \delta))
&= \gamma\left(\alpha + \psi(\alpha) - \frac{(1-\alpha)}{2}\delta +o(\delta), \psi(\alpha)- \frac{(1-\alpha)}{2}\delta +o(\delta)\right)\\
&=\left( \sqrt{\alpha + \psi(\alpha) - \frac{(1-\alpha)}{2}\delta +o(\delta)} + \sqrt{\psi(\alpha)- \frac{(1-\alpha)}{2}\delta +o(\delta)} \right)^2\\
&\le \left( \sqrt{\alpha + \psi(\alpha)} - \frac{(1-\alpha)}{4 \sqrt{\alpha + \psi(\alpha)}}(\delta +o(\delta)) + \sqrt{\psi(\alpha)}- \frac{(1-\alpha)}{4\sqrt{\psi(\alpha)}}(\delta +o(\delta)) \right)^2\\
&\le 1 - \frac{(1-\alpha)}{ \sqrt{\alpha + \psi(\alpha)}}\delta +o(\delta) = 1- \frac{2(1-\alpha)}{1+\alpha} \delta + o(\delta). \qedhere
\end{align*}
\end{proof}

The first thing we need to argue is that the continuity points $(x,t)$ of $\tilde c(x, -v(x,t))$ have full measure in any compact set, and are dense in $\R^2$.

\begin{lemma}[Continuity points of $\tilde c(x,-v(x,t))$]
\label{lem:continuity-points-full-measure}
Under Assumption \ref{ass:c0strong}
there exists a full measure and dense set $D\subset \mathbb R\times(0,\infty)$ such that the composition
\[
(x,t)\longmapsto \tilde c(x,-v(x,t))
\]
is continuous at every point $(x_0,t_0)\in D$.
\end{lemma}

\begin{proof}
Fix a compact rectangle
\[
K=[a,b]\times[\tau_1,\tau_2]\subset \mathbb R\times(0,\infty).
\]

Let
\[
\Phi(x,t):=(x,-v(x,t)).
\]
Since $v$ is continuous, $\Phi(K)$ is compact. By the local finiteness of the
discontinuity curves of $\widetilde c$, the discontinuity set of $\widetilde c$
inside $\Phi(K)$ is contained in the union of finitely many non-vertical graph
segments
\[
\Sigma_j=\{(x,h_j(x)):x\in I_j\}, \qquad j=1,\ldots,N,
\]
finitely many vertical segments
\[
V_\ell=\{x_\ell\}\times J_\ell,\qquad \ell=1,\ldots,M,
\]
and a countable exceptional set $P$ consisting of terminal and intersection
points.

For each fixed $t\in[\tau_1,\tau_2]$, define
\[
f_t(x):=-v(x,t), \qquad x\in[a,b].
\]
Since $v$ is jointly locally Lipschitz, each $f_t$ is Lipschitz on $[a,b]$, hence absolutely continuous. Moreover, by Theorem \ref{thm:LLN},
$
f_t'(x)\in[-1,0]
\text{ for a.e. }x\in[a,b].
$

We first deal with the finitely many discontinuity curves that are graphs. For each $j=1,\dots,N$ and each fixed $t$, define
\[
E_{j,t}
=
\{x\in[a,b]:(x,f_t(x))\in\Sigma_j\}
=
\{x\in[a,b]: f_t(x)=h_j(x)\}.
\]
We claim that
$|E_{j,t}|=0$
for every fixed $t\in[\tau_1,\tau_2].$
Assume by contradiction that $|E_{j,t}|>0$ for some $j$ and some $t$. Consider
\[
\phi(x):=f_t(x)-h_j(x).
\]
Since both $f_t$ and $h_j$ are absolutely continuous, so is $\phi$. Moreover, $\phi=0$ on the set $E_{j,t}$ of positive measure.

Now let $x_0$ be a density point of $E_{j,t}$ such that $\phi'(x_0)$ exists. Since $\phi$ is differentiable at $x_0$ and there exists a sequence $x_n\in E_{j,t}$ with $x_n\to x_0$, we have
\[
\phi'(x_0)
=
\lim_{n\to\infty}\frac{\phi(x_n)-\phi(x_0)}{x_n-x_0}
=
0.
\]
Thus $\phi'(x)=0$ for a.e. $x\in E_{j,t}$. Equivalently,
\[
f_t'(x)=h_j'(x)
\qquad\text{for a.e. }x\in E_{j,t}.
\]
But $f_t'(x)\in[-1,0]$ a.e., whereas by assumption
$
h_j'(x)\notin[-1,0]$
for a.e. $x\in I_j.
$
This is impossible on a set of positive measure. Therefore $|E_{j,t}|=0$.

Now define the subset of bad points in $K$ associated with the curves $\Sigma_j$ by
\[
E_j
=
\{(x,t)\in K : (x,-v(x,t))\in \Sigma_j\}
=
\{(x,t)\in K : -v(x,t)=h_j(x)\}.
\]
Since the section of $E_j$ at time $t$ is exactly $E_{j,t}$, Fubini's theorem gives
\[
|E_j|
=
\int_{\tau_1}^{\tau_2} |E_{j,t}|\,dt
=
0.
\]
Since there are only finitely many curves on $K$,
$
\left|\bigcup_{j=1}^N E_j\right|=0.$

We next treat the vertical discontinuity segments. For each $\ell=1,\ldots,M$, set
\[
E^{\mathrm{vert}}_\ell
:=
\{(x,t)\in K:\Phi(x,t)\in V_\ell\}.
\]
Since $V_\ell\subset \{x_\ell\}\times\mathbb R$, we have
$
E^{\mathrm{vert}}_\ell
\subset
\{x_\ell\}\times[\tau_1,\tau_2],
$
and therefore
$
|E^{\mathrm{vert}}_\ell|=0.
$
As there are only finitely many vertical segments meeting $\Phi(K)$, it follows that
\[
\left|\bigcup_{\ell=1}^M E^{\mathrm{vert}}_\ell\right|=0.
\]

Next we deal with the countable exceptional set $P=\{p_m\}_{m\ge 1}$, where $p_m=(x_m,y_m)$. Define
\[
E_P
=
\{(x,t)\in K : (x,-v(x,t))\in P\}.
\]
For each $m$,
\[
\{(x,t)\in K:(x,-v(x,t))=p_m\}
\subseteq
\{x_m\}\times[\tau_1,\tau_2],
\]
which is a vertical line segment and therefore has two-dimensional Lebesgue measure zero. Since $P$ is countable,
$
|E_P|=0.
$

Set
\[
E
:=
\left(\bigcup_{j=1}^N E_j\right)
\cup
\left(\bigcup_{\ell=1}^M E^{\mathrm{vert}}_\ell\right)
\cup
E_P .
\]
Then $|E|=0$.

If $(x,t)\in K\setminus E$, then
\[
(x,-v(x,t))\notin \Big(\bigcup_{j=1}^N \Sigma_j\Big)\cup P,
\]
so $\tilde c$ is continuous at the point $(x,-v(x,t))$. Since $(x,t)\mapsto (x,-v(x,t))$ is continuous, it follows that the composition
$
(x,t)\longmapsto \tilde c(x,-v(x,t))$
is continuous at $(x,t)$.

Hence every point of $K\setminus E$ is a continuity point of the composition, and the set of bad points in $K$ has measure zero. 
Finally, since $K$ was an arbitrary compact rectangle in $\mathbb R\times(0,\infty)$, the set of bad points has
measure zero in every such $K$. Hence its complement is dense in $\mathbb R\times(0,\infty)$ and has full
measure there.
\end{proof}

\begin{proof}[Proof of Theorem \ref{thm: existence}.]
Let $(x_0,t_0)$ be a point of differentiability of $v(x,t)$ so that $\tilde c(x, -v(x,t) )$ is also continuous at $(x_0,t_0)$. 

Find a $z_0$ so that $v(x_0,t_0) = v_0(z_0)-g^{z_0}(x_0-z_0, t_0)$ so that it also satisfies \eqref{eq:maximiser-levelcurve}, and consider the point $(x_0, g^{z_0}(x_0-z_0, t_0) - v_0(z_0)) = (x_0, -v(x_0,t_0))$. 

By Lemma \ref{lem:recbelowlc} the parallelogram $A_{z_0, x_0,t_0}$ is below the curve $g^{z_0}(\cdot-z_0, t_0)-v_0(z_0)$ except when they intersect at point $(x_0, -v(x_0,t_0))$. As such, for any $h>0$  we can find a point $A = (x_A, -v(x_0,t_0))$ on the north boundary and a point $B = (x_B, y_B)$ on the east boundary of $A_{z_0, x_0,t_0}$ so that both $A$ and $B$ belong on the curve $g^{z_0}(\cdot-z_0, t_0-h )-v_0(z_0)$ by equation \eqref{eq:tinyhcurve}.
\begin{figure}[ht]

\tikzset{every picture/.style={line width=0.75pt}} %set default line width to 0.75pt        

\begin{tikzpicture}[x=0.75pt,y=0.75pt,yscale=-1,xscale=1]
%uncomment if require: \path (0,861); %set diagram left start at 0, and has height of 861

%Straight Lines [id:da3491239725745282] 
\draw  [dash pattern={on 4.5pt off 4.5pt}]  (364.25,657) -- (366.5,747) ;
%Curve Lines [id:da7858377845599059] 
\draw [line width=2.25]    (238.25,747) .. controls (203.5,631) and (455.5,694) .. (411.25,597) ;
%Straight Lines [id:da07984451360043854] 
\draw    (105.5,588) -- (244,753) ;
%Straight Lines [id:da8639191035208056] 
\draw    (232,747) -- (535.5,747) ;
%Straight Lines [id:da6955730914788979] 
\draw  [dash pattern={on 4.5pt off 4.5pt}]  (109,597) -- (453.5,597) ;
%Straight Lines [id:da8281511936800956] 
\draw  [dash pattern={on 4.5pt off 4.5pt}]  (397,582) -- (535.5,747) ;
%Curve Lines [id:da45538215803657567] 
\draw [color={rgb, 255:red, 208; green, 2; blue, 27 }  ,draw opacity=1 ][line width=1.5]    (357,541) .. controls (359.5,602) and (455.5,595) .. (538.5,661) ;
%Curve Lines [id:da7178147316156217] 
\draw [color={rgb, 255:red, 144; green, 19; blue, 254 }  ,draw opacity=1 ][line width=2.25]    (198,569) .. controls (200.5,630) and (465.5,670) .. (548.5,736) ;
%Shape: Circle [id:dp44320123720132143] 
\draw  [fill={rgb, 255:red, 245; green, 166; blue, 35 }  ,fill opacity=1 ] (210.75,598.5) .. controls (210.75,595.05) and (213.55,592.25) .. (217,592.25) .. controls (220.45,592.25) and (223.25,595.05) .. (223.25,598.5) .. controls (223.25,601.95) and (220.45,604.75) .. (217,604.75) .. controls (213.55,604.75) and (210.75,601.95) .. (210.75,598.5) -- cycle ;
%Shape: Circle [id:dp1323500259489765] 
\draw  [fill={rgb, 255:red, 245; green, 166; blue, 35 }  ,fill opacity=1 ] (496.75,710.5) .. controls (496.75,707.05) and (499.55,704.25) .. (503,704.25) .. controls (506.45,704.25) and (509.25,707.05) .. (509.25,710.5) .. controls (509.25,713.95) and (506.45,716.75) .. (503,716.75) .. controls (499.55,716.75) and (496.75,713.95) .. (496.75,710.5) -- cycle ;
%Shape: Circle [id:dp6713719285657432] 
\draw  [fill={rgb, 255:red, 245; green, 166; blue, 35 }  ,fill opacity=1 ] (232,747) .. controls (232,743.55) and (234.8,740.75) .. (238.25,740.75) .. controls (241.7,740.75) and (244.5,743.55) .. (244.5,747) .. controls (244.5,750.45) and (241.7,753.25) .. (238.25,753.25) .. controls (234.8,753.25) and (232,750.45) .. (232,747) -- cycle ;
%Shape: Circle [id:dp39349277947754213] 
\draw  [fill={rgb, 255:red, 245; green, 166; blue, 35 }  ,fill opacity=1 ] (405,597) .. controls (405,593.55) and (407.8,590.75) .. (411.25,590.75) .. controls (414.7,590.75) and (417.5,593.55) .. (417.5,597) .. controls (417.5,600.45) and (414.7,603.25) .. (411.25,603.25) .. controls (407.8,603.25) and (405,600.45) .. (405,597) -- cycle ;
%Shape: Circle [id:dp25390931843376296] 
\draw  [fill={rgb, 255:red, 245; green, 166; blue, 35 }  ,fill opacity=1 ] (358,657) .. controls (358,653.55) and (360.8,650.75) .. (364.25,650.75) .. controls (367.7,650.75) and (370.5,653.55) .. (370.5,657) .. controls (370.5,660.45) and (367.7,663.25) .. (364.25,663.25) .. controls (360.8,663.25) and (358,660.45) .. (358,657) -- cycle ;

% Text Node
\draw (340,539) node [anchor=north west][inner sep=0.75pt]    {$t$};
% Text Node
\draw (180,539) node [anchor=north west][inner sep=0.75pt]    {$t-h$};
% Text Node
\draw (223,570) node [anchor=north west][inner sep=0.75pt]    {$( x_{A} ,\ y_{A})$};
% Text Node
\draw (511,683) node [anchor=north west][inner sep=0.75pt]    {$( x_{B} ,\ y_{B})$};
% Text Node
\draw (418,567) node [anchor=north west][inner sep=0.75pt]    {$( x_{0} ,\ -v( x_{0} ,t_{0}))$};
% Text Node
\draw (180,757) node [anchor=north west][inner sep=0.75pt]    {$( z_{0} ,\ -v_{0}( z_{0}))$};
% Text Node
\draw (357,757) node [anchor=north west][inner sep=0.75pt]    {$x_{h} \ $};
\end{tikzpicture}
\caption{Schematic of the points at the first part of the proof of Theorem \ref{thm: existence}. The boldfaced increasing path is a near optimiser for the passage time from $(z_0, -v_0(z_0))$ up to the point $(x_0, -v(x_0,t_0))$ and it intersects the $t-h$ level curve at point $(x_h, y_h)$.}
\end{figure}
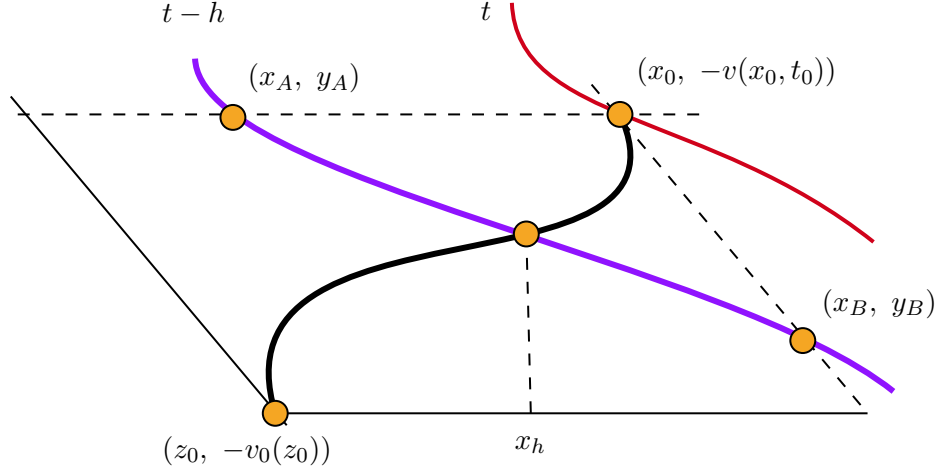

Consider the closed region $D= D_{x_0,t_0,h}$ bounded by the curves $g^{z_0}(\cdot-z_0, t_0-h )-v_0(z_0)$, $g^{z_0}(\cdot-{z_0}, t_0)-v_0(z_0)$ and $A_{z_0,x_0,t_0}$ and define $\e_h >0$ by 
\be\label{eq:eh}
\e_h = \max_{(x,y) \in D} \tilde c(x,y) - \min_{(x,y) \in D} \tilde c(x,y) : = \tilde c_D^{\rm high}(h) - \tilde c_D^{\rm low}(h), 
\ee
and note that as $h \to 0$, $\e_h \to 0$ by our assumptions on the point $(x_0,t_0)$ since 
\be \label{eq:localconv}
\lim_{h \to 0} \tilde c_D^{\rm low}(h) = \lim_{h \to 0} \tilde c_D^{\rm high}(h) = \tilde c(x_0, -v(x_0,t_0)).
\ee

By the superadditivity of the passage times we have that 
\begin{align*} 
t &= \widetilde \Gamma ^{z_0}_{\tilde c}(x_0-z_0, g^{z_0}(x_0-z_0,t_0))\\
&\ge \widetilde \Gamma ^{z_0}_{\tilde c}(x_A-z_0, g^{z_0}(x_A-z_0,t_0-h)) + \widetilde \Gamma ^{z_0}_{\tilde c}((x_A-z_0, g^{z_0}(x_A-z_0,t_0-h)),(x_0-z_0, g^{z_0}(x_0-z_0,t_0)) \\
&\ge t-h + \frac{1}{\tilde c_D^{\rm high}(h)}(x_0-x_A).
\end{align*}
This gives us that 
\be\label{eq:xislope1}
1 \ge \frac{x_0-x_A}{h \tilde c_D^{\rm high}(h)}.
\ee
Similarly, by using point $B$ we have 
\[
t\ge t-h + \frac{1}{\tilde c_D^{\rm high}(h)}(y_0-y_B) = t-h + \frac{1}{\tilde c_D^{\rm high}(h)}(x_B-x_0),
\]
giving that 
\be\label{eq:xislope2}
-1 \le \frac{x_0-x_B}{h \tilde c_D^{\rm high}(h)}.
\ee
For any $h>0$ consider a near optimizer $\bf w$ for 
$\widetilde \Gamma ^{z_0}_{\tilde c}(x_0-z_0, g^{z_0}(x_0-z_0,t_0))=t$ so that 
\[
I({\bf w}) := \int_0^1 \frac{\gamma(w_1'(s)+w_2'(s), w_2'(s))}{\tilde c (w_1(s)+z_0, w_2(s)-v_0(z_0))}\,ds \ge \widetilde \Gamma ^{z_0}_{\tilde c}(x_0-z_0, g^{z_0}(x_0-z_0,t_0)) - h\e_h =t- h\e_h.
\]
The path ${\bf w}$ will cross the curve $g^{z_0}(\cdot-z_0, t_0-h)-v_0(z_0)$
at some point $(x_h,y_h) = (x_h, g^z(x_h-z_0, t-h)-v_0(z_0) )$ inside $A_{z_0, x_0,t_0}$. 

Let $I_1, I_2$ denote the contributions of the path $\bf w$ before and after $(x_h,y_h)$. Then we can estimate 
\begin{align*}
t &= \widetilde \Gamma ^{z_0}_{\tilde c}(x_0-z_0, g^{z_0}(x_0-z_0,t_0)) \le I({\bf w}) +h\e_h=I_1 + I_2 + h\e_h\\
&\le\widetilde\Gamma ^{z_0}_{\tilde c}(x_h-z_0, g^{z_0}(x_h-z_0,t_0-h)) \\
&\phantom{xxxxxxx}+  \widetilde \Gamma ^{z_0}_{\tilde c}(x_h-z_0, g^{z_0}(x_h-z_0,t_0-h),(x_0-z_0, g^{z_0}(x_0-z_0,t_0)))+h\e_h\\
&=t-h +  \widetilde \Gamma ^{z_0}_{\tilde c}(x_h-z_0, g^{z_0}(x_h-z_0,t_0-h),(x_0-z_0, g^{z_0}(x_0-z_0,t_0)))+h\e_h,
\end{align*}
which in turn gives 
\be\label{eq:weightinD1} 
\widetilde\Gamma ^{z_0}_{\tilde c}(x_h-z_0, g^{z_0}(x_h-z_0,t_0-h),(x_0-z_0, g^{z_0}(x_0-z_0,t_0))) \ge h -h\e_h.
\ee
For the opposite bound use superadditivity of the passage times,
\begin{align*}
t &= \widetilde \Gamma ^{z_0}_{\tilde c}(x_0-z_0, g^{z_0}(x_0-z_0,t_0))\\
&\ge\widetilde \Gamma ^{z_0}_{\tilde c}(x_h-z_0, g^{z_0}(x_h-z_0,t_0-h)) +  \widetilde \Gamma ^{z_0}_{\tilde c}(x_h-z_0, g^{z_0}(x_h-z_0,t_0-h),(x_0-z_0, g^{z_0}(x_0-z_0,t_0)))\notag\\
&=t-h + \widetilde\Gamma ^{z_0}_{\tilde c}(x_h-z_0, g^{z_0}(x_h-z_0,t_0-h),(x_0-z_0, g^{z_0}(x_0-z_0,t_0))),\notag
\end{align*}
giving that 
\be\label{eq:weightinD2} 
\widetilde\Gamma ^{z_0}_{\tilde c}(x_h-z_0, g^{z_0}(x_h-z_0,t_0-h),(x_0-z_0, g^{z_0}(x_0-z_0,t_0))) \le h.
\ee
From equations \eqref{eq:weightinD1}, \eqref{eq:weightinD2} we can write 
\be \label{eq:smallpathwiggle}
h - o(h) = \widetilde\Gamma ^{z_0}_{\tilde c}(x_h-z_0, g^{z_0}(x_h-z_0,t_0-h),(x_0-z_0, g^{z_0}(x_0-z_0,t_0))).
\ee

By the monotonicity properties of $g^{z_0}$, we have that 
$x_A \le x_h \le x_B $%, \quad y_A \ge y_h \ge y_B$, 
and therefore 
\[
x_0 - x_A \ge x_0 - x_h \ge x_0- x_B, %\quad 0= y_0 - y_A \le y_0 - y_h \le y_0- y_B. 
\]
Then equations \eqref{eq:xislope1},\eqref{eq:xislope2} give us 
\be \label{eq:xislope3}
1 \ge \frac{x_0-x_h}{h\tilde c_D^{\rm high}(h)} \ge -1.
\ee
Use this to define $\xi_h$ implicitly by
\be\label{eq:xislope4}
x_0 = x_h + h \xi_h, \quad \xi_h \in [-\tilde c_D^{\rm high}(h), \tilde c_D^{\rm high}(h)].
\ee 
Define $\tilde \gamma(x,y) := \gamma(x+y,y)$.
We can now estimate, from equation \eqref{eq:smallpathwiggle} that
\begin{align}
\frac{1}{\tilde c_D^{\rm high}(h)}&\tilde \gamma(x_0-x_h, g^{z_0}(x_0-z_0, t_0) - g^{z_0}(x_h-z_0, t_0-h)) \label{eq:A}\\
&\le h \notag\\
&\phantom{xxx}\le \frac{1}{\tilde c_D^{\rm low}(h)}\tilde \gamma(x_0-x_h, g^{z_0}(x_0-z_0, t_0) - g^{z_0}(x_h-z_0, t_0-h)) + o(h) \label{eq:UBV}
\end{align} 
Note that, by construction,
\be \label{eq:reachability}
x_0-x_h + g^{z_0}(x_0-z_0, t_0) - g^{z_0}(x_h-z_0, t_0-h) \ge 0,
\ee
by the construction. 

We only present the remaining argument starting from inequality \eqref{eq:UBV}. This will give that $v$ will be a classical subsolution. The bound for \eqref{eq:A}, and reaching a supersolution inequality follows in similar steps though you can also find a proof in Appendix \ref{app:C} using a short path construction from later parts in the article.
Set
\[
\alpha_h:=\frac{\xi_h}{\widetilde c_D^{\rm high}(h)},\qquad
Y_h:=
\frac{
g^{z_0}(x_0-z_0,t_0)-g^{z_0}(x_h-z_0,t_0-h)
}{
h\widetilde c_D^{\rm high}(h)
}.
\]
By wedge admissibility of the terminal segment from
$(x_h,g^{z_0}(x_h-z_0,t_0-h))$ to
$(x_0,g^{z_0}(x_0-z_0,t_0))$, we have
\[
Y_h\ge 0,
\qquad
\alpha_h+Y_h\ge 0.
\]
Passing to a subsequence, assume $\alpha_h\to\bar\alpha\in[-1,1]$.

If $\bar\alpha=-1$, then the second inequality gives
\be \label{eq:-1issue}
Y_h\ge -\alpha_h = 1+o(1)=\psi(\alpha_h)+o(1).
\ee

If $\bar\alpha=1$, then the first inequality gives
\be \label{eq:1issue}
Y_h\ge 0=\psi(1)=\psi(\alpha_h)+o(1).
\ee

Otherwise, consider all other possible accumulation points of $\frac{x_0 -x_h}{h \tilde c^{\rm high}_D(h)}$, away from $-1,1$.
We may assume without loss of generality that the subsequences remain in a compact subset of
$(-1,1)$, so the constant $C_{\alpha_h}$ in Lemma \ref{lem:psitaylor} is bounded below by a positive constant. 

Write \eqref{eq:UBV} as
\begin{align*}
1 - \e_h&\le \frac{1}{h\cdot \tilde c_D^{\rm low}(h)}\tilde \gamma(x_0-x_h, g^{z_0}(x_0-z_0, t_0) - g^{z_0}(x_h-z_0, t_0-h))\\
&=  \frac{\tilde c_D^{\rm high}(h)}{\tilde c_D^{\rm low}(h)}\tilde \gamma\left(\frac{x_0-x_h}{h \tilde c_D^{\rm high}(h)},\frac{g^{z_0}(x_0-z_0, t_0) - g^{z_0}(x_h-z_0, t_0-h)}{h \tilde c_D^{\rm high}(h)}\right), \text{by homogeneity of $\tilde \gamma$},\\
&\le \frac{\e_h + \tilde c_D^{\rm low}(h)}{ \tilde c_D^{\rm low}(h)} \tilde \gamma\left(\frac{\xi_h}{\tilde c_D^{\rm high}(h)},\frac{g^{z_0}(x_0-z_0, t_0) - g^{z_0}(x_h-z_0, t_0-h)}{h \tilde c_D^{\rm high}(h)}\right), \text{by \eqref{eq:eh},\eqref{eq:xislope4}},\\
&=\left(1 + \frac{\e_h}{\tilde c_D^{\rm low}(h)}\right)\tilde \gamma\left(\frac{\xi_h}{\tilde c_D^{\rm high}(h)},\frac{g^{z_0}(x_0-z_0, t_0) - g^{z_0}(x_h-z_0, t_0-h)}{h \tilde c_D^{\rm high}(h)}\right).
\end{align*} 
Let $\delta_h = \sqrt{\e_h/\tilde c_D^{\rm low}(h)}$ and assume by way of contradiction that 
\[
\frac{g^{z_0}(x_0-z_0, t_0) - g^{z_0}(x_h-z_0, t_0-h)}{h \tilde c_D^{\rm high}(h)} < \psi\left(\frac{\xi_h}{\tilde c_D^{\rm high}(h)} +\delta_h\right).
\]
By the monotonicity of $\tilde\gamma$ and Lemma \ref{lem:psitaylor} for $h$ small enough, we have that 
\be 
1-\e_h \le \left(1 + \frac{\e_h}{\tilde c_D^{\rm low}(h)}\right)\left(1 - C_{\xi_h, \tilde c} \delta_h \right) \le 1 - C \sqrt{\e_h}, 
\ee
which cannot hold as $\e_h \to 0$. Therefore, we have that 
\be \label{eq:psibound1}
\frac{g^{z_0}(x_0-z_0, t_0) - g^{z_0}(x_h-z_0, t_0-h)}{h \tilde c_D^{\rm high}(h)} \ge \psi\left(\frac{\xi_h}{\tilde c_D^{\rm high}(h)} +\delta_h\right).
\ee
Now we can prove that $v(x,t)$ is a classical subsolution to the PDE at differentiability points $(x_0,t_0)$.

We already have by construction that $v(x_0,t_0)=v_0(z_0) - g^{z_0}(x_0-z_0,t_0)$. Then
\begin{align}
v(x_0,t_0) &- v(x_0- h\xi_h, t_0-h) = v(x_0,t_0) - v(x_h, t_0-h) \notag \\
&= v_0(z_0) - g^{z_0}(x_0-z_0,t_0) - \sup_{z}\{ v_0(z) - g^{z}(x_h-z,t_0-h)\}\notag \\
&\le  -g^{z_0}(x_0-z_0, t_0) + g^{z_0}(x_0 -h\xi_h-z_0, t_0-h)\le -h\tilde c_D^{\rm high}(h)\psi\left(\frac{\xi_h}{\tilde c_D^{\rm high}(h)} +\delta_h\right).\label{eq:diffquo}
\end{align}
In the last inequality we used \eqref{eq:psibound1} but for accumulations at $-1,1$, use \eqref{eq:-1issue} or \eqref{eq:1issue} instead.
Work along a subsequence of $\xi_h \to \bar \xi$, where we still have \eqref{eq:localconv}. Equation \eqref{eq:xislope4} holds for $\bar \xi$ by compactness. From the differentiability assumption of $v$ at $x_0,t_0$, dividing by $h$ in \eqref{eq:diffquo} and taking the limit along the subsequence, we have that 
\be
v_t(x_0,t_0) + \bar \xi v_x(x_0,t_0) + \tilde c(x_0,-v(x_0,t_0))\psi\left(\frac{\bar \xi}{\tilde c(x_0,-v(x_0,t_0))} \right) \le 0.
\ee
Optimise over $\bar \xi/\tilde c(x_0,-v(x_0,t_0)) \in [-1,1]$ to eventually conclude
\be
v_t(x_0,t_0) + \tilde c(x_0,-v(x_0,t_0)) v_x(x_0,t_0)(1-v_x(x_0,t_0)) \le 0,
\ee
as required. 

The proof works at all points of differentiability of $v$ and continuity of $\tilde c(x,-v)$.  
From Lemma \ref{lem:continuity-points-full-measure} we have that these points have full measure and are dense in $\R^2$ which gives the second part of the theorem.
\end{proof}

\begin{remark}
The technicality with the accumulation point of speeds at $-1,1$ arises because the parallelogram in Lemma \ref{lem:recbelowlc} may be degenerate.
\qed \end{remark}

\subsection{The Lax-Oleinik representation} 
\label{sec:altrep}
In this section we prove that $v(x,t)$ is a viscosity solution for the PDE. We begin with the alternative representation for the scaling limit $v(x,t)$ given in Proposition \ref{thm:lang}, that mimics a Lax-Oleinik formula, as this will be the bridge towards the PDE considerations. 

\begin{proof}[Proof of Proposition \ref{thm:lang}.]
We start the proof by finding a way to map the paths $\mathbf{x}(\cdot)$ used in the definition of $\Gamma_c$ (see \eqref{eq:pathsH}) onto the paths ${\bf w}(\cdot)$ used in the definition of $v$ (see \eqref{eq:pathsw}). 
We first prove the inequality
\[
\sup_{q\in\mathbb R}\{v_0(q)-g^q(x-q,t)\}
\le
\sup_{{\bf w}\in \mathcal H^{v_0}_{x,t}}
\left\{
v_0(w_1(0))
-
\int_0^t
\tilde c({\bf w}(s))
\psi\!\left(\frac{w_1'(s)}{\tilde c({\bf w}(s))}\right)\,ds
\right\}.
\]
Fix $q\in\mathbb R$ and set
\[
y_q:=g^q(x-q,t).
\]
We first assume that $y_q$ is not on the boundary of the admissible wedge, and that
\[
\Gamma_c^q(x-q+y_q,y_q)=t.
\]
The boundary case is treated at the end of the argument.

Let $\varepsilon>0$. By the variational formula for $\Gamma_c^q$, and by the path-space reduction to paths
which lie in continuity regions of the coefficient for a.e. parameter value, we may choose a path
\[
{\bf x}(s)=(x_1(s),x_2(s))\in \mathcal H^{\rm cont}(x-q+y_q,y_q)
\]
such that
\[
\int_0^1
\frac{\gamma({\bf x}'(s))}
{c(x_1(s)+q-v_0(q), x_2(s)-v_0(q))}
\,ds
\ge t-\varepsilon .
\]
Put
\[
u_1(s):=x_1(s)-x_2(s),
\qquad
u_2(s):=x_2(s),
\]
which gives 
\[
u(0)=(0,0),
\qquad
u(1)=(x-q,y_q).
\]
Using $c(a,b)=\tilde c(a-b,b)$,
\[
c(x_1(s)+q-v_0(q),x_2(s)-v_0(q))
=
\tilde c(u_1(s)+q,u_2(s)-v_0(q)).
\]

The weight collected by the path up to $s$ is defined by 
\[
\tau(s):=
\int_0^s
\frac{\gamma(u_1'(\theta)+u_2'(\theta),u_2'(\theta))}
{\tilde c(u_1(\theta)+q,u_2(\theta)-v_0(q))}
\,d\theta .
\]
Let $T_\varepsilon:=\tau(1)$. Then $T_\varepsilon\ge t-\varepsilon$.
On every interval on which $u$ is $C^1$, and for a.e. $s$, we have
\be \label{eq:tauder}
\tau'(s)=
\frac{\gamma(u_1'(s)+u_2'(s),u_2'(s))}
{\tilde c(u_1(s)+q,u_2(s)-v_0(q))}.
\ee
At points where $\tau'(s)>0$, consider the inverse parametrisation  $s=s(\tau)$ and set
\be \label{eq:zd}
{\bf z}(\tau)= (z_1(\tau), z_2(\tau)) :=
\bigl(u_1(s(\tau))+q,\ u_2(s(\tau))-v_0(q)\bigr),
\qquad 0\le \tau\le T_\varepsilon .
\ee

The differentiation of $\tau$ in \eqref{eq:tauder} is meant in the a.e. sense. Since the integrand is locally integrable,
$\tau$ is absolutely continuous and $\tau'$ is equal to the integrand at Lebesgue points. To avoid any
ambiguity caused by paths travelling along discontinuity curves, we use the path-space reduction from
Lemma \ref{lem:betterpaths}, after applying the shear $(x_1,x_2)\mapsto(x_1-x_2,x_2)$. Thus the relevant paths may be
chosen so that they lie in continuity regions of $\tilde c$ for a.e. parameter value. 

Then
\[
{\bf z}(0)=(q,-v_0(q)),
\qquad
{\bf z}(T_\varepsilon)=(x,y_q-v_0(q))
\]
and moreover,
\[
z_1'(\tau)=u_1'(s)s'(\tau),\qquad
z_2'(\tau)=u_2'(s)s'(\tau).
\]
From this, for a.e.~ $\tau\in[0,T_\varepsilon]$ and using the homogeneity of $\gamma$ we can rescale \eqref{eq:tauder} to get
\begin{align}
    \gamma &\left( \frac{u_{1}'(s) + u_{2}'(s)}{\tau'(s)\tilde c(u_{1}(s) + q, u_{2}(s) -  v_{0}(q))}, \frac{u_{2}'(s)}{\tau'(s) \tilde c(u_{1}(s) + q, u_{2}(s) - v_{0}(q))} \right) \label{eq:taulc}\\
    &\phantom{xxxxxxxxxxxxxxxxxxxxxxxxxxxxxxxxxx}= \gamma\!\left(
\frac{z_1'(\tau)+z_2'(\tau)}
{\tilde c({\bf z}(\tau))},
\frac{z_2'(\tau)}
{\tilde c({\bf z}(\tau))}
\right)= 1.\notag 
\end{align}

Since both coordinates of $\bf x$ are nondecreasing, the two arguments of $\gamma$ above are nonnegative.
Using the parametrisation of the unit level curve of $\gamma$,
\[
\gamma(a,b)=1,\quad a,b\ge0
\quad\Longleftrightarrow\quad
b=\psi(a-b),
\]
we obtain
\[
\frac{z_2'(\tau)}{\tilde c({\bf z}(\tau))}
=
\psi\!\left(
\frac{z_1'(\tau)}{\tilde c({\bf z}(\tau))}
\right), 
\]
or equivalently,
\[
z_2'(\tau)
=
\tilde c({\bf z}(\tau))
\psi\!\left(
\frac{z_1'(\tau)}{\tilde c({\bf z}(\tau))}
\right)
\quad\text{for a.e. }\tau\in[0,T_\varepsilon].
\]
This is the desired second-coordinate relation for membership in $\mathcal H^{v_0}_{x,t}$.

If $T_\varepsilon=t$, then $z\in \mathcal H^{v_0}_{x,t}$, and
\[
v_0(q)-y_q
=
v_0(z_1(0))
-
\int_0^t
\tilde c({\bf z}(\tau))
\psi\!\left(
\frac{z_1'(\tau)}{\tilde c({\bf z}(\tau))}
\right)\,d\tau .
\]

If $T_\varepsilon<t$, extend ${\bf z}$ from $T_\varepsilon$ to $t$ by keeping the first coordinate equal to
$x$ and solving
\[
z_2'(\tau)
=
\tilde c(x,z_2(\tau))
\psi(0)
=
\frac14\tilde c(x,z_2(\tau)).
\]
Since $\tilde c$ is locally bounded away from 0 and $\infty$ on the compact region under consideration, this extension increases the terminal second coordinate by at most $C\varepsilon$, after replacing the approximating path if necessary.
Thus we obtain a path ${\bf z}^\varepsilon\in \mathcal H^{v_0}_{x,t}$ with
\[
z^\varepsilon_2(t)\le y_q-v_0(q)+C\varepsilon .
\]
Consequently,
\[
v_0(q)-y_q-C\varepsilon
\le
v_0(z^\varepsilon_1(0))
-
\int_0^t
\tilde c({\bf z}^\varepsilon(\tau))
\psi\!\left(
\frac{(z^\varepsilon_1)'(\tau)}
{\tilde c({\bf z}^\varepsilon(\tau))}
\right)\,d\tau .
\]
First take a supremum over paths $\mathcal H^{v_0}_{x,t}$ on the right-hand side and then let $\varepsilon\downarrow0$ to obtain 
\[
v_0(q)-g^q(x-q,t)
\le
\sup_{w\in \mathcal H^{v_0}_{x,t}}
\left\{
v_0(w_1(0))
-
\int_0^t
\tilde c({\bf w}(s))
\psi\!\left(\frac{w_1'(s)}{\tilde c({\bf w}(s))}\right)\,ds
\right\}.
\]
Finally take the supremum over $q$.

It remains to consider the case where the infimum is attained on the wedge boundary. Write
$r=x-q$ and $y_q=g^q(r,t)$. Since
\[
y_q=\beta(r)=(-r)^+,
\]
one of the two LPP coordinates $r+y_q$ and $y_q$ is zero.

First suppose $r\ge0$. Then $y_q=0$, and the LPP endpoint is $(r,0)$. By the definition of $g^q$, $\Gamma_c^q(r,0)\ge t$. Parametrise the horizontal LPP-axis path from $(0,0)$ to
$(r,0)$ by physical time $s\in[0,t]$ so that
\[
\frac{w_1'(s)}{\tilde c({\bf w}(s))}=\frac{\Gamma_c^q(r,0)}{t}\ge1,
\qquad
w_2(s)\equiv -v_0(q).
\]
Then $\psi(w_1'/\tilde c)=0$, and hence
\[
w_2'(s)=0
=
\tilde c({\bf w}(s))
\psi\!\left(\frac{w_1'(s)}{\tilde c({\bf w}(s))}\right).
\]
Thus ${\bf w}\in \mathcal H^{v_0}_{x,t}$, and $w_2(t)=-v_0(q)=y_q-v_0(q)$.

Now suppose $r<0$. Then $y_q=-r$, and the LPP endpoint is $(0,y_q)$. Parametrise the vertical LPP-axis path by physical time $s\in[0,t]$ so that
\[
\frac{w_1'(s)}{\tilde c({\bf w}(s))}=-\frac{\Gamma_c^q(0,y_q)\ge t}{t}\le -1,
\qquad
w_2'(s)=-w_1'(s).
\]
Since $\psi(\lambda)=-\lambda$ for $\lambda\le -1$, we obtain
\[
w_2'(s)
=
-w_1'(s)
=
\tilde c({\bf w}(s))
\psi\!\left(\frac{w_1'(s)}{\tilde c({\bf w}(s))}\right).
\]
Again ${\bf w}\in \mathcal H^{v_0}_{x,t}$, and its terminal second coordinate is
\[
w_2(t)=y_q-v_0(q).
\]
Thus the boundary case gives the same contribution $v_0(q)-y_q$.

Conversely, let ${\bf w}\in \mathcal H^{v_0}_{x,t}$. Put $q=w_1(0)$ and define
\[
u_1(s):=w_1(s)-q,\qquad
u_2(s):=w_2(s)+v_0(q),\qquad 0\le s\le t.
\]
Set
\[
x_1(s):=u_1(s)+u_2(s),\qquad x_2(s):=u_2(s),
\]
and let ${\bf x}(s) = (x_1(s), x_2(s))$.
Then ${\bf x}(0)=(0,0)$ and
\[
{\bf x}(t)=\bigl(x-q+y,y\bigr),
\quad \text{ where }
y:=w_2(t)+v_0(q).
\]
Moreover, for a.e.~ $s$
\[
x_2'(s)=w_2'(s)
=
\tilde c({\bf w}(s))
\psi\!\left(
\frac{w_1'(s)}{\tilde c({\bf w}(s))}
\right)\ge0,
\]
and
\[
x_1'(s)=w_1'(s)+w_2'(s)
=
\tilde c({\bf w}(s))
\left[
\frac{w_1'(s)}{\widetilde c(w(s))}
+
\psi\!\left(
\frac{w_1'(s)}{\tilde c({\bf w}(s))}
\right)
\right]\ge0.
\]
Thus ${\bf x}$ is an admissible LPP path. Furthermore, for a.e. $s$,
\[
\frac{\gamma({\bf x}'(s))}
{c(x_1(s)+q-v_0(q),x_2(s)-v_0(q))}
=
\gamma\!\left(
\lambda(s)+\psi(\lambda(s)),\psi(\lambda(s))
\right),
\]
where
\[
\lambda(s):=\frac{w_1'(s)}{\tilde c({\bf w}(s))}.
\]
By the definition of $\psi$, $
\gamma(\lambda+\psi(\lambda),\psi(\lambda))\ge1$,
for all $\lambda\in\mathbb R$ and hence $\Gamma_c^q(x-q+y,y)\ge t.$

Therefore  $y\ge g^q(x-q,t)$, and so
\[
v_0(w_1(0))-
\int_0^t
\tilde c({\bf w}(s))
\psi\!\left(
\frac{w_1'(s)}{\tilde c({\bf w}(s))}
\right)\,ds
=
v_0(q)-y
\le
v_0(q)-g^q(x-q,t)
\le
\sup_{r\in\mathbb R}\{v_0(r)-g^r(x-r,t)\}.
\]
Taking the supremum over ${\bf w}\in \mathcal H^{v_0}_{x,t}$ gives the reverse inequality.
\end{proof}

\begin{remark}\label{rem:reparapath}
The first part of the previous proof shows how an LPP-admissible path in $\mathcal H^{\rm cont}$ of weight $T$ can be reparametrised and transformed to a path belonging to $\mathcal H^{v_0}_{\cdot, T}$.   \qed \end{remark}

\subsection{Short paths}
Define $\bar v_0(x) = v(x,t - h)$ and consider the LLN starting from this initial condition but running only up to time $h$

\be\label{eq:vbardefq}
\bar v(x, h) = \sup_{q \in \R}\{ \bar v_0(q) - \bar g^q(x-q,h) \} =  \sup_{q \in \R}\{v(q, t-h) - \bar g^q(x-q,h) \},
\ee
where here 
\[\bar g^q(x-q, h) = \inf\{ y: \widetilde \Gamma_{\tilde c} ((q,  - \bar v_0(q)) , (x,  y- \bar v_0(q))) \ge h \}.\]
As this is the LLN starting from a different initial condition, all previous variational characterisations work for $\bar v (x,h)$. In the homogeneous environment, one can show that $\bar v(x,h) = v(x,t)$ which simplifies a lot of the following arguments. Here, the inhomogeneous $\tilde c$ does not allow for path shifting, so we do need different comparisons. 
We can indeed obtain a comparison between $\bar v(x,h)$ and $v(x,t)$. 
We first need a preliminary comparison that shows a super-additivity of level curves.

\begin{lemma}[Subadditivity of level curves]
\label{lem:lcsup}
Fix $q_0\in \mathbb{R}$, $x\in\mathbb{R}$, $t>0$, and
$0<h<t$. For $q\in\mathbb{R}$ set
\[
    y_1 := g^{q_0}(q-q_0,t-h).
\]
and define
\[
\bar g^{q,q_0}(x-q,h)
:=
\inf\Big\{
    y\geq \beta(x-q):
    \widetilde\Gamma_{\tilde c}
    \Big(
        (q, y_1-v_0(q_0)),
        (x, y+y_1-v_0(q_0))
    \Big)
    \geq h
\Big\}.
\]
Then
\be
\label{eq:supeqbar}
    g^{q_0}(q-q_0,t-h)
    +
    \bar g^{q,q_0}(x-q,h)
    \geq
    g^{q_0}(x-q_0,t).
\ee
\end{lemma}

\begin{proof}
Write
\[
    y_1 := g^{q_0}(q-q_0,t-h),
    \qquad
    y_2 := \bar g^{q,q_0}(x-q,h).
\]
By the definition of $g^{q_0}$,
\[
    \widetilde\Gamma_{\tilde c}
    \Big(
        (q_0,-v_0(q_0)),
        (q,y_1-v_0(q_0))
    \Big)
    \geq t-h.
\]
Moreover, by the definition of $y_2$, for every $\varepsilon>0$ there
exists $y_\varepsilon\geq \beta(x-q)$ such that
\[
    y_\varepsilon \leq y_2+\varepsilon
\]
and
\[
    \widetilde\Gamma_{\tilde c}
    \Big(
        (q,y_1-v_0(q_0)),
        (x,y_1+y_\varepsilon-v_0(q_0))
    \Big)
    \geq h.
\]
The endpoint
$
    (x,y_1+y_\varepsilon-v_0(q_0))
$
is admissible from $(q_0,-v_0(q_0))$. Indeed,
\[
    y_1\geq \beta(q-q_0),
    \qquad
    y_\varepsilon\geq \beta(x-q),
\]
and therefore
\[
    y_1+y_\varepsilon\geq 0,
    \qquad
    x-q_0+y_1+y_\varepsilon
    =
    (q-q_0+y_1)+(x-q+y_\varepsilon)
    \geq 0.
\]
By superadditivity of the last-passage functional under concatenation,
\[
\begin{aligned}
    &\widetilde\Gamma_{\tilde c}
    \Big(
        (q_0,-v_0(q_0)),
        (x,y_1+y_\varepsilon-v_0(q_0))
    \Big)
    \\
    &\qquad\geq
    \widetilde\Gamma_{\tilde c}
    \Big(
        (q_0,-v_0(q_0)),
        (q,y_1-v_0(q_0))
    \Big)
    +
    \widetilde\Gamma_{\tilde c}
    \Big(
        (q,y_1-v_0(q_0)),
        (x,y_1+y_\varepsilon-v_0(q_0))
    \Big)
    \\
    &\qquad\geq
    (t-h)+h=t.
\end{aligned}
\]
Hence $y_1+y_\varepsilon$ is an admissible value in the defining
infimum for $g^{q_0}(x-q_0,t)$. Consequently,
\[
    g^{q_0}(x-q_0,t)
    \leq y_1+y_\varepsilon
    \leq y_1+y_2+\varepsilon.
\]
Letting $\varepsilon\downarrow0$ gives
\[
    g^{q_0}(x-q_0,t)
    \leq
    g^{q_0}(q-q_0,t-h)
    +
    \bar g^{q,q_0}(x-q,h),
\]
which is the desired inequality.
\end{proof}
\tikzset{every picture/.style={line width=0.75pt}} %set default line width to 0.75pt        
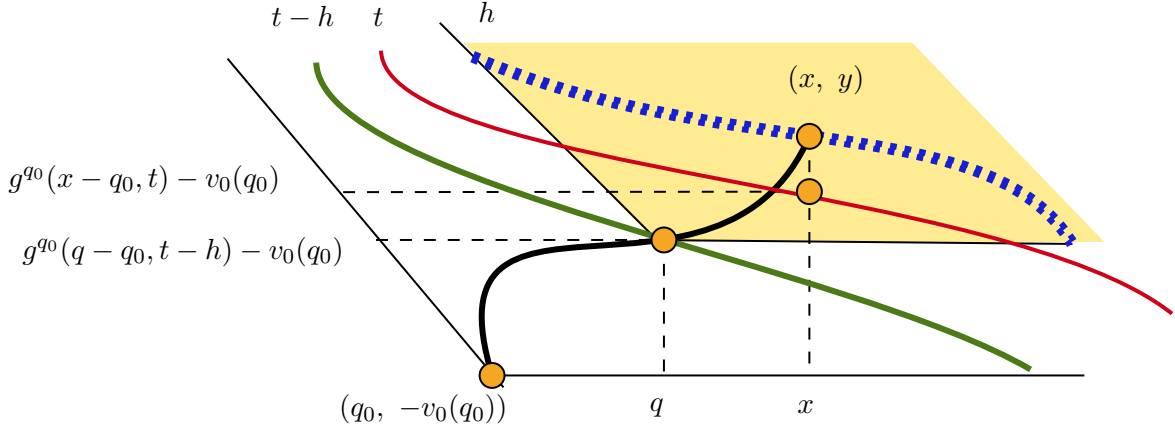
\begin{figure}[ht]
\begin{tikzpicture}[x=0.75pt,y=0.75pt,yscale=-1,xscale=1]
%uncomment if require: \path (0,1203); %set diagram left start at 0, and has height of 1203

%Shape: Parallelogram [id:dp9699292286110465] 
\draw  [draw opacity=0][fill={rgb, 255:red, 255; green, 236; blue, 146 }  ,fill opacity=1 ] (511.05,976) -- (286,976) -- (382.45,1076) -- (607.5,1076) -- cycle ;
%Straight Lines [id:da5865400871623307] 
\draw  [dash pattern={on 4.5pt off 4.5pt}]  (225.5,1051) -- (459.5,1050.75) ;
%Straight Lines [id:da012044699908861145] 
\draw    (377,1075) -- (591.5,1077) ;
%Straight Lines [id:da9854829592866671] 
\draw    (274,966) -- (389,1081) ;
%Straight Lines [id:da36731565113649545] 
\draw  [dash pattern={on 4.5pt off 4.5pt}]  (459.5,1023) -- (459.5,1143) ;
%Straight Lines [id:da27212666422994647] 
\draw  [dash pattern={on 4.5pt off 4.5pt}]  (386.5,1075) -- (386.5,1143) ;
%Curve Lines [id:da6673824424719812] 
\draw [line width=2.25]    (300.25,1143) .. controls (265.5,1027) and (414.5,1126) .. (459.5,1023) ;
%Straight Lines [id:da3459518597727794] 
\draw    (167.5,984) -- (306,1149) ;
%Straight Lines [id:da49364851770792784] 
\draw    (294,1143) -- (597.5,1143) ;
%Straight Lines [id:da21023496368801897] 
\draw  [dash pattern={on 4.5pt off 4.5pt}]  (242,1075) -- (386.5,1075) ;
%Curve Lines [id:da5392343521393179] 
\draw [color={rgb, 255:red, 208; green, 2; blue, 27 }  ,draw opacity=1 ][line width=1.5]    (244.5,980) .. controls (247,1041) and (558.5,1046) .. (641.5,1112) ;
%Curve Lines [id:da11840674975852328] 
\draw [color={rgb, 255:red, 75; green, 122; blue, 22 }  ,draw opacity=1 ][line width=2.25]    (212,986) .. controls (214.5,1047) and (479.5,1087) .. (570.5,1140) ;
%Shape: Circle [id:dp5150994280808098] 
\draw  [fill={rgb, 255:red, 245; green, 166; blue, 35 }  ,fill opacity=1 ] (453.25,1050.75) .. controls (453.25,1047.3) and (456.05,1044.5) .. (459.5,1044.5) .. controls (462.95,1044.5) and (465.75,1047.3) .. (465.75,1050.75) .. controls (465.75,1054.2) and (462.95,1057) .. (459.5,1057) .. controls (456.05,1057) and (453.25,1054.2) .. (453.25,1050.75) -- cycle ;
%Shape: Circle [id:dp8789419902130196] 
\draw  [fill={rgb, 255:red, 245; green, 166; blue, 35 }  ,fill opacity=1 ] (294,1143) .. controls (294,1139.55) and (296.8,1136.75) .. (300.25,1136.75) .. controls (303.7,1136.75) and (306.5,1139.55) .. (306.5,1143) .. controls (306.5,1146.45) and (303.7,1149.25) .. (300.25,1149.25) .. controls (296.8,1149.25) and (294,1146.45) .. (294,1143) -- cycle ;
%Shape: Circle [id:dp816065694950747] 
\draw  [fill={rgb, 255:red, 245; green, 166; blue, 35 }  ,fill opacity=1 ] (380.25,1075) .. controls (380.25,1071.55) and (383.05,1068.75) .. (386.5,1068.75) .. controls (389.95,1068.75) and (392.75,1071.55) .. (392.75,1075) .. controls (392.75,1078.45) and (389.95,1081.25) .. (386.5,1081.25) .. controls (383.05,1081.25) and (380.25,1078.45) .. (380.25,1075) -- cycle ;
%Curve Lines [id:da3470044500696592] 
\draw [color={rgb, 255:red, 23; green, 30; blue, 211 }  ,draw opacity=1 ][line width=2.25]  [dash pattern={on 2.53pt off 3.02pt}]  (291.09,981.62) .. controls (295.25,983.41) and (299.41,985.1) .. (303.57,986.72) .. controls (358.45,1007.99) and (413.34,1014.72) .. (460.95,1021.46) .. controls (521.36,1030.01) and (569.99,1038.78) .. (592.78,1076.22)(289.91,984.38) .. controls (294.1,986.18) and (298.29,987.89) .. (302.49,989.51) .. controls (357.6,1010.88) and (412.71,1017.66) .. (460.53,1024.43) .. controls (519.72,1032.81) and (567.81,1040.97) .. (590.22,1077.78) ;
%Shape: Circle [id:dp47869073770252246] 
\draw  [fill={rgb, 255:red, 245; green, 166; blue, 35 }  ,fill opacity=1 ] (453.25,1023) .. controls (453.25,1019.55) and (456.05,1016.75) .. (459.5,1016.75) .. controls (462.95,1016.75) and (465.75,1019.55) .. (465.75,1023) .. controls (465.75,1026.45) and (462.95,1029.25) .. (459.5,1029.25) .. controls (456.05,1029.25) and (453.25,1026.45) .. (453.25,1023) -- cycle ;

% Text Node
\draw (239,957) node [anchor=north west][inner sep=0.75pt]    {$t$};
% Text Node
\draw (188,956) node [anchor=north west][inner sep=0.75pt]    {$t-h$};
% Text Node
\draw (447,986) node [anchor=north west][inner sep=0.75pt]    {$( x,\ y)$};
% Text Node
\draw (63,1072) node [anchor=north west][inner sep=0.75pt]    {$g^{q_{0}}( q-q_{0} ,t-h)-v_{0}( q_{0})$};
% Text Node
\draw (222,1151) node [anchor=north west][inner sep=0.75pt]    {$( q_{0} ,\ -v_{0}( q_{0}))$};
% Text Node
\draw (378,1153) node [anchor=north west][inner sep=0.75pt]    {$q\ $};
% Text Node
\draw (452,1154) node [anchor=north west][inner sep=0.75pt]    {$x$};
% Text Node
\draw (56,1036) node [anchor=north west][inner sep=0.75pt]    {$g^{q_{0}}( x-q_{0} ,t)-v_{0}( q_{0})$};
% Text Node
\draw (292,954) node [anchor=north west][inner sep=0.75pt]    {$h$};
\end{tikzpicture}
\caption{Schematic proof of the argument in Lemma \ref{lem:lcsup}. The two thick level curves are those for the passage time starting from $(q_0,-v_0(q_0))$ reaching level $t-h$ and $t$. The thick dashed level curve is an $h$-level curve starting from $Q =(q, g^{q_{0}}( q-q_{0} ,t-h)-v_{0}( q_{0}))$. 
The thick path is a concatenation of an optimal path from $(q_0,-v_0(q_0))$ to $Q$ with weight $t-h$ and an optimal path from $Q$ to $(x,y)$ with weight $h$. As such, its weight $t$, is not larger than the optimal weight up to $(x,y)$, giving that $y \ge g^{q_{0}}( x-q_{0} ,t) \ -v_{0}( q_{0})$.}
\end{figure}

\begin{lemma}\label{lem:vbarvcomp} Let $v(x,t)$ be given by \eqref{eq:goodv} and let $\bar v(x,h)$ be given by \eqref{eq:vbardefq} (equiv. \eqref{eq:vbardef} below). Then 
\[
\bar v(x,h) \le v(x,t).
\]
\end{lemma}

\begin{proof}[Proof of Lemma \ref{lem:vbarvcomp}]
Let $q$ be the point that achieves the supremum in \eqref{eq:vbardefq}
so that $
\bar v(x,h) = v(q,t-h)- \bar g^q(x-q,h)$
and $q_0$ be such that $v(q,t-h) = v_0(q_0) -  g^{q_0}(q-q_0, t-h)$. 

Note that $\bar g^q(x-q,h)$ is a level curve for a last passage time starting at location $(q,- v(q,t-h)) = (q,g^{q_0}(q-q_0,t-h) - v_0(q_0))$ which means that $\bar g^q(x-q,h) = \bar g^{q,q_0}(x-q,h)$, as needed from Lemma \ref{lem:lcsup}. The admissibility conditions are automatically satisfied by the way we found $q, q_0$.
Then we can estimate
\begin{align*}
\bar v(x, h) &= v(q,t-h) - \bar g^{q,q_0}(x-q,h) =  v_0(q_0) -  g^{q_0}(q-q_0, t-h)-\bar g^{q,q_0}(x-q,h) \\
&= v_0(q_0) -  (g^{q_0}(q-q_0, t-h) +\bar g^{q,q_0}(x-q,h)) \le v_0(q_0) -  g^{q_0}(x-q_0, t), \quad \text{ by  Lemma \ref{lem:lcsup}}\\
&\le v(x,t).\qedhere
\end{align*}
\end{proof}

 Moreover, by Proposition \ref{thm:lang}, we can also write 
\begin{align}
v(x,t) \ge \bar v(x, h)  &= \sup_{{\bf w}(\cdot) \in \mathcal{H}^{\bar v_0}_{x,h}} \left\{\bar v_{0}(w_{1}(0)) - \int^{h}_{0}\tilde c({\bf w}(s))\psi\left( \frac{w_{1}'(s)}{\tilde c({\bf w}(s))} \right) ds \right\} \notag\\
&=\sup_{{\bf w}(\cdot) \in \mathcal{H}^{\bar v_0}_{x,h}} \left\{v(w_{1}(0), t-h) - \int^{h}_{0}\tilde c({\bf w}(s))\psi\left( \frac{w_{1}'(s)}{\tilde c({\bf w}(s))} \right) ds \right\} \label{eq:vbardef}
\end{align}

In particular, for any $\xi \in \R$, we can start a path ${\bf x}$ so that $x_1(0)=x -\xi h$ and $x_1(h)=x$, to have that 
\be\label{eq:lb4visco}
v(x,t) \ge v( x-\xi h, t-h) - \int_{0}^h \tilde c({\bf x}(s))\psi\left( \frac{x_{1}'(s)}{\tilde c({\bf x}(s))} \right) ds.
\ee

\begin{lemma}\label{lem:pathapprox} For any fixed $(x,t)$ and $h >0$ we can find functions $d(h) \le 2e(h)=o(h)$ and a path ${\bf w} \in \mathcal H^{v_0}_{x,t- e(h)}$ such that: 
\begin{enumerate}
\item ${\bf w}(0) = (w_1(0), -v_0(w_1(0)))$ and ${\bf w}(t-e(h)) = (x, -v(x,t)),$
\item there exists a constant $C_{\max}(x,t) \in (0,\infty)$ and a $\xi_h \in [-C_{\max}(x,t), C_{\max}(x,t)]$ so that $${\bf w}(t-h-d(h)) := (x_h, y_h) = (x - h\xi_h, g^{z_0}(x_h-z_0, t-h)-v_0(z_0)),$$
where $-v(x,t) = g^{z_0}(x-z_0, t)-v_0(z_0)$,
\item we have the following equality
\be\label{eq:438}
v(x,t) \le v(x_h, t-h-d(h)) - \int_{t-h-d(h)}^{t-e(h)}\tilde c({\bf w}(s))\psi \left( \frac{w_1'(s)}{\tilde c({\bf w}(s))} \right)\,ds.
\ee
\item and 
\be \label{eq:localwedgeh}
        \sup_{s\in[t-h-d(h),\,t-e(h)]}
        \left|
        {\bf w}(s)-(x,-v(x,t))
        \right|
        \longrightarrow 0
        \qquad\text{as }h\downarrow0.
\ee
\end{enumerate}
\end{lemma}
\begin{proof}
%From the proof of Theorem \ref{thm: existence} and 
For any $(x,t)$ we can find a $z_0$ so that $v_0(z_0) - g^{z_0}(x-z_0,t) = v(x,t)$ and that satisfies \eqref{eq:maximiser-levelcurve}. 

Once $z_0$ is specified, from the variational description for the wedge shape function \eqref{eq:LLPP}, we can find a path ${\bf w}$ starting from $(z_0, -v_0(z_0))$ and ending at $(x,-v(x,t))$ with weight $t-e(h)$. Moreover, with lower-left corner $(z_0, -v_0(z_0))$ consider the level curve $g^{z_0}(\cdot, t-h) - v_0(z_0)$. Let $(x_h,y_h)$ denote its intersection with ${\bf w}$.

Let $I_1({\bf w})$ denote the weight of the path up to $(x_h,y_h)$ and $I_2({\bf w})$ the weight collected on it from 
$(x_h,y_h)$ to $(x,-v(x,t))$.

 Since $(x_h,y_h)$ is on the lifted level curve, we have that $I_1({\bf w}) = t-h - d_1(t,h)$, where $d_1(t,h)$ is some positive quantity. Then we have that $I_2({\bf w}) = h- e(h) +d_1(t,h)$ since the total weight is $t-e(h)$.
 
Now, by possibly altering the first part of the path to a better maximizer $\hat{\bf w}$, we can also guarantee that we can find a function $d_2(h) \le e(h)$ so that %\note{Check this claim.}
\[
{\hat{\bf w}}(t-h-d_2(h)) = (x_h, y_h) = (x - h\xi_h, g^{z_0}(x_h-z_0, t-h)-v_0(z_0)),
\]
for some $\xi_h \in [-C_{\max}(x,t), C_{\max}(x,t)]$, where $C_{\max}(x,t)$ is the largest value of $\tilde c$ in a sufficiently large compact set containing $(x,-v(x,t))$, as in the construction in the proof of Theorem \ref{thm: existence}. 

But then the concatenated path $\hat{\bf w} \sqcup \bf w $ has total weight $t - d_2(h) + d_1(t,h) - e(h)$. Since this must still be $\le t$, we have that the original error $d_1(t,h) \le 2e(h) = o(h)$ and we do not need to alter the original path $\bf w$ and we set $d(h) = d_1(t,h)$.  
See also Fig.\ref{fig:pathapprox} for a diagram containing this construction.
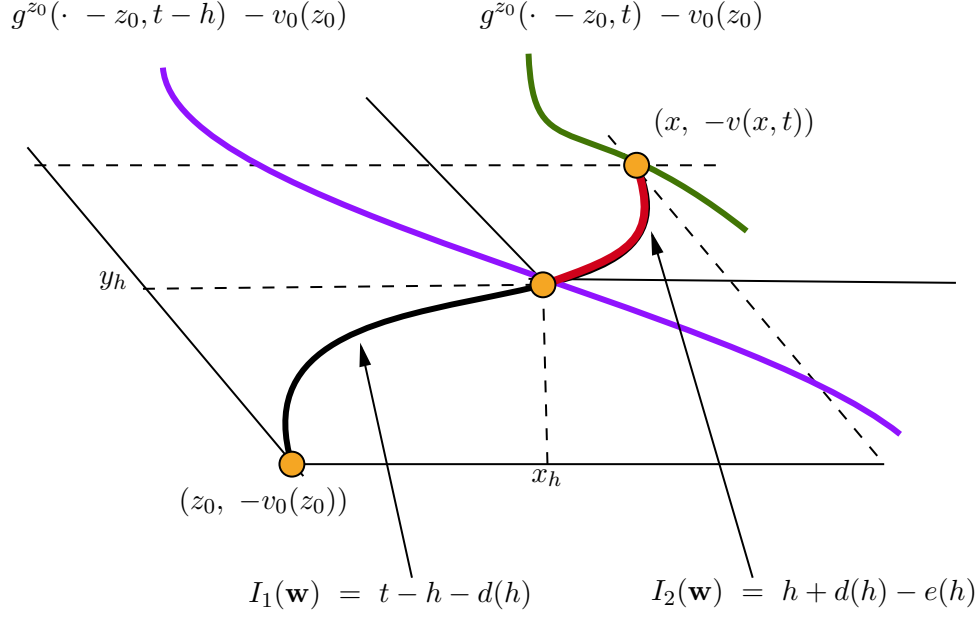
\begin{figure}[ht]
\tikzset{every picture/.style={line width=0.75pt}} %set default line width to 0.75pt        

\begin{tikzpicture}[x=0.75pt,y=0.75pt,yscale=-1,xscale=1]
%uncomment if require: \path (0,2645); %set diagram left start at 0, and has height of 2645

%Straight Lines [id:da7515796990872653] 
\draw    (316,2438) -- (530.5,2440) ;
%Straight Lines [id:da8191575369041483] 
\draw    (234.5,2347) -- (328,2444) ;
%Straight Lines [id:da99191229232572] 
\draw  [dash pattern={on 4.5pt off 4.5pt}]  (323.25,2441) -- (325.5,2531) ;
%Curve Lines [id:da418606019944906] 
\draw [line width=2.25]    (197.25,2531) .. controls (162.5,2415) and (414.5,2478) .. (370.25,2381) ;
%Straight Lines [id:da7753962693181496] 
\draw    (64.5,2372) -- (203,2537) ;
%Straight Lines [id:da1727809527014087] 
\draw    (191,2531) -- (494.5,2531) ;
%Straight Lines [id:da6640235099234436] 
\draw  [dash pattern={on 4.5pt off 4.5pt}]  (68,2381) -- (412.5,2381) ;
%Straight Lines [id:da5547422215356665] 
\draw  [dash pattern={on 4.5pt off 4.5pt}]  (356,2366) -- (494.5,2531) ;
%Curve Lines [id:da573156285378803] 
\draw [color={rgb, 255:red, 65; green, 117; blue, 5 }  ,draw opacity=1 ][line width=2.25]    (316,2325) .. controls (318.5,2386) and (342.5,2348) .. (425.5,2414) ;
%Curve Lines [id:da046581797373653755] 
\draw [color={rgb, 255:red, 144; green, 19; blue, 254 }  ,draw opacity=1 ][line width=2.25]    (132.5,2332) .. controls (138,2402) and (419.5,2450) .. (502.5,2516) ;
%Shape: Circle [id:dp6111313819422785] 
\draw  [fill={rgb, 255:red, 245; green, 166; blue, 35 }  ,fill opacity=1 ] (191,2531) .. controls (191,2527.55) and (193.8,2524.75) .. (197.25,2524.75) .. controls (200.7,2524.75) and (203.5,2527.55) .. (203.5,2531) .. controls (203.5,2534.45) and (200.7,2537.25) .. (197.25,2537.25) .. controls (193.8,2537.25) and (191,2534.45) .. (191,2531) -- cycle ;
%Straight Lines [id:da5235371145218722] 
\draw  [dash pattern={on 4.5pt off 4.5pt}]  (122.5,2443) -- (317,2441) ;
%Straight Lines [id:da22644399907610757] 
\draw    (257.5,2588) -- (231.93,2471.95) ;
\draw [shift={(231.5,2470)}, rotate = 77.57] [fill={rgb, 255:red, 0; green, 0; blue, 0 }  ][line width=0.08]  [draw opacity=0] (12,-3) -- (0,0) -- (12,3) -- cycle    ;
%Straight Lines [id:da7818662279461956] 
\draw    (428.5,2588) -- (378.06,2415.92) ;
\draw [shift={(377.5,2414)}, rotate = 73.66] [fill={rgb, 255:red, 0; green, 0; blue, 0 }  ][line width=0.08]  [draw opacity=0] (12,-3) -- (0,0) -- (12,3) -- cycle    ;
%Curve Lines [id:da8504376269241103] 
\draw [color={rgb, 255:red, 208; green, 2; blue, 27 }  ,draw opacity=1 ][line width=3]    (323.25,2441) .. controls (366.5,2426) and (382.5,2419) .. (370.25,2381) ;
%Shape: Circle [id:dp7559019801425302] 
\draw  [fill={rgb, 255:red, 245; green, 166; blue, 35 }  ,fill opacity=1 ] (364,2381) .. controls (364,2377.55) and (366.8,2374.75) .. (370.25,2374.75) .. controls (373.7,2374.75) and (376.5,2377.55) .. (376.5,2381) .. controls (376.5,2384.45) and (373.7,2387.25) .. (370.25,2387.25) .. controls (366.8,2387.25) and (364,2384.45) .. (364,2381) -- cycle ;
%Shape: Circle [id:dp15013697067414966] 
\draw  [fill={rgb, 255:red, 245; green, 166; blue, 35 }  ,fill opacity=1 ] (317,2441) .. controls (317,2437.55) and (319.8,2434.75) .. (323.25,2434.75) .. controls (326.7,2434.75) and (329.5,2437.55) .. (329.5,2441) .. controls (329.5,2444.45) and (326.7,2447.25) .. (323.25,2447.25) .. controls (319.8,2447.25) and (317,2444.45) .. (317,2441) -- cycle ;

% Text Node
\draw (377,2351) node [anchor=north west][inner sep=0.75pt]    {$( x,\ -v( x,t))$};
% Text Node
\draw (139,2541) node [anchor=north west][inner sep=0.75pt]    {$( z_{0} ,\ -v_{0}( z_{0}))$};
% Text Node
\draw (316,2532) node [anchor=north west][inner sep=0.75pt]    {$x_{h} \ $};
% Text Node
\draw (99,2433) node [anchor=north west][inner sep=0.75pt]    {$y_{h} \ $};
% Text Node
\draw (55,2297) node [anchor=north west][inner sep=0.75pt]    {$g^{z_{0}}( \cdot \ -z_{0} ,t-h) \ -v_{0}( z_{0})$};
% Text Node
\draw (290,2297) node [anchor=north west][inner sep=0.75pt]    {$g^{z_{0}}( \cdot \ -z_{0} ,t) \ -v_{0}( z_{0})$};
% Text Node
\draw (174,2588) node [anchor=north west][inner sep=0.75pt]    {$I_{1}( {\bf w}) \ =\ t-h-d( h) \ $};
% Text Node
\draw (376,2587) node [anchor=north west][inner sep=0.75pt]    {$I_{2}( {\bf w}) \ =\ h+d( h) -e( h) \ $};
\end{tikzpicture}
\caption{Diagram for the proof of Lemma \ref{lem:pathapprox}. The increasing path is a near optimiser for the wedge last passage time shape function \eqref{eq:LLPP} up to $(x,-v(x,t))$. Since that point belongs on the $t$-level curve, a near optimiser will have weight $t - e(h)$ with $e(h)$ as near zero as we wish. It will also intersect the $t-h$ level curve at some point $(x_h, y_h)$. This naturally splits $\bf w$ into two pieces, where the weight of each path segment are also near-optimal, respectively $t-h-d(h)$ and $h+d(h)-e(h)$. }
\label{fig:pathapprox}
\end{figure}

This path can be reparametrised as a path ${\bf w} (s), s\in [0,t]$ using the inverse last passage time as in the first part of the proof of Proposition \ref{thm:lang} (see also Remark \ref{rem:reparapath}), to a path ${\bf w} \in \mathcal H^{v_0}_{x,t-e(h)}$, so that the reparametrised path $\bf w$ satisfies 
\[
{\bf w}(0) = (w_1(0), -v_0(w_1(0)))\textrm{ \, and \, }{\bf w}(t-e(h)) = (x, -v(x,t)),
\]

Then we can write
\begin{align*}
v(x_h, t-h-d(h))&\ge-w_2(t-h-d(h)) = -w_2(t-h-d(h)) + w_2(t-e(h)) - w_2(t-e(h)) \\
&=-w_2(t-e(h)) + \left(w_2(t-e(h))-w_2(t-h-d(h)) \right)\\
&=-w_2(t-e(h)) + \int_{t-h-d(h)}^{t-e(h)}\tilde c({\bf w}(s))\psi \left(\frac{w'_1(s)}{\tilde c({\bf w}(s))} \right)\,ds, \text{\; since \;} {\bf w} \in \mathcal H^{v_0}_{x,t-e(h)} \\
&=v(x,t)+ \int_{t-h-d(h)}^{t-e(h)}\tilde c({\bf w}(s))\psi \left(\frac{w'_1(s)}{\tilde c({\bf w}(s))} \right)\,ds.
\end{align*}
Rearrange to obtain \eqref{eq:438}.

Finally, for the local compactness property of the terminal segment,
consider for the moment the first coordinate. By construction, $ w_1(t-e(h))=x,$ 
$ w_1(t-h-d(h))=x_h=x-h\xi_h,$
where $\xi_h$ is bounded uniformly for $h$ small. Hence
\[
        |x_h-x|\le h|\xi_h|\to0.
\]
Moreover the terminal segment of $\bf w$ is contained in the monotone wedge
between its endpoints, as it is a segment of an admissible path. 
We now consider the second coordinate. The terminal endpoint is 
$w_2(t-e(h))=-v(x,t)$ and the initial point of the terminal segment is by construction
\[
        w_2(t-h-d(h))=g^{z_0}(x_h-z_0, t-h)-v_0(z_0)\ge -v(x_h,t-h-d(h)).
\]
Then we have 
\[
0 \le w_2(t-e(h)) - w_2(t-h-d(h)) \le v(x_h,t-h-d(h))- v(x,t)
\]
 Since $v$ is locally Lipschitz and $ x_h\to x$, $t-h-d(h)\to t,$ as $h\to 0$, the RHS in the inequality above converges to 0 and therefore
we obtain
\[
        w_2(t-h-d(h))
        \longrightarrow
        -v(x,t).
\]
Finally, the second coordinate of the terminal segment is monotone between its
endpoint values, therefore for every
        $s\in[t-h-d(h),t-e(h)],$
\[
        \sup_{s\in[t-h-d(h),t-e(h)]}
        |w_2(s)+v(x,t)|
        \longrightarrow0. 
\]
Now let us return to the first coordinate. Because the path segment is in the parallelogram defined by $(x_h,y_h)$ and $(x,-v(x,t))$, we have for 
$ s\in[t-h-d(h),t-e(h)]$
\[
        |w_1(s)-x|
        \le
        |x_h-x| + |y_h +v(x,t)|
        \longrightarrow0.
\]
Take a supremum over $s$ to conclude.
\end{proof}

\subsection{Envelope-selected Viscosity Solution to the PDE}

Using Lemma \ref{lem:pathapprox} we can now argue that $v(x,t)$ is a viscosity solution to the PDE 
\be \label{eq:ThePDE}
v_t + \tilde c(x,-v(x,t))(v_x(1-v_x)) = 0, 
\ee
according to Definition \ref{defn: classical viscosity solution}.
To this end, let $v(x,t)$ be given by the law of large numbers \eqref{eq:vdef}.
Let $\psi$ be the negative Legendre dual of $f$ from \eqref{Legendre dual relation}, so that
\begin{equation}\label{eq:duality-cf}
\tilde c\, f(p)=\inf_{\xi\in\mathbb{R}}
\left\{\xi p+\tilde c\,\psi\!\left(\frac{\xi}{\tilde c}\right)\right\}
\qquad (\tilde c>0,\ p\in\mathbb{R}),
\end{equation}
where we recall that $H_{\mathrm{sel}}(x, v(x,t), \phi_x (x,t)) = \tilde c(x,-v) f(\phi_x)$ for any functions $v,\phi$.

\begin{proof}[Proof of Theorem \ref{thm:visco}]

First pick a point $(x_0,t_0)$. 

We begin with the viscosity sub-solution part.

Now pick any $C^1$ function $\phi$ for which $v-\phi$ has a local maximum at $(x_0,t_0)$. After subtracting a constant we may assume
\[
(v-\phi)(x_0,t_0)=0,
\qquad
v(x,t)\le \phi(x,t).
\]
for all $(x,t)$ in a neighborhood of $(x_0,t_0)$.

Since $v$ is nondecreasing and Lipschitz-$1$ in $x$, any touching test function
has slope in $[0,1]$ at the touching point. Indeed, for $h>0$ small, from
$v-\phi$ having a local maximum at $(x_0,t_0)$ we obtain
\[
v(x_0+h,t_0)-v(x_0,t_0)
\le \phi(x_0+h,t_0)-\phi(x_0,t_0),
\]
and
\[
v(x_0-h,t_0)-v(x_0,t_0)
\le \phi(x_0-h,t_0)-\phi(x_0,t_0).
\]
Using the inequalities above and since $v$ is Lipschitz-$1$ in space, 
\[
0\le \phi(x_0+h,t_0)-\phi(x_0,t_0),
\qquad
-h\le \phi(x_0-h,t_0)-\phi(x_0,t_0),
\]
and dividing by $h$ and $-h$ respectively, then letting $h\downarrow 0$, yields
$
0\le \phi_x(x_0,t_0)\le 1.
$

Now choose admissible paths ${\bf w}^{h,\xi}$ as in Lemma \ref{lem:pathapprox}. Then we can write
\[
v(x_0,t_0)\le
v\bigl(x_0 - h\xi_h,t_0-h-d(h)\bigr)
-\int^{t_0-e(h)}_{t_0-h-d(h)} \tilde c({\bf w}^{h,\xi}(s))
\psi\!\left(\frac{w_1'(s)}{\tilde c({\bf w}^{h,\xi}(s))}\right)\,ds.
\]
Since $v\le \phi$ near $(x_0,t_0)$, for $h$ small enough,
\[
v(x_0,t_0)\le
\phi\bigl(x_0-h\xi_h,t_0-h-d(h)\bigr)
-\int^{t_0-e(h)}_{t_0-h-d(h)} \tilde c({\bf w}^{h,\xi}(s))
\psi\!\left(\frac{w_1'(s)}{\tilde c({\bf w}^{h,\xi}(s))}\right)\,ds.
\]
Using $v(x_0,t_0)=\phi(x_0,t_0)$, we obtain
\begin{equation}\label{eq:pre-taylor-visc}
\phi(x_0,t_0)\le
\phi\bigl(x_0-\xi_h h,t_0-h-d(h)\bigr)
-\int_{t_0-h-d(h)}^{t_0-e(h)} \tilde c({\bf w}^{h,\xi}(s))
\psi\!\left(\frac{w_1'(s)}{\tilde c({\bf w}^{h,\xi}(s))}\right)\,ds.
\end{equation}
 The $C^1$ expansion of $\phi$ gives
\[
\phi\bigl(x_0-\xi_h h,t_0-h-d(h)\bigr)
=
\phi(x_0,t_0)-h\bigl(\phi_t(x_0,t_0)+\xi_h \phi_x(x_0,t_0)\bigr)+o(h).
\]
Substituting into \eqref{eq:pre-taylor-visc},
\[
0\le
-h\bigl(\phi_t(x_0,t_0)+\xi_h \phi_x(x_0,t_0)\bigr)
-\int_{t_0-h-d(h)}^{t_0-e(h)}  \tilde c({\bf w}^{h,\xi}(s))
\psi\!\left(\frac{w_1'(s)}{\tilde c({\bf w}^{h,\xi}(s))}\right)\,ds
+o(h).
\]
Divide by $h$ to get:
\begin{equation}\label{eq:after-divide-visc}
\phi_t(x_0,t_0)+\xi_h \phi_x(x_0,t_0)
+\frac1h\int_{t_0-h-d(h)}^{t_0-e(h)}  \tilde c({\bf w}^{h,\xi}(s))
\psi\!\left(\frac{w_1'(s)}{\tilde c({\bf w}^{h,\xi}(s))}\right)\,ds
\le o(1).
\end{equation}
From \eqref{eq:after-divide-visc}, and using
\[
        x_0-x_h=h\xi_h
        =
        \int_{t_0-h-d(h)}^{t_0-e(h)} w_1'(s)\,ds,
\]
we may rewrite the inequality as
\begin{equation}
\label{eq:subsolution-pre-fenchel}
        \phi_t(x_0,t_0)
        +
        \frac1h
        \int_{t_0-h-d(h)}^{t_0-e(h)}
        \left[
        \phi_x(x_0,t_0)w_1'(s)
        +
        \tilde c({\bf w}^{h,\xi}(s))
        \psi\!\left(
        \frac{w_1'(s)}{\tilde c({\bf w}^{h,\xi}(s))}
        \right)
        \right]ds
        \le o(1).
\end{equation}
Set
\[
        p:=\phi_x(x_0,t_0) \in [0,1],
        \qquad
        c_0:=\tilde c(x_0,-v(x_0,t_0)).
\]
Recall the Fenchel relation
\[
        cH(p)
        =
        \inf_{\xi\in\mathbb R}
        \left\{
        \xi p+c\psi\!\left(\frac{\xi}{c}\right)
        \right\},
        \qquad
        H(p)=p(1-p),
\]
valid for $0\le p\le1$. Hence, pointwise in $s$,
\[
        \phi_x(x_0,t_0)w_1'(s)
        +
        \tilde c({\bf w}^{h,\xi}(s))
        \psi\!\left(
        \frac{w_1'(s)}{\tilde c({\bf w}^{h,\xi}(s))}
        \right)
        \ge
        \tilde c({\bf w}^{h,\xi}(s))H(p).
\]
Substituting this into \eqref{eq:subsolution-pre-fenchel}, we obtain
\begin{equation}
\label{eq:subsolution-liminf-bound}
        \phi_t(x_0,t_0)
        +
        \frac1h
        \int_{t_0-h-d(h)}^{t_0-e(h)}
        \tilde c({\bf w}^{h,\xi}(s))H(p)\,ds
        \le o(1).
\end{equation}
We now use only the lower semicontinuity of $\tilde c$. By the
construction in Lemma \ref{lem:pathapprox}, the terminal segment
\[
        \{{\bf w}^{h,\xi}(s):s\in[t_0-h-d(h),t_0-e(h)]\}
\]
shrinks to the point $(x_0,-v(x_0,t_0))$ as $h\downarrow0$ by Lemma \ref{lem:pathapprox}-$(4)$. Therefore, for
every $\varepsilon>0$, there exists $h_\varepsilon>0$ such that, for all
$0<h<h_\varepsilon$ and all
$s\in[t_0-h-d(h),t_0-e(h)]$,
\[
        \tilde c({\bf w}^{h,\xi}(s))\ge c_0-\varepsilon.
\]
Since $H(p)\ge0$, this gives
\[
        \frac1h
        \int_{t_0-h-d(h)}^{t_0-e(h)}
        \tilde c({\bf w}^{h,\xi}(s))H(p)\,ds
        \ge
        \frac{h+d(h)-e(h)}{h}(c_0-\varepsilon)H(p).
\]
As $d(h),e(h)=o(h)$,
\[
        \liminf_{h\downarrow0}
        \frac1h
        \int_{t_0-h-d(h)}^{t_0-e(h)}
        \tilde c({\bf w}^{h,\xi}(s))H(p)\,ds
        \ge
        (c_0-\varepsilon)H(p).
\]
Letting $\varepsilon\downarrow0$, we obtain
\[
        \liminf_{h\downarrow0}
        \frac1h
        \int_{t_0-h-d(h)}^{t_0-e(h)}
        \tilde c({\bf w}^{h,\xi}(s))H(p)\,ds
        \ge
        c_0H(p).
\]
Combining this with \eqref{eq:subsolution-liminf-bound} yields
\[
        \phi_t(x_0,t_0)+c_0H(p)\le0.
\]
That is,
\[
        \phi_t(x_0,t_0)
        +
        \tilde c(x_0,-v(x_0,t_0))
        \phi_x(x_0,t_0)\bigl(1-\phi_x(x_0,t_0)\bigr)
        \le0.
\]
This is the selected viscosity subsolution inequality of Definition \ref{defn: classical viscosity solution}.
\[
\phi_t(x_0,t_0)+H_{\rm low}\bigl(x_0, v(x_0,t_0),\phi_x(x_0,t_0)\bigr)\le 0,
\]

For the viscosity super-solution, we use \eqref{eq:lb4visco}. Start from any point $(x_0,-v(x_0,t_0))$.
Pick any $\xi \in \R$ and choose any path that is admissible for \eqref{eq:lb4visco}.
Since $v-\phi$ has a local minimum at $(x_0,t_0)$, after subtracting a constant we may assume
\[
(v-\phi)(x_0,t_0)=0,
\qquad
v(x,t)\ge \phi(x,t).
\]
for all $(x,t)$ in a neighborhood of $(x_0,t_0)$ and as in the previous case, the touching test functions have $\phi_x(x_0,t_0) \in [0,1]$. By \eqref{eq:lb4visco} we can write
\[
v(x_0,t_0)\ge
v\bigl(x_0 - h\xi,t_0-h\bigr)
-\int^{h}_{0} \tilde c({\bf x}(s))
\psi\!\left(\frac{x_1'(s)}{\tilde c({\bf x}(s))}\right)\,ds.
\]
Since $v\ge \phi$ near $(x_0,t_0)$, for $h$ small enough, and by using $v(x_0,t_0)=\phi(x_0,t_0)$ we have
\be\label{eq:pre-taylor-visc0}
\phi(x_0,t_0)\ge
\phi\bigl(x_0-h\xi,t_0-h\bigr)
-\int^{h}_{0} \tilde c({\bf x}(s))
\psi\!\left(\frac{x_1'(s)}{\tilde c({\bf x}(s))}\right)\,ds.
\ee

 The $C^1$ expansion of $\phi$ gives
\[
\phi\bigl(x_0-\xi h,t_0-h\bigr)
=
\phi(x_0,t_0)-h\bigl(\phi_t(x_0,t_0)+\xi \phi_x(x_0,t_0)\bigr)+o(h).
\]
Substituting into \eqref{eq:pre-taylor-visc0},
\[
0\ge
-h\bigl(\phi_t(x_0,t_0)+\xi \phi_x(x_0,t_0)\bigr)
-\int^{h}_{0} \tilde c({\bf x}(s))
\psi\!\left(\frac{x_1'(s)}{\tilde c({\bf x}(s))}\right)\,ds
+o(h).
\]
Divide by $h$:
\begin{equation}\label{eq:after-divide-visc3}
\phi_t(x_0,t_0)+\xi \phi_x(x_0,t_0)
+\frac1h\int^{h}_{0} \tilde c({\bf x}(s))
\psi\!\left(\frac{x_1'(s)}{\tilde c({\bf x}(s))}\right)\,ds
\ge o(1).
\end{equation}
The inequality above is over an arbitrary path. Specialise to an admissible short path
$\bf x$ with
$x_{1}(0)=x_0-\xi h,\,x_{1}(h)=x_0.$
For example, take
\[
        x_{1}(s)=x_0-\xi h+\xi s,
\]
and define $x_{2}$ by the admissibility relation
\[
        x'_{2}(s)
        =
        \tilde c({\bf x}(s))
        \psi\!\left(
        \frac{\xi}{\tilde c({\bf x}(s))}
        \right),
        \qquad
        x_{2}(0)
        =
        -v(x_0-\xi h,t_0-h).
\]
Note that for this path $
 \xi = x_1'(s).$
 
Set
\[
        p:=\phi_x(x_0,t_0),
        \qquad
        c_0:=\tilde c^*(x_0,-v(x_0,t_0)).
\]
The path $\bf x$ remains in an $o(1)$-neighbourhood of
$(x_0,-v(x_0,t_0))$. Since $\tilde c$ is lower semicontinuous, for
every $\varepsilon>0$, and all sufficiently small $h$,
\[
        \tilde c({\bf x}(s))\le c_0+\varepsilon,
        \qquad 0\le s\le h.
\]
The function $c\mapsto L(c,\xi) = c \psi(\xi/c)$ is non-decreasing in $c$, therefore
\[
        \tilde c({\bf x}(s))
        \psi\!\left(
        \frac{\xi}{\tilde c({\bf x}(s))}
        \right)
        =
        L(\tilde c({\bf x}(s)),\xi)
        \le
        L(c_0+\varepsilon,\xi).
\]
Substitute in the integral in \eqref{eq:after-divide-visc3} to obtain
\[
        \phi_t(x_0,t_0)+\xi p+L(c_0+\varepsilon,\xi)\ge0.
\]
Letting $\varepsilon\downarrow0$, we get
\[
        \phi_t(x_0,t_0)
        +
        \xi p
        +
        c_0\psi\!\left(\frac{\xi}{c_0}\right)
        \ge0.
\]
Since $\xi\in\mathbb R$ was arbitrary, taking the infimum over $\xi$ gives
\[
        \phi_t(x_0,t_0)
        +
        c_0 p(1-p)
        \ge0.
\]
This is precisely
\[
        \phi_t(x_0,t_0)
        +
       H_{\rm up}(x_0,v(x_0,t_0),\phi_x(x_0,t_0))
        \ge0. \qedhere
\]
\end{proof}

\section{Uniqueness of solutions in the absence of temporal discontinuities}
\label{sec:uniqueness}

In this section we prove the uniqueness part of Theorem \ref{thm:uniqueness} in the case where
the TASEP speed is spatially inhomogeneous but has no height dependence. Thus
we assume throughout this section that
\[
        \tilde c(x,y)=\tilde c(x).
\]
The Hamilton--Jacobi equation becomes
\[
        v_t+\tilde c(x)v_x(1-v_x)=0,
\]
and the associated scalar conservation law is
\[
        \rho_t+\bigl(\tilde c(x)\rho(1-\rho)\bigr)_x=0.
\]
We write
\[
        f(\rho)=\rho(1-\rho),\qquad
        G(x,\rho)=\tilde c(x)f(\rho),\qquad 0\leq \rho\leq 1.
\]

The viscosity notion used in this section remains the selected one from
Definition \ref{defn: classical viscosity solution}. In the present spatial setting this means that the
Hamiltonian is evaluated with the lower-semicontinuous representative of
$\tilde c$:
\[
        H_{\mathrm{low}}(x,r,p)
        =
        \tilde c(x)p(1-p) = G(x,p),
        \qquad 0\leq p\leq 1.
\]
Since $\tilde c$ is assumed lower semicontinuous in the spatial
case, this is simply the coefficient $\tilde c(x)$ itself. 
Also note that in this setup, Assumption \ref{ass:c0strong} is automatically satisfied, and $v(x,t)$ is an a.e. classical solution. 

For the comparison argument, however, we shall use a cone-determination
result of \cite{ChenHu2008}. That result is formulated in terms of the essential
upper and lower semicontinuous envelopes of the coefficient. Therefore, for a
locally bounded measurable spatial coefficient $\tilde c$, identified
with the space--time function $(x,t)\mapsto \tilde c(x)$, define
\[
\tilde c^{*}({\bf z}):=\lim_{\delta\downarrow0}
        \operatorname*{ess\,sup}_{B_\delta(\bf z)}\tilde c,
        \qquad
\tilde c_{*}({\bf z}):=\lim_{\delta\downarrow0}
        \operatorname*{ess\,inf}_{B_\delta({\bf z})}\tilde c,
        \qquad {\bf z}=(x,t)\in \mathbb R\times(0,\infty),
\]
where $B_\delta({\bf z})$ denotes the Euclidean ball in $\R^2$ centered at ${\bf z}$.
 
 Since $\tilde c$ is independent of $t$, these envelopes reduce to the
corresponding one-dimensional essential envelopes in the spatial variable:
\[
\tilde c^{*}(x,t)=\tilde c^{*}(x),\qquad
\tilde c_{*}(x,t)=\tilde c_{*}(x).
\]
 
Under our lower-semicontinuity convention,
\[
        \tilde c_{*}(x)=\tilde c(x),
\]
while $\tilde c^{\,*}$ records the possible upper value across a spatial
jump. Thus the selected Hamiltonian uses the lower representative, whereas
$\tilde c^{\,*}$ enters only through the comparison criterion below.

 Moreover, since  $\tilde c(x)$ is locally piecewise continuous,
with only locally finitely many jump discontinuities, we see that $\tilde c$ satisfies the \emph{monotonicity property at $x_*$} as defined in \cite{ChenHu2008} since there exists a sign
$\sigma=\sigma(x_*)\in\{-1,+1\}$ such that for every $M>1$,
\begin{equation}\label{eq:mono-1d}
\lim_{\rho\downarrow 0}\ 
\sup_{\substack{|y-x_*|\le M\rho}}
\Big[\tilde c^{*}(y)-\tilde c_{*}(y-\sigma\rho)\Big]\ \le\ 0.
\end{equation}
Away from discontinuities this condition is immediate, since
$\tilde c^{\,*}=\tilde c_{*}=\tilde c$ locally.

Fix ${\bf z}_0=(x_0,t_0)\in\R\times(0,\infty)$ and a constant $K>0$.
The backward cone with vertex ${\bf z}_0$ and slope $K$ is
\be\label{eq:cone-def}
C({\bf z}_0,K):=C(x_0,t_0,K)
:=\Big\{(x,t): 0\le t\le t_0,\ |x-x_0|\le K(t_0-t)\Big\}.
\ee

The next lemma is the one-dimensional version of Chen--Hu's cone-determination
Lemma 8 \cite{ChenHu2008}, rewritten for the sign convention of \eqref{eq:no-t-HJ}.

\begin{lemma}[Cone determination for \eqref{eq:no-t-HJ}]
\label{lem:cone-determination-1d}
Let $H:\mathbb R\to [0,\infty)$ satisfy
\[
|H(p)-H(q)|\le |p-q| \qquad \forall p,q\in \mathbb R,
\]
and assume in addition that
\[
H(1-p)=H(p)\qquad \forall p\in \mathbb R.
\]
Let $\tilde c:\mathbb R\to [0,\infty)$, so that $ \tilde c\in L^{\infty}_{\rm loc}(\R)$. Fix
${\bf z}_0=(x_0,t_0)\in \mathbb R\times (0,\infty)$ and $K>1$ such that
\begin{equation}
\label{eq:K-dominates}
K>\sup_{C({\bf z}_0,K)} \tilde c^{*}.
\end{equation}
Assume that for every ${\bf z}^*=(x^*,t^*)\in C({\bf z}_0,K)$ there exists a sign
$\sigma=\sigma(z^*)\in\{-1,+1\}$ such that, for all $M>1$,
\begin{equation}
\label{eq:cone-mono}
\lim_{\rho\downarrow 0}
\sup_{\substack{|y-x^*|\le K(t^*-\tau)\le M\rho}}
\bigl[\tilde c^{*}(y)-\tilde c_{*}(y-\sigma\rho)\bigr]\le 0.
\end{equation}

Let $v_1,v_2\in C(C({\bf z}_0,K))$ (the continuous functions on the cone) be such that, on $C({\bf z}_0,K)$,
\begin{enumerate}
\item $v_1$ is an envelope-selected viscosity subsolution of
$
v_t+\tilde c(x)H(v_x)=0,
$
\item $v_2$ is a envelope-selected viscosity supersolution of
$
v_t+\tilde c(x)H(v_x)=0,
$
\item both $v_1$ and $v_2$ are nondecreasing and Lipschitz-$1$ in the $x$-variable,
\item at least one of $v_1,v_2$ is Lipschitz continuous on $C({\bf z}_0,K)$.
\end{enumerate}
Then
\begin{equation}
\label{eq:cone-max}
\max_{C({\bf z}_0,K)}(v_1-v_2)
=
\max_{C({\bf z}_0,K)\cap\{t=0\}}(v_1-v_2).
\end{equation}
\end{lemma}

\begin{proof}
Define
\[
u_1(x,t):=x-v_2(x,t),
\qquad
u_2(x,t):=x-v_1(x,t).
\]
Then
\[
u_1-u_2=v_1-v_2.
\]
Also, if one of $v_1,v_2$ is Lipschitz on $C({\bf z}_0,K)$, then the corresponding
$u_i$ is Lipschitz as well. Moreover, since $v_1$ and $v_2$ are nondecreasing and
Lipschitz-1 in $x$, the same is true of $u_1$ and $u_2$.

We claim that $u_1$ is a viscosity \emph{subsolution} and $u_2$ is a viscosity
\emph{supersolution} of
\begin{equation}
\label{eq:positive-sign-HJ}
u_t=\tilde c(x)H(u_x)
\qquad\text{on } C({\bf z}_0,K),
\end{equation}
in the sense of \cite{ChenHu2008}.
Here the coefficient $\tilde c$ is understood in the following sense: subsolutions
are tested with $\tilde c^*$, while supersolutions are tested with $\tilde c_*$.

First let $\psi\in C^1$, and suppose that $u_1-\psi$ has a local maximum
at $(x_*,t_*)$. Set
\[
        \phi(x,t):=x-\psi(x,t).
\]
Then
\[
        v_2-\phi=-(u_1-\psi),
\]
so $v_2-\phi$ has a local minimum at $(x_*,t_*)$. Since $v_2$
is a selected supersolution,
\[
        \phi_t(x_*,t_*)
        +
        \tilde c^*(x_*) H(\phi_x(x_*,t_*))
        \ge0.
\]
Using
\[
        \phi_t=-\psi_t,
        \qquad
        \phi_x=1-\psi_x,
\]
and the symmetry $ H(1-p)= H(p)$, we get
\[
        -\psi_t(x_*,t_*)
        +
        \tilde c^*(x_*) H(\psi_x(x_*,t_*))
        \ge0.
\]
Hence
\[
        \psi_t(x_*,t_*)
        \le
        % \tilde c_*(x_*) H(\psi_x(x_*,t_*))
        % \le
        \tilde c^*(x_*) H(\psi_x(x_*,t_*)).
\]
Thus $u_1$ is a subsolution in the Chen--Hu sense.

Now let $\psi\in C^1$, and suppose that $u_2-\psi$ has a local minimum at
$(x_*,t_*)$. Again set
\[
        \phi(x,t):=x-\psi(x,t).
\]
Then
\[
        v_1-\phi=-(u_2-\psi),
\]
so $v_1-\phi$ has a local maximum at $(x_*,t_*)$. Since $v_1$
is a selected subsolution,
\[
        \phi_t(x_*,t_*)
        +
        \tilde c_*(x_*) H(\phi_x(x_*,t_*))
        \le0.
\]
Therefore
\[
        -\psi_t(x_*,t_*)
        +
        \tilde c_*(x_*) H(1-\psi_x(x_*,t_*))
        \le0,
\]
and by the symmetry of $H$,
\[
        \psi_t(x_*,t_*)
        \ge
        \tilde c_*(x_*) H(\psi_x(x_*,t_*)).
\]
Thus $u_2$ is a supersolution in the Chen--Hu sense.

We may therefore apply the cone-determination Lemma~8 of \cite{ChenHu2008} to $u_1$ and
$u_2$, with $\tilde c$ as the coefficient. The monotonicity hypothesis in Lemma~8 is
precisely equation \eqref{eq:mono-1d}. 
Hence
\[
        \max_{C(z_0,K)}(u_1-u_2)
        =
        \max_{C(z_0,K)\cap\{t=0\}}(u_1-u_2).
\]
Since $u_1-u_2=v_1-v_2$, this is precisely \eqref{eq:cone-max}.
\end{proof}

\begin{remark} In \cite{ChenHu2008} there is also the assumption that $\tilde c$ is bounded. Indeed, this is not necessary and local boundedness is enough, as only the restrictions of $v$ to compact cones are used. 
\qed \end{remark}

The uniqueness of the viscosity solution is essentially contained in the previous statement. We just need to be careful to satisfy the assumptions of Lemma \ref{lem:cone-determination-1d}.

\begin{proof}[Proof of Theorem \ref{thm:uniqueness}-HJ part] 
We first replace the Hamiltonian by a globally Lipschitz, non-negative extension on all
of $\R$.
To this end,
define the projection
\begin{equation}
\label{eq:proj}
\Pi(p):=\min\{1,\max\{0,p\}\},
\qquad p\in \mathbb R,
\end{equation}
and the clamped Hamiltonian
\begin{equation}
\label{eq:H-clamped}
\bar H(p):=H(\Pi(p))=\Pi(p)\bigl(1-\Pi(p)\bigr) = \{p(1-p)\}\vee0,
\qquad p\in \mathbb R.
\end{equation}
Then $\bar H\equiv H$ on $[0,1]$, $\bar H(p)\in [0,1/4]$ for all $p$, $\bar H$ is globally
Lipschitz, and
\begin{equation}
\label{eq:H-symmetry}
\bar H(1-p)=\bar H(p)\qquad \forall p\in \mathbb R.
\end{equation}

Let $v^{(1)},v^{(2)}$ be two viscosity solutions of \eqref{eq:no-t-HJ} with the same
initial datum $v_0$, and assume that both are nondecreasing and Lipschitz-$1$ in the
space variable and at least one of them is Lipschitz on compact subsets of $\mathbb R\times [0,\infty)$.
These are the natural bounds satisfied by the TASEP current from
Theorem~\ref{thm:LLN}.

We first note that each $v^{(i)}$ is also an envelope-selected viscosity solution of the clamped equation
\begin{equation}
\label{eq:HJ-clamped}
v_t+\tilde c(x)\bar H(v_x)=0,
\qquad (x,t)\in \mathbb R\times (0,\infty).
\end{equation}
Indeed, if $\phi$ touches $v^{(i)}$ from above or below at some point $(x_*,t_*)$, then
the same argument as in the proof of Theorem~\ref{thm:visco} gives
\[
0\le \phi_x(x_*,t_*)\le 1.
\]
Hence at every touching point,
\[
\bar H(\phi_x(x_*,t_*))=H(\phi_x(x_*,t_*)),
\]
so the viscosity inequalities for \eqref{eq:no-t-HJ} and \eqref{eq:HJ-clamped} coincide.

Now fix an arbitrary point ${\bf z}_0=(x_0,t_0)\in \mathbb R\times (0,\infty)$ and choose
$K>1$ so large that \eqref{eq:K-dominates} holds on $C({\bf z}_0,K)$. Since $\tilde c$ is
lower semicontinuous, locally piecewise continuous, and depends only on $x$, the monotonicity condition
\eqref{eq:cone-mono} holds on this cone. Applying
Lemma~\ref{lem:cone-determination-1d} to $v_1=v^{(1)}$ and $v_2=v^{(2)}$ with
$H=\bar H$, we obtain
\[
\max_{C({\bf z}_0,K)}\bigl(v^{(1)}-v^{(2)}\bigr)
=
\max_{C({\bf z}_0,K)\cap\{t=0\}}\bigl(v^{(1)}-v^{(2)}\bigr).
\]
Since $v^{(1)}(\cdot,0)=v^{(2)}(\cdot,0)=v_0$, the right-hand side is $0$. Thus
\[
v^{(1)}\le v^{(2)} \qquad \text{on } C({\bf z}_0,K).
\]
Because ${\bf z}_0$ was arbitrary, this holds on all of $\mathbb R\times (0,\infty)$ and $t=0$ equality holds from the initial conditions.
Exchanging the roles of $v^{(1)}$ and $v^{(2)}$ yields the reverse inequality.
Hence
\[
v^{(1)}\equiv v^{(2)}
\qquad \text{on } \mathbb R\times [0,\infty).
\]
This suffices for the uniqueness of the current $v(x,t)$ as a envelope-selected viscosity solution of the Hamilton--Jacobi equation, since by Theorem \ref{thm:LLN} it is locally Lipschitz on $\R\times[0,\infty)$. 
\end{proof}

We now turn to the scalar conservation law. The Hamilton--Jacobi comparison
argument gives uniqueness of the integrated current in the selected viscosity
class. To obtain the corresponding statement for particle densities, we use
the relation $\rho=v_x$ and show that the TASEP density is the weak solution
selected by maximal time-integrated current. The argument follows the same maximal current principle as in \cite{Seppalainen1999KExclusion}, with the proof adjusted to the current setting.

\begin{proof}[Proof of Theorem \ref{thm:uniqueness}-SCL part]
Let $\lambda:\mathbb R\times \mathbb R_+\to [0,1]$ be a locally bounded measurable weak solution of
\[
\lambda_t + \bigl(G(x,\lambda)\bigr)_x = 0
\]
in the sense of \eqref{weak soln integral form}. We will prove that for every fixed $t>0$,
\[
\int_0^t G(x,\lambda(x,s))\,ds
\leq
\int_0^t G(x,\rho(x,s))\,ds
\qquad\text{for a.e.\ }x\in\mathbb R,
\]
where $\rho=v_x$ and $v$ is the TASEP current from Theorem \ref{thm:LLN}.

We begin by recording two consequences of the weak formulation. First, let
\[
\varphi(x,t)=h(x)g_\varepsilon(t),
\]
where $h\in C_c^1(\mathbb R)$ and $g_\varepsilon$ are chosen so that 
\[
g_\varepsilon(0)=1,
\qquad
g_\varepsilon'(t)=-\frac1\varepsilon \ \text{on }(0,\varepsilon),
\qquad
g_\varepsilon(t)=0 \ \text{for }t\geq \varepsilon.
\]
Substituting $\varphi$ into \eqref{weak soln integral form} gives
\begin{equation}\label{eq:5.1-new}
\int_0^\varepsilon \frac1\varepsilon \int_{\mathbb R}
\Bigl(\lambda(x,t)h(x)+\varepsilon G(x,\lambda(x,t))h'(x)g_\varepsilon(t)\Bigr)\,dx\,dt
=
\int_{\mathbb R}\rho_0(x)h(x)\,dx.
\end{equation}
Since $h'$ is compactly supported and $G(x,\lambda(x,t))$ is locally bounded, the second term vanishes as
$\varepsilon\downarrow 0$, and therefore
\begin{equation}\label{eq:5.2-new}
\lim_{\varepsilon\downarrow 0}
\int_0^\varepsilon \frac1\varepsilon \int_{\mathbb R}\lambda(x,t)h(x)\,dx\,dt
=
\int_{\mathbb R}\rho_0(x)h(x)\,dx.
\end{equation}

Next define
\[
\Leb(\lambda)=\{\text{Lebesgue points of }\lambda\}
\]
and
\begin{equation}\label{eq:T-new}
T=\{t>0:(x,t)\in \Leb(\lambda)\text{ for a.e.\ }x\in\mathbb R\}.
\end{equation}
Since the complement of $\Leb(\lambda)$ has zero two-dimensional Lebesgue measure, $T$ is dense and has full measure in
$(0,\infty)$.

Fix $t\in T$ and choose
\[
\varphi(x,s)=h(x)g_{t,\varepsilon}(s),
\]
where
\[
g_{t,\varepsilon}(s)=1 \ \text{for }0\le s\le t,
\qquad
g_{t,\varepsilon}'(s)=-\frac1\varepsilon \ \text{on }(t,t+\varepsilon),
\qquad
g_{t,\varepsilon}(s)=0 \ \text{for }s\ge t+\varepsilon.
\]
Substituting into \eqref{weak soln integral form}, letting $\varepsilon\downarrow 0$, and using that $t\in T$, we obtain
\begin{equation}\label{eq:5.6-new}
\int_{\mathbb R} h(x)\lambda(x,t)\,dx
-
\int_{\mathbb R} h(x)\rho_0(x)\,dx
=
\int_0^t\!\!\int_{\mathbb R} G(x,\lambda(x,s))h'(x)\,dx\,ds.
\end{equation}
By approximation, \eqref{eq:5.6-new} remains valid for compactly supported piecewise $C^1$ functions $h$.

We now prove the time-regularity estimate used later.

\begin{lemma}\label{lem:5.1-new}
\leavevmode
\begin{enumerate}
\item
Fix finite $A<B$ and $\tau>0$. Then there exists a constant $C=C(A,B,\tau)$ such that for all
$A\le a<b\le B$ and all $s,t\in T\cup\{0\}$ with $0\le s,t\le \tau$,
\begin{equation}\label{eq:5.7-new}
\left|\int_a^b \lambda(x,t)\,dx-\int_a^b \lambda(x,s)\,dx\right|
\le C|t-s|.
\end{equation}

\item
If $t\in T\cup\{0\}$ and $(a,t),(b,t)\in {\rm Leb}(G\circ \lambda)$, then
\begin{equation}\label{eq:5.8-new}
\lim_{s\to t,\ s\in T}
\frac{1}{t-s}
\left(
\int_a^b \lambda(x,t)\,dx-\int_a^b \lambda(x,s)\,dx
\right)
=
G(a,\lambda(a,t))-G(b,\lambda(b,t)).
\end{equation}
\end{enumerate}
\end{lemma}

\begin{proof}
Choose $h\in C_c^1(\mathbb R)$ such that $h=0$ outside $[a-\varepsilon,b+\varepsilon]$, $h'=1/\varepsilon$ on
$(a-\varepsilon,a)$, $h=1$ on $[a,b]$, and $h'=-1/\varepsilon$ on $(b,b+\varepsilon)$. Apply
\eqref{eq:5.6-new} first at time $t$ and then at time $s$, and subtract. This gives
\[
\int_a^b \lambda(x,t)\,dx-\int_a^b \lambda(x,s)\,dx + O(\varepsilon)
=
\frac1\varepsilon \int_s^t\!\!\int_{a-\varepsilon}^a G(x,\lambda(x,\tau))\,dx\,d\tau
-
\frac1\varepsilon \int_s^t\!\!\int_b^{b+\varepsilon} G(x,\lambda(x,\tau))\,dx\,d\tau.
\]
Since $\lambda$ takes values in $[0,1]$ and $\tilde c$ is locally bounded, $G(x,\lambda(x,\tau))$ is locally
bounded on $[A,B]\times[0,\tau]$, and \eqref{eq:5.7-new} follows by letting $\varepsilon\downarrow 0$.

For \eqref{eq:5.8-new}, choose $\varepsilon=\delta|t-s|$, divide through by $t-s$, let $s\to t$ with
$s\in T$, and then let $\delta\downarrow 0$. The Lebesgue-point assumption at $(a,t)$ and $(b,t)$ yields
the stated limit.
\end{proof}

We now return to the proof of the theorem. Since $v$ is locally Lipschitz and solves
\[
v_t + G(x,v_x)=0
\]
for almost every $(x,t)$, with $v_x=\rho$ a.e., we have for almost every $x_0$ and almost every $t>0$,
\begin{equation}\label{eq:5.11-new}
\int_0^t G(x_0,\rho(x_0,s))\,ds
=
-\int_0^t v_t(x_0,s)\,ds
=
v_0(x_0)-v(x_0,t).
\end{equation}
Fix $x_0$ such that \eqref{eq:5.11-new} holds for almost every $t>0$, and such that
$(x_0,t)\in \Leb(G\circ\lambda)$ for almost every $t>0$.

For $(x,t)\in \mathbb R\times T$, define
\begin{equation}\label{eq:5.12-new}
u(x,t)
=
\int_{x_0}^x \lambda(y,t)\,dy
-
\int_0^t G(x_0,\lambda(x_0,s))\,ds
+
v_0(x_0).
\end{equation}
By Lemma~\ref{lem:5.1-new} and local boundedness of $\lambda$ and $G\circ\lambda$, the function $u$ is
locally Lipschitz on $\mathbb R\times T$, and therefore extends uniquely to a locally Lipschitz function on
$\mathbb R\times [0,\infty)$. Moreover, using bounded convergence, we obtain
\begin{equation}\label{eq:5.13-new}
\lim_{t\to 0+,\ t\in T} u(x,t)
=
\int_{x_0}^x \rho_0(y)\,dy + v_0(x_0)
=
v_0(x).
\end{equation}

We claim that
\begin{equation}\label{eq:5.14-new}
u(x,t)\ge v(x,t)\qquad \text{for all }(x,t)\in \mathbb R\times [0,\infty).
\end{equation}
Once this is proved, setting $x=x_0$ in \eqref{eq:5.12-new} and using \eqref{eq:5.11-new} gives, for every
$t\in T$,
\[
-\int_0^t G(x_0,\lambda(x_0,s))\,ds + v_0(x_0)
=
u(x_0,t)
\ge v(x_0,t)
=
v_0(x_0)-\int_0^t G(x_0,\rho(x_0,s))\,ds.
\]
Hence
\[
\int_0^t G(x_0,\lambda(x_0,s))\,ds
\le
\int_0^t G(x_0,\rho(x_0,s))\,ds
\qquad\text{for all }t\in T.
\]
Since both time integrals are continuous in $t$, the same inequality holds for all $t\ge 0$. As the choice of $x_0$
was valid for almost every spatial point, this proves \eqref{maximum principle}.

It remains to prove \eqref{eq:5.14-new}. For this, fix $(x,t)$ and let
\[
\bm \xi=(\xi_1,\xi_2)\in \mathcal H^{v_0}_{x,t}
\]
be any admissible path. Since $u$ is locally Lipschitz and $\bm \xi$ is piecewise $C^1$, the composition
\[
\gamma(s):=u(\xi_1(s),s),\qquad s\in [0,t],
\]
is absolutely continuous. By Rademacher's theorem and the chain rule for Lipschitz functions,
for almost every $s\in [0,t]$ at which $u$ is differentiable and $\xi$ is differentiable,
\begin{equation}\label{eq:gamma-chain}
\gamma'(s)=u_x(\xi_1(s),s)\,\xi_1'(s)+u_t(\xi_1(s),s).
\end{equation}

We identify the partial derivatives of $u$ almost everywhere. From \eqref{eq:5.12-new}, for almost every
$(x,t)$,
\begin{equation}\label{eq:ux-new}
u_x(x,t)=\lambda(x,t).
\end{equation}
Next, fix $x$ such that $(x,t)\in \Leb(G\circ\lambda)$ for almost every $t$, and consider
\[
u(x,t)-u(x,s)
=
\int_{x_0}^x \bigl(\lambda(y,t)-\lambda(y,s)\bigr)\,dy
-
\int_s^t G(x_0,\lambda(x_0,\tau))\,d\tau.
\]
Divide by $t-s$ and use Lemma~\ref{lem:5.1-new}(2) with $a=x_0$ and $b=x$ to conclude that, for almost
every $(x,t)$,
\begin{equation}\label{eq:ut-new}
u_t(x,t)
=
-G(x,\lambda(x,t)).
\end{equation}
Substituting \eqref{eq:ux-new} and \eqref{eq:ut-new} into \eqref{eq:gamma-chain} gives
\[
\gamma'(s)
=
\lambda(\xi_1(s),s)\,\xi_1'(s)-G(\xi_1(s),\lambda(\xi_1(s),s))
\]
for almost every $s$.

Now recall that
\[
G(x,r)=\tilde c(x)f(r),
\]
and that $\psi=-f^*$ is the negative Legendre transform of $f$. Hence, for every $x\in \mathbb R$,
every $p\in \mathbb R$, and every $r\in [0,1]$,
\[
rp-\tilde c(x)f(r)\ge -\tilde c(x)\psi\!\left(\frac{p}{\tilde c(x)}\right).
\]
Applying this with $r=\lambda(\xi_1(s),s)$ and $p=\xi_1'(s)$ yields
\[
\gamma'(s)
\ge
-\tilde c(\xi_1(s))\psi\!\left(\frac{\xi_1'(s)}{\tilde c(\xi_1(s))}\right).
\]
Since $\xi\in \mathcal H^{v_0}_{x,t}$, we have by definition of admissible paths that
\[
\xi_2'(s)=\tilde c(\xi_1(s))\psi\!\left(\frac{\xi_1'(s)}{\tilde c(\xi_1(s))}\right)
\qquad\text{for a.e.\ }s\in [0,t].
\]
Therefore,
\[
\gamma'(s)\ge -\xi_2'(s)
\qquad\text{for a.e.\ }s\in [0,t].
\]
Integrating from $0$ to $t$, and using $\gamma(t)=u(x,t)$ together with
$\xi_2(0)=-v_0(\xi_1(0))$, we obtain
\[
u(x,t)-u(\xi_1(0),0)
=
\gamma(t)-\gamma(0)
\ge
-\int_0^t \xi_2'(s)\,ds
=
-\xi_2(t)+\xi_2(0).
\]
Using \eqref{eq:5.13-new} at time $0$, this becomes
\[
u(x,t)
\ge
v_0(\xi_1(0))-\int_0^t \tilde c(\xi_1(s))\psi\!\left(\frac{\xi_1'(s)}{\tilde c(\xi_1(s))}\right)\,ds.
\]
Since this inequality holds for every admissible path $\bm \xi\in \mathcal H^{v_0}_{x,t}$, Proposition \ref{thm:lang} yields
\[
u(x,t)\ge v(x,t).
\]
This proves \eqref{eq:5.14-new}, and hence also \eqref{maximum principle}.

Finally, assume that equality holds in \eqref{maximum principle} for almost every $(x,t)$, namely
\begin{equation}\label{eq:5.28-new}
\int_0^t G(x,\lambda(x,s))\,ds
=
\int_0^t G(x,\rho(x,s))\,ds
\qquad\text{for a.e.\ }(x,t).
\end{equation}
For almost every $t>0$, both $\lambda$ and $\rho$ satisfy \eqref{eq:5.6-new}. Using
\eqref{eq:5.28-new}, Fubini's theorem, and the compact support of $h'$, we may replace
\[
\int_0^t\!\!\int_{\mathbb R} G(x,\lambda(x,s))h'(x)\,dx\,ds
\]
by the corresponding expression with $\rho$ in place of $\lambda$. Therefore, for almost every $t>0$,
\[
\int_{\mathbb R} h(x)\lambda(x,t)\,dx
=
\int_{\mathbb R} h(x)\rho(x,t)\,dx
\qquad\text{for all }h\in C_c^1(\mathbb R).
\]
Hence $\lambda(\cdot,t)=\rho(\cdot,t)$ almost everywhere for almost every $t$, and so
\[
\lambda(x,t)=\rho(x,t)
\qquad\text{for a.e.\ }(x,t)\in \mathbb R\times \mathbb R_+.
\]
This proves the final assertion.
\end{proof}

\appendix

\section{Generalised shape function for the space-time discontinuous LPP}
\label{app:LLN}

In this appendix we prove Proposition \ref{lem:SPEP}. We need a couple of preliminary lemmas. 

\begin{lemma}[Shift of the speed function]\label{lem:gammashift}
Suppose that we have a speed function $c(x,y) : \bR \times \bR \to \bR_{+}$ that satisfies Assumptions \ref{ass:c}, \ref{ass:c2}, \ref{ass:c0}. 
Fix $(a,b) \in \bR^{2}$ and consider the shifted $(T_{a,b}c)(s,t) = c(s+a, t+b)$.
Suppose that we then discretize both $c$ and $T_{a,b}c$ separately as 
 \[ c^{n}(i,j) = c\left(i/n, j/n \right) \text{ and } (T_{a,b}c)^n(i,j) = T_{a,b}c\left(i/n, j/n \right).\]

Denote by $G^{(n),c}$ and $G^{(n), (T_{a,b}c)}$ the last passage times using the speed functions $c$ and $T_{a,b}c$ respectively. We then have the scaling limits
\begin{equation*}
    \lim_{n \to \infty}n^{-1}G^{(n),c}_{{\bf 0},(\ce{nx}, \ce{ny})} = \Gamma_{c}(x,y), \hspace{15.pt} \bP \text{ - a.s,}
\end{equation*}
and 
\begin{equation*}
    \lim_{n \to \infty}n^{-1}G^{(n),T_{a,b}c}_{{\bf 0}, (\ce{nx}, \ce{ny})} = \Gamma_{T_{a,b}c}(x,y), \hspace{15.pt} \bP \text{ - a.s.}
\end{equation*}
Furthermore, we have that
\begin{equation}
    \Gamma_{T_{a,b}c}(x,y) = \Gamma_{c}((a,b),(x+a,y+b)), \hspace{15.pt} \bP \text{ - a.s.}
\end{equation}
\end{lemma}

\begin{proof}
The first two limits follow directly from the proof of Theorem 2.5 in \cite{ciech2021}; since if $c$ satisfies Assumptions \ref{ass:c}, \ref{ass:c2}, \ref{ass:c0}, it will still satisfy those assumptions after any horizontal or vertical shift. Let us then deal with the final equality. The goal is to show the equality
\begin{equation}
    \sup_{w(\cdot) \in \mathcal{H}(x,y)}\left\{\int^1_0 \frac{\gamma(w'(s))}{(T_{a,b}c)(w_{1}(s), w_2(s))} ds \right\} = \sup_{w(\cdot) \in \mathcal{H}((a,b)(x + a,y + b))}\left\{\int^1_0 \frac{\gamma(w'(s))}{c(w_{1}(s), w_2(s))} ds \right\}.
\end{equation}
To do this, we can find a simple mapping between the two sets of paths. For any path $w(\cdot) \in \mathcal{H}(x,y)$, we can map onto a path  $u(\cdot) \in \mathcal{H}((a, b),(x + a, y + b))$ by setting
\begin{equation*}
    u(s) = (w_{1}(s) + a, w_{2}(s) + b).
\end{equation*}
From this, $u'(s) = (w_1'(s), w_2'(s))$ and $c(u_1(s), u_2(s)) = c(w_1(s) + a, w_2(s) + b) = (T_{a,b}c)(w_1(s), w_2(s))$. Thus it is immediate that 
\begin{equation*}
    \int^1_0 \frac{\gamma(u'(s))}{c(u_1(s), u_2(s))} = \int^1_0 \frac{\gamma(w'(s))}{(T_{a,b}c)(w_1(s), w_2(s))}.
\end{equation*}
Since this equality is path-by-path, we also have equality in the supremum. 
\end{proof}

\begin{lemma}[Shift of the starting point] \label{lem:LPPshift}
	Assume $c(x,y)$ satisfies Assumptions \ref{ass:c}, \ref{ass:c2}, \ref{ass:c0} with discretisation $c^n= c(in^{-1}, jn^{-1})$. Then we have  the limit
	\[
	\lim_{n\to \infty} n^{-1}G^{(n)}_{(\fl{na}, \fl{nb}), (\fl{ns}, \fl{nt})} = \Gamma_c((a,b), (s,t)).
	\]
\end{lemma}

\begin{proof}
Without loss of generality, assume $(x,y) \in \R^2_+$. For any $\eta> 0$, define step functions which have the same discontinuity curves as $c(x,y)$ and potentially further discontinuities that can be either vertical or horizontal, so that for any $(s,t)$ in the rectangle $R_{(0,0), (x,y)} = \{ (s,t): (0,0) \le (s,t) \le (x,y)\}$ we have 
\[ c^{\eta}_{\text{ low}}(s,t) \le c(s,t) \le c^{\eta}_{\text{high}}(s,t) 
\quad \text{ and }\quad
\sup_{(s,t) \in R_{(0,0), (x,y)}} \Big\{ \frac{1}{c^{\eta}_{\text{low}}(s,t)} - \frac{1}{c^{\eta}_{\text{ high}}(s,t)}\Big\} \le \eta.\]
This guarantees that 
\be \label{eq:etanew}
|\Gamma_{c^{\eta}_{\text{low}}}(x,y) - \Gamma_{c}(x,y)| + |\Gamma_{c}(x,y)- \Gamma_{c^{\eta}_{\text{high}}}(x,y)| < C_{x,y} \eta. 
\ee

We discretise $c^{\eta}_{\text{low}}(x,y)$ to assign rates for the exponential variables using a fixed integer $K$ in the following way: 
\[
c^{\eta,K}_{\text{low}}(in^{-1},jn^{-1}) := c^{n,\eta,K}_{\text{low}}(i,j) := \min_{(s_n,t_n) \in [-n^{-1}K, n^{-1}K]^2} c^{\eta}_{\text{low}}\Big(\frac{i}{n}+s_n,\frac{j}{n}+t_n\Big) \le c^{n,\eta}_{\text{low}}(i,j) = c^{\eta}_{\text{low}}(in^{-1},jn^{-1}).
\]
Note that because $c^{\eta}_{\text{low}}$ is a step function, this discretisation affects values for the exponential rates only in a macroscopic neighbourhood of diameter $2Kn^{-1}$ around the finitely many discontinuity curves. If $(in^{-1}, jn^{-1})$ is more than $Kn^{-1}$ away from a discontinuity curve then
\[
 c^{\eta,K}_{\text{low}}(in^{-1},jn^{-1}) = c^{\eta}_{\text{low}}(in^{-1}, jn^{-1}).
\]
The proof of Theorem 2.4. in \cite{ciech2021} now follows for this discretisation, and the limiting LLN starting from 0 does not change due to the continuity of the limit $\Gamma$, i.e. 
\be\label{eq:diffapp}
\Gamma_{c^\eta_{\rm low}}(x,y) = \Gamma_{ c^{\eta,K}_{\rm low}}(x,y). 
\ee

Now, fix $(a,b) \in \bR^{2}$ and consider the shifted $T_{a,b}c(s,t) = c(s+a, t+b)$ for $(s,t) \in \R^2_+$. We have from Lemma \ref{lem:gammashift} that 
\[
\Gamma_{T_{a,b}c}(s,t) = \Gamma_{c}((a,b)(s+a, t+b)),
\]
and we define 
\be \label{A6}
(T_{a,b}c)^{\eta}_{\rm low} := T_{a,b}(c^{\eta}_{\rm low}).
\ee
but we want this to be reflected in the passage times. So consider the usual discretisation of $c^n(i,j) = c\big(\frac{i}{n}, \frac{j}{n}\Big)$ and note 
\begin{align*}
(T_{\fl{na},\fl{nb}} c)^{n}(i,j) = c^n(\fl{na} + i, \fl{nb} + j)& = c\Big(\frac{\fl{na} + i}{n}, \frac{\fl{nb} + j}{n}\Big)\\
&=c\Big(a + \frac{i}{n}+ \frac{s_n}{n}, b +\frac{j}{n}+\frac{t_n}{n}\Big) \ge T_{a,b}( c^{\eta, K}_{\text{low}})(in^{-1},jn^{-1}). 
\end{align*}
Moreover, 
\be \label{Shifdisc}
T_{a,b}( c^{\eta, K}_{\text{low}})(in^{-1},jn^{-1}) = ({T_{a,b} c})^{\eta, K}_{\text{low}}(in^{-1},jn^{-1})
\ee
This allows us to couple the exponential weights and obtain the passage-time inequality:
\[
G^{(n), c }_{(\fl{na}, \fl{nb}), (\fl{ns}, \fl{nt})} = G^{(n), T_{\fl{na}, \fl{nb}}c^n}_{{\bf 0}, (\fl{ns}-\fl{na}, \fl{nt}- \fl{nb})} \le G^{(n), T_{a,b}( c^{\eta,K}_{\text{low}})}_{{\bf 0}, (\fl{ns}-\fl{na}, \fl{nt}- \fl{nb})}. 
\]
The superscripts in the passage times above signify which discretisation or shift we are using. 
Dividing by $n$ and letting $n \to \infty$ we have 
\begin{align*}
\varlimsup_{n\to \infty} \frac{1}{n} G^{(n), c }_{(\fl{na}, \fl{nb}), (\fl{ns}, \fl{nt})} &\le \lim_{n\to \infty} \frac{1}{n} G^{(n), T_{a,b}( c^{ \eta,K}_{\text{low}})}_{{\bf 0}, (\fl{ns}-\fl{na}, \fl{nt}- \fl{nb})} \\
&= \Gamma_{T_{a,b}( c^{\eta,K}_{\text{low}})}(s-a, t-b) =  \Gamma_{ ({T_{a,b} c})^{\eta, K}_{\text{low}}}((s-a,t-b)), \quad \text{ by }\, \eqref{Shifdisc} \\
&= \Gamma_{ ({T_{a,b} c})^{\eta}_{\text{low}}}((s-a,t-b)) \quad \text{ by }\, \eqref{eq:diffapp}\\
&= \Gamma_{c^\eta_{\text{low}}}((a,b),(s,t)), \quad \text{ by \eqref{A6}, and then Lemma \ref{lem:gammashift}}.
\end{align*}

Similar arguments can be used to argue that using the discretisation 
\[
 c^{n, \eta,K}_{\text{high}}(i,j) := \max_{(s_n,t_n) \in [-n^{-1}K, n^{-1}K]^2} c^{\eta}_{\text{high}}\Big(\frac{i}{n}+s_n,\frac{j}{n}+t_n\Big) \ge c^{\eta}_{\text{high}}\Big(\frac{i}{n},\frac{j}{n}\Big),
\]
does not alter the limiting shape for $c^{\eta}_{\text{high}}$ and give the lower bound 
\[
\varliminf_{n\to \infty} \frac{1}{n} G^{(n), c^n }_{(\fl{na}, \fl{nb}), (\fl{ns}, \fl{nt})} \ge  \Gamma_{T_{a,b}(c^{\eta}_{\text{high}})}(s-a, t-b) =  \Gamma_{c^\eta_{\text{high}}}((a,b),(s,t)).
\]

Use \eqref{eq:etanew} and let $\eta \to 0$ to obtain 
\be
\lim_{n\to \infty} \frac{1}{n} G^{(n), c^n }_{(\fl{na}, \fl{nb}), (\fl{ns}, \fl{nt})} = \Gamma_c((a,b),(s,t)).  \qedhere
\ee
\end{proof}

We can now prove Proposition \ref{lem:SPEP}. 

\begin{proof}[Proof of Proposition \ref{lem:SPEP}]Fix a $\delta>0$ so that in vector ordering and for $n$ sufficiently large 
\[
(a-\delta, b-\delta)  < n^{-1}(a_n, b_n) < (a+\delta, b + \delta) \text { and } (s-\delta, t-\delta)  < n^{-1}(s_n, t_n) < (s+\delta, t + \delta).
\]
Therefore, by Lemma \ref{lem:LPPshift},
\begin{align*}
\Gamma_{c}((a-\delta, b-\delta), (s+\delta, t+\delta)) 
%&=  \Gamma_{T_{a-\delta, b-\delta}c^{\eta}_{\text{low}}} (c-a+2\delta, d-b+2\delta) \\
%&=\lim_{n\to \infty}\frac{1}{n}G^{(n), T_{a-\delta, b-\delta}c^{\eta}_{\text{low}}}_{(\ce{n(c -a+2\delta)}, \ce{n(d-b+2\delta)})} \\
&= \lim_{n\to \infty} \frac{1}{n}G^{(n)}_{\fl{n(a-\delta, b-\delta)}, \fl{n(s+\delta, t+\delta)}} \ge
\varlimsup_{n\to \infty}\frac{1}{n}G^{(n),}_{(a_n, b_n), (s_n, t_n)}.
\end{align*}
This is true for any $\delta$. Then let $\delta \to 0$, and use Theorem 2.4 in \cite{ciech2021} for the upper bound.
The lower bound follows in a similar manner, approximating from inside; this step requires the continuity of $\Gamma$. 

The last equality in the Proposition follows from Lemma \ref{lem:gammashift}. 
\end{proof}

Finally, we also prove Lemma \ref{lem:betterpaths} that the LLN holds over a smaller set of paths that are a.s.~in continuity regions of $c(x,y)$. In the proof below we assume wihtout loss of generality that the lower left corner of the LPP is $(0,0)$.

\begin{proof}[Proof of Lemma \ref{lem:betterpaths}] Consider any path $\bf x$ from $\mathcal H(x,y)$. The path may or may not follow a non-decreasing discontinuity curve for some time interval $[s_1,s_2] \subseteq [0,1]$. If the path does not follow any discontinuity curve for an interval of time, then it only crosses them, but since it can only find a finite number of curves in up to $(x,y)$ the set of crossing points has measure 0. At all other points the path is in continuity regions of $c$, so it can remain in $\mathcal{H}^{\rm cont}(x,y) $. 

Now assume the path does follow at least discontinuity curve for some interval of time. As such, it cannot be included in the set $\mathcal{H}^{\rm cont}(x,y)$. The path $\bf x$ has at most a countable such parts ${\bf x}_i(s)= {\bf x}(s)\big|_{[s_i,t_i]}$ that follow discontinuity curves $h_i$ for intervals. The weight collected on this interval is 
\[ I({\bf x}_i) = \int_{s_i}^{t_i} \frac{\gamma ({\bf x}'_i(s))}{ c({\bf x}_i(s)))}\,ds.\]
By the lower semicontinuity Assumption \ref{ass:c2},  we know that we can find a set $U$ bounded by the discontinuity curve for which 
\[ c\Big|_{h_i}(x,y)  = \lim_{(w,z) \in U: (w,z) \to (x,y)} c(w,z).\]
Assumption \ref{ass:c0} guarantees that the interior of the set $U$ in the limit above intersects $R_{(0,0),(x,y)}$. This allows for a small perturbation of the path $\bf x$ to still be an admissible path. 

To see this, let $\e > 0$ and consider the set 
\[ C_i = U^{\circ} \cap\{ (w,z): \|c(w,z) - c(x,y)\|_1 < \e 2^{-i} c_{\min}^2(U)/4({\bf x}(t_i) - {\bf x}(s_i)) \} \subseteq \mathcal C. \]
Take any path $\tilde{ \bf x}$ starting at ${\bf x}(s_i)$ and ending at ${\bf x}(t_i)$ but remains in this region. Moreover, assume that except in a set of arbitrary small measure $\delta_i = \e2^{-i}$ near the endpoints, ${\bf x}'(s) = \tilde {\bf x}'(s) = h_i'(s)$  Then the weight collected on $\tilde{\bf x}$ is 
\begin{align*}
\Big|I({\bf x}) - I(\tilde{\bf x}) \Big| &\le  \int_{s_i}^{t_i}  \gamma ({\bf x}'_i(s)) \Big|\left(\frac{1}{ c({\bf x}_i(s))} - \frac{1}{ c(\tilde{\bf x}_i(s))}\right)\Big| + C\delta_i \\
&\le \e2^{-i} +C\delta_i.
\end{align*}
By performing such modifications on each part of ${\bf x}_i$ and switching it to $\tilde{\bf  x}_i$, we construct a path $\tilde{\bf x}$ in $\mathcal{H}^{\rm cont}(x,y)$ with weight $I(\tilde{\bf x}) \ge I({\bf x})- C_0\e$, where $\e$ is arbitrary and $C_0$ is a uniformly bounded constant by the maximum value of $c$ in a compact set, depending only on the endpoints $x,y$. 
Therefore, 
\[
\sup_{{\bf x} \in \mathcal H(x,y)} I({\bf x}) \ge \sup_{{\bf x} \in \mathcal H^{\rm cont }(x,y)} I({\bf x}) \ge \sup_{{\bf x} \in \mathcal H(x,y)} I({\bf x}) -C_0\e.
\]
Letting $\e \to 0$, we obtain the result.
\end{proof}

\section{The proof of the inhomogeneous envelope property} 
\label{app:B}

\begin{proof}[Proof of Lemma \ref{lemma: Envelope property}]

From the construction of the processes, we already have that for any finite time horizon $T$ we can find arbitrarily large $-L$ and $U$ such that $L \ll 0 \ll U$ and no events occurred for $N_{L,-z_{L}^n}$ and  $N_{U,-z_{U}^n}$. 

Since the rates of the Poisson processes are bounded away from zero and infinity, there is only a finite set of Poisson jumps for any finite time horizon $T$ and for any indices in $[L,U]$. In particular, the probability of any two of the processes jumping simultaneously is 0, due to this and their mutual independence.

To prove the statement of the lemma, we will reason by induction on the jump events in $[L, U]  \times[0, T]$. \

It is straightforward to show that the Lemma holds at time $t = 0$ by the construction of our processes.

Now assume the lemma holds up to and including the first $m$ event times. Say the $m$-th event happens at time $t_{m} \in [0,T]$. Until the next event time $t_{m+1} \in [0,T]$  nothing happens, and therefore the relationship still holds. So focus on time $t_{m+1}$. First we want to check that the inequality
\begin{equation} \label{Lemma 3.1 jump inequality}
    z^{n}_{i}(t) \geq z^{n}_{k}(0) - \xi^{n,\bf k}_{i-k}(t)
\end{equation}
continues to hold for all $k$s at time $t_{m+1}$. The only way this can fail to hold is if $z^{n}_{i}$ jumps down at $t_{m+1}$ and $z^{n}_{i}(t_{m+1}^-) = z^{n}_{k}(0) - \xi^{n,\bf k}_{i-k}(t_{m+1}^-)$ for some $k \in \mathbb{Z}$ while $ \xi^{n,\bf k}_{i-k}$ does not jump at $t_{m+1}$. 

From the construction of the processes, we have that $ \xi^{n,\bf k}_{i-k}$ will attempt a jump according to the Poisson process $N_{i, \xi^{n,\bf k}_{i-k}((t_{m+1}^-) - z^{n}_{k}(0) } = N_{i, -z^{n}_{i}(t_{m+1}^-)}$ and in particular, it will attempt to jump together with $z^{n}_{i}$. 
Since $z^{n}_{i}$ is able to jump at time $t_{m+1}$, we must have 
\begin{equation*}
    z^{n}_{i+1}(t_{m+1}^-) = z^{n}_{i}(t_{m+1}^-), \hspace{20.pt} z^{n}_{i-1}(t_{m+1}^-) + 1 = z^{n}_{i}(t_{m+1}^-),
\end{equation*}
from which we obtain
\begin{align*}
    z^{n}_{k}(0) - \xi^{n,\bf k}_{i-k}(t_{m+1}^-) = z^{n}_{i}(t_{m+1}^-) &= z^{n}_{i+1}(t_{m+1}^-) \geq z^{n}_{k}(0) - \xi^{n,\bf k}_{i+1-k}(t_{m+1}^-), \\
    z^{n}_{k}(0) - \xi^{n,\bf k}_{i-k}(t_{m+1}^-) = z^{n}_{i}(t_{m+1}^-) &= z^{n}_{i-1}(t_{m+1}^-) + 1 \geq z^{n}_{k}(0) - \xi^{n,\bf k}_{i-1-k}(t_{m+1}^-) + 1.
\end{align*}
These two equations then yield
\begin{equation*} 
    \xi^{n,\bf k}_{i-k}(t_{m+1}^-) \leq \xi^{n,\bf k}_{i-k+1}(t_{m+1}^-),
    \hspace{20.pt} \xi^{n,\bf k}_{i-k}(t_{m+1}^-) + 1 \leq \xi^{n,\bf k}_{i-k-1}(t_{m+1}^-).
\end{equation*}
These imply that $ \xi^{n,\bf k}_{i-k}$ is also able to jump at $t_{m+1}^-$. Since both $z^{n}_{i}(t)$ and $z^{n}_{k}(0) - \xi^{n,\bf k}_{i-k}(t)$ decrease by 1 at time $t_{m+1}$, \eqref{Lemma 3.1 jump inequality} must hold at the same time. Since we assumed the contrary, we reached a contradiction and therefore the original inequality holds for $t_{m+1}$.

Now it only remains to check that the supremum continues to be attained at time $t_{m+1}$. Suppose that we have $z^{n}_{i}(t_{m+1}^-) = z^{n}_{k}(0) - \xi^{n,\bf k}_{i-k}(t_{m+1}^-)$ for some $k \in \mathbb{Z}$. Then if $z^{n}_{i}$ is able to jump at time $t_{m+1}$, by the first part of this proof,  $ \xi^{n,\bf k}_{i-k}$ is also able to jump at time $t_{m+1}$ and thus the supremum is preserved. 

Next suppose that $z^{n}_{i}$ is unable to jump at time $t_{m+1}$. Then it was blocked by the dynamics and therefore we must have either of
\begin{equation} \label{eq:2cases}
    z^{n}_{i+1}(t_{m+1}^-) = z^{n}_{i}(t_{m+1}^-) + 1 \hspace{15.pt} \text{or} \hspace{15.pt} z^{n}_{i-1}(t_{m+1}^-) = z^{n}_{i}(t_{m+1}^-).
\end{equation}
In the first case of \eqref{eq:2cases}, we have $z^{n}_{i+1}(t_{m+1}^-) = z^{n}_{k_{1}}(0) - \xi^{n, \bf k_{1}}_{i+1-k_{1}}(t_{m+1}^-)$ for some $k_{1} \in \bZ$ by the induction hypothesis. Thus
\begin{align*}
      z^{n}_{i}(t_{m+1}^-) &= z^{n}_{i+1}(t_{m+1}^-) - 1 = z^{n}_{k_{1}}(0) - \xi^{n, \bf k_{1}}_{i+1-k_{1}}(t_{m+1}^-) - 1 \\
    &\leq z^{n}_{k_{1}}(0) - \xi^{n,\bf k_{1}}_{i-k_{1}}(t_{m+1}^-) \leq z^{n}_{i}(t_{m+1}^-).
\end{align*}
In particular,  $z^{n}_{k_{1}}(0) - \xi^{n,{\bf k_{1}}}_{i-k_{1}}(t_{m+1}^-) = z^{n}_{i}(t_{m+1}^-)$ and $\xi^{n,\bf k_{1}}_{i-k_{1}}(t_{m+1}^-) = \xi^{n,\bf k_{1}}_{i+1-k_{1}}(t_{m+1}^-)$ $- 1$. 
This means that $\xi^{n, \bf k_{1}}_{i-k_{1}}$ cannot jump at time $t_{m+1}$ and hence $z^{n}_{k_{1}}(0) - \xi^{n, \bf k_{1}}_{i-k_{1}}(t_{m+1}) = z^{n}_{i}(t_{m+1})$ and the supremum continues to be attained at time $t_{m+1}$. 

In the latter case of \eqref{eq:2cases}, we have $z^{n}_{i-1}(t_{m+1}^-) = z^{n}_{k_{2}}(0) - \xi^{n, \bf k_{2}}_{i-1-k_{2}}(t_{m+1}^-)$ for some $k_{2} \in \bZ$ by the induction hypothesis. Thus
\begin{align*}
      z^{n}_{i}(t_{m+1}^-) &= z^{n}_{i-1}(t_{m+1}^-)  = z^{n}_{k_{2}}(0) - \xi^{n, \bf  k_{2}}_{i-1-k_{2}}(t_{m+1}^-) \\
    &\leq z^{n}_{k_{2}}(0) - \xi^{n, \bf k_{2}}_{i-k_{2}}(t_{m+1}^-) \leq z^{n}_{i}(t_{m+1}^-),
\end{align*}
which implies that $\xi^{n, \bf k_{2}}_{i-1-k_{2}}(t_{m+1}^-)= \xi^{n, \bf k_{2}}_{i-k_{2}}(t_{m+1}^-) $. This similarly implies that $\xi^{n, \bf k_{2}}_{i-k_{2}}$ cannot jump at time $t_{m+1}$. Therefore, we have $z^{n}_{k_{2}}(0) - \xi^{n, \bf k_{2}}_{i-k_{2}}(t_{m+1}) = z^{n}_{i}(t_{m+1})$ and so the supremum once again continues to hold at time $t_{m+1}$. 
\end{proof}

\section{An alternative method for the classical solution using short paths}
\label{app:C}

\begin{proof}[Proof of the supersolution part of Theorem \ref{thm: existence}]
 Let $(x_0,t_0)\in \mathbb{R}\times(0,\infty)$. Assume that $v$ is differentiable at 
$(x_0,t_0)$ and that $\tilde c (\cdot, -v (\cdot, \cdot))$ is continuous at $(x_0,t_0)$. By Lemma \ref{lem:continuity-points-full-measure}, the set for which this does not hold has measure 0.

Set
\[
    c_0:=\tilde c(x_0,-v(x_0,t_0)).
\]

Fix $\xi\in\mathbb{R}$. By the short-time dynamic-programming
inequality \eqref{eq:lb4visco}, applied with $x=x_0$ and $t=t_0$,
we have
\[
    v(x_0,t_0)
    \geq
    v(x_0-\xi h,t_0-h)
    -
    \int_0^h
    \tilde c({\bf w}^{h,\xi}(s))
    \psi\!\left(
        \frac{(w^{h,\xi}_1)'(s)}
             {\tilde c({\bf w}^{h,\xi}(s))}
    \right)\,ds,
\]
where ${\bf w}^{h,\xi}$ is any admissible short path in the restarted
problem with initial profile $v(\cdot,t_0-h)$, chosen so that
\[
    w^{h,\xi}_1(0)=x_0-\xi h,
    \qquad
    w^{h,\xi}_1(h)=x_0.
\]
For instance, take
\[
    w^{h,\xi}_1(s)=x_0-\xi h+\xi s,
\]
and let $w^{h,\xi}_2$ be defined by the Lagrangian relation
\[
    (w^{h,\xi}_2)'(s)
    =
    \tilde c({\bf w}^{h,\xi}(s))
    \psi\!\left(
        \frac{\xi}{\tilde c({\bf w}^{h,\xi}(s))}
    \right),
    \qquad
    w^{h,\xi}_2(0)=-v(x_0-\xi h,t_0-h).
\]
Notice that we do not require
$
    w^{h,\xi}_2(h)\approx-v(x_0,t_0).
$
The path is admissible for the restarted variational problem, and this is
all that \eqref{eq:lb4visco} requires.

Since $v$ is continuous and ${\bf w}^{h,\xi}$ remains in an $O(h)$
neighbourhood of $(x_0,-v(x_0,t_0))$, the continuity of
$\tilde c$ at this point gives
\[
    \frac1h
    \int_0^h
    \tilde c({\bf w}^{h,\xi}(s))
    \psi\!\left(
        \frac{\xi}{\tilde c({\bf w}^{h,\xi}(s))}
    \right)\,ds
    \longrightarrow
    c_0\psi\!\left(\frac{\xi}{c_0}\right).
\]
Therefore
\[
    v(x_0,t_0)-v(x_0-\xi h,t_0-h)
    \geq
    -h c_0\psi\!\left(\frac{\xi}{c_0}\right)+o(h).
\]
Using differentiability of $v$ at $(x_0,t_0)$, a Taylor expansion gives 
\[
    v(x_0-\xi h,t_0-h)
    =
    v(x_0,t_0)
    -
    h\bigl(v_t(x_0,t_0)+\xi v_x(x_0,t_0)\bigr)
    +
    o(h).
\]
Substituting this into the previous inequality gives
\[
    v_t(x_0,t_0)
    +
    \xi v_x(x_0,t_0)
    +
    c_0\psi\!\left(\frac{\xi}{c_0}\right)
    \geq 0.
\]
Since this holds for every $\xi\in\mathbb{R}$, taking the infimum over
$\xi$ and using the duality relation
\[
    c_0 f(p)
    =
    \inf_{\xi\in\mathbb{R}}
    \left\{
        \xi p+
        c_0\psi\!\left(\frac{\xi}{c_0}\right)
    \right\}
\]
we obtain
\[
    v_t(x_0,t_0)+c_0 f(v_x(x_0,t_0))\geq 0.
\]
This proves the supersolution inequality.

\end{proof}

\end{document}